\documentclass[fleqn,11pt]{article}
\usepackage{amssymb,amsmath,amsfonts}
\usepackage[margin=1in]{geometry}
\usepackage{dsfont}
\usepackage{color}
\usepackage{graphicx}
\usepackage{latexsym}
\usepackage{amsmath}
\usepackage{amssymb}
\usepackage{graphics}
\usepackage[dvips]{epsfig}
\usepackage{mathrsfs}
\usepackage{leqno}
\definecolor{refkey}{rgb}{1,0,0.5}
\definecolor{labelkey}{rgb}{0,0.4,1}

\numberwithin{equation}{section}
\newtheorem{thm}{Theorem}[section]

\newtheorem{lem}[thm]{Lemma}
\newtheorem{prop}[thm]{Proposition}
\newtheorem{rmk}[thm]{Remark}
\newtheorem{defn}[thm]{Definition}
\newcommand{\ea}{\epsilon}
\newcommand{\ta}{\theta}
\newcommand{\ka}{\kappa}
\newcommand{\da}{\delta}

\newcommand{\al}{\alpha}
\newcommand{\za}{\zeta}

\newcommand{\pl}{\partial}

\newcommand{\iy}{\infty}

\newcommand{\Da}{\Delta}

\newcommand{\lt}{\left}
\newcommand{\rt}{\right}
\newcommand{\les}{\lesssim}

\newcommand{\be}{\begin{equation}}
\newcommand{\bs}{\begin{split}}
\newcommand{\es}{\end{split}}
\newcommand{\ee}{\end{equation}}
\newcommand{\bee}{\begin{equation*}}
\newcommand{\eee}{\end{equation*}}

\newcommand{\ef}{\eqref}

\allowdisplaybreaks

\begin{document}
\begin{center}
\large{ \bf On  the Free Boundary Problems of 3-D  Compressible  Euler Equations Coupled or Uncoupled With a Nonlinear Poisson Equation}
\end{center}
\centerline{ Tao Luo, Konstantina Trivisa, Huihui Zeng}
\begin{abstract} For the problem of the non-isentropic compressible  Euler Equations coupled with a nonlinear Poisson equation with the electric potential satisfying the Dirichlet boundary condition in three spatial dimensions with a general free boundary not restricting to a graph,   we identify  suitable  stability conditions on  the electric potential
and the pressure  under which we obtain a priori  estimates  on the Sobolev norms of the fluid  and electric variables  and bounds for geometric quantities of free surface. The stability conditions in this case for a general variable entropy are that the  outer normal derivative of the electric potential is positive on the free surface, whereas that  on the pressure is negative. In the isentropic case, the stability condition reduces to a single one, the  outer normal derivative of the difference of the enthalpy and the electric potential is negative on the free surface. For the free boundary problem of the  non-isentropic compressible Euler equations with variable entropy without coupling with the nonlinear Poisson equation, the corresponding higher-order estimates are also obtained under the Taylor sign condition. It is also found that one less derivative is needed to close the energy estimates for the problem for the non-isentropic compressible  Euler Equations coupled with a nonlinear Poisson equation when the electric potential satisfies  the Dirichlet boundary condition under the stability conditions on the electric potential and the pressure, compared with the problem of the non-isentropic compressible Euler equations.
\end{abstract}

\tableofcontents
\section{Introduction}

\subsection{Background and motivations}\label{sec1.1}

Consider the following free boundary problem of non-isentropic compressible Euler-Poisson equations or Euler equations in three dimensions:
\begin{subequations}\label{1}
\begin{align}
&  D_t \rho+\rho {\rm div } v=0   &  {\rm in} \ \  \mathscr{D}_t, \label{1-a}\\
& \rho D_t v+\pl p(\rho,s) =\ka \rho \pl \phi & {\rm in} \ \  \mathscr{D}_t, \label{1-b}\\
&D_t s=0 &   {\rm in} \ \  \mathscr{D}_t, \label{1-c}\\
& p=\bar p  &  {\rm on} \ \   \pl \mathscr{D}_t,  \label{s110} \\
&    v\cdot N =\mathscr{V}(\pl\mathscr{D}_t)
&  {\rm on} \ \  \pl \mathscr{D}_t,\\
&  (\rho, v, s)=(\rho_0, v_0, s_0)  &  {\rm on}  \ \mathscr{D}_0.
\end{align} \end{subequations}
Here the velocity field $v=(v_1, v_2, v_3)$,  density $\rho$,   entropy  $s$, and changing  domain
$\mathscr{D}_t \subset  \mathbb{R}^3$ are the unknowns to be determined;  $\bar p $ is a positive constant, $N $  is the exterior unit normal to the free surface $\pl \mathscr{D}_t$, and $\mathscr{V}(\pl\mathscr{D}_t) $ is the normal velocity of $\pl \mathscr{D}_t$;
$\mathscr{D}_0\subset \mathbb{R}^3$ is a given simply connected bounded domain, and  $(\rho_0, v_0,s_0)$ are given functions.
We use the Einstein summation convention to denote $D_t=\pl_t+v^k\pl_k$ with $v^k=v_k$.
The constant $\kappa$ in \ef{1-b} is either $1$ or $0$. We assume that the pressure function  $p(\rho, s) $
satisfies
\be\label{conP}
 p(\rho, s)\in C^{7-\ka}  \ {\rm and}  \  p_{\rho}(\rho, s)>0 \ {\rm for} \ \rho>0, \  s\in \mathbb{R}.\ee
When $\kappa=0$, system \ef{1} is the celebrated non-isentropic compressible Euler equations. When $\kappa=1$, the potential $\phi$ is determined by
\be\label{1-d}
\Delta \phi+e^{-\phi}=\rho \  \  {\rm in}  \ \mathscr{D}_t,
\ \     \phi=0 \ \  {\rm on} \ \ \pl \mathscr{D}_t.
\ee
The model of \ef{1} and \ef{1-d} describes the motion of a plasma consisting of  cold ions and hot electrons, which arises extensively in astrophysics, plasma physics and semiconductors.
In this context, $\rho$ denotes the density of ions, $e^{-\phi}$ the density of electrons (the electrons follow the classical Maxwell-Boltzmann relation under the massless assumption, cf. \cite{28}),  $\phi$ the electric potential field, $v$ the velocity of the ions,  $p(\rho, s)$ the pressure and $s$ the entropy.
The plasma considered in this paper is unmagnetized, consisting of free electrons and a single species of ions that form a compressible charged  fluid.

The present work devotes to identifying the suitable stability conditions to obtain the a priori estimates on the Sobolev norms of the fluid variables and bounds for geometric quantities of free surface.
System \ef{1} with $\kappa=0$ is just the  non-isentropic compressible Euler equations, for which the local-in-time well-posedness of the free boundary problem was proved  in \cite{35} under  the  Taylor sign condition
\be\label{tls-1}
-\pl_N p(\rho,s)>0 \  \ {\rm on} \ \ \pl \mathscr{D}_t
\ee
for the case  when the free boundary is an unbounded graph and the gravity effect is taken into consideration in the momentum equations, based on the approach of symmetrization of hyperbolic systems and the techniques developed in the study of  weakly stable shock waves and characteristic discontinuities (see \cite{BT, majda,met} for instance).  Here and thereafter, we use the notation $\pl_N=N^i\pl_i$.
The assumption that the free boundary is an unbounded graph in \cite{35} is crucially used to flatten the boundary. For a problem in a general  domain whose boundary is not a graph,  it may be feasible to reduce the problem into the case when the boundary is a graph by using multiple coordinates charts.
However, it is quite technically involved because multiple of free boundary problems have to be solved simultaneously.
Another issue in the local-in-time well-posedness theory in \cite{35} is
the loss of derivatives,
which does not contain full a priori estimates,
since the iteration schemes based on the linearization lose  the  regularity on the moving boundary, and the linearized problems do not preserve the full estimates of the nonlinear problems of which the full symmetry of the problems provided by the physical laws (for example, conservation laws) is used.
In fact, it is proved in \cite{35} that  when the initial data of  the fluid variables $(\rho_0, v_0, s_0)\in H^{m+7}$ and $\pl\mathscr{D}_0\in H^{m+7}$ for $m\ge 6$, there is a local-in-time  solution with
$ (\rho, v, s)(\cdot, t)\in H^{m}$ and $\pl\mathscr{D}_t\in H^{m}$ for $t\in (0, T]$ for some $T>0$. The solution loses
$7$-derivatives.

The first motivation of this paper is to derive the nonlinear higher-order estimates without loss of derivatives for the free boundary problem  of non-isentropic compressible Euler equations, \ef{1} with $\ka=0$, when the free boundary is a general closed surface not restricting to
graphs.
It should be noted that this was achieved in \cite{lindluo} for isentropic compressible  Euler equations which extends the estimates of \cite{CL00} for the incompressible Euler equations.  (See also \cite{lcy} for the a priori estimates on an unbounded domain for isentropic Euler equations with gravity when the free boundary is a graph.)
In the present work, we want to illustrate the role of variations of entropy to the free surface motions. In the isentropic case (that is, the entropy $s$ is a constant),  the pressure $p$ is a sole strictly increasing function of density $\rho$,  so that the enthalpy $h=h(\rho)=\int_1^\rho \lambda^{-1} p'(\lambda)d\lambda$,  pressure $p$ and density $\rho$ are equivalent. One may take either one of them as an independent thermal dynamical variable.  This is an advantage taken in the estimates in \cite{lindluo,lcy}.  Indeed, the enthalpy $h$ is used in \cite{lindluo,lcy} as an independent thermal dynamical variable which satisfies a nice wave equation, since the fact that
$\rho^{-1}\pl p(\rho)=\pl h(\rho)$ is used when  the momentum equation is divided  by $\rho$. However, this does not hold anymore for a variable entropy $s$ for $p=p(\rho, s)$.

When the compressible Euler equations couple with the electric potential field $\phi$, that is, $\kappa=1$ in \ef{1},  and $\phi$ is determined by the Dirichlet problem of a nonlinear Poisson equation \ef{1-d}, the situation becomes more intricate and interesting. New phenomena occur in this case compared with the case of $\kappa=0$ in \ef{1}, in which the usual Taylor sign condition  \ef{tls-1}  is sufficient for our higher-order estimates.
However, for the case of $\kappa=1$  in \ef{1} with $\phi$ satisfying \ef{1-d},
we  need the following stability condition:
\be\label{tls-2}
-\pl_N p(\rho,s)>0 \ \ {\rm and} \  \  \pl_N \phi>0 \  \ {\rm on} \ \ \pl \mathscr{D}_t.
\ee
In the isentropic case for the Euler-Poisson equations (that is, $\kappa=1$ and $s$ being a constant):
\begin{subequations}\label{22.1}
\begin{align}
&  D_t \rho+\rho {\rm div } v=0   &  {\rm in} \ \  \mathscr{D}_t, \label{22.1-a}\\
& \rho D_t v+\pl p(\rho) = \rho \pl \phi & {\rm in} \ \  \mathscr{D}_t, \label{22.1-b}\\
& p=\bar p  &  {\rm on} \ \   \pl \mathscr{D}_t,  \label{22.s110} \\
&    v\cdot N =\mathscr{V}(\pl\mathscr{D}_t)
&  {\rm on} \ \  \pl \mathscr{D}_t,\\
&  (\rho, v)=(\rho_0, v_0)  &  {\rm on}  \ \mathscr{D}_0,
\end{align} \end{subequations}
with potential $\phi$ being determined by \ef{1-d} and pressure $p$ satisfying
\be\label{22.conP}
 p(\rho)\in C^7  \ \ {\rm and}  \ \   p_{\rho}(\rho)>0 \ {\rm for} \ \rho>0,
 \ee
condition \ef{tls-2} for the non-isentropic case can be reduced to a simple one:
 \be\label{23.8}
  \pl_N (\phi-h)>0  \ \ {\rm on} \ \  \pl\mathscr{D}_t, \ee
due to $p=p(\rho)$ and
$h=h(\rho)=\int_1^\rho \lambda^{-1} p'(\lambda)d\lambda$.
We will discuss  these stability conditions later in more detail.
It is also found that, for  the non-isentropic Euler-Poisson equations, under the stability condition \ef{tls-2}, the Dirichlet boundary condition $\phi=0$ on $\pl\mathscr{D}_t$ enables us to use one less
derivative to close the higher-order  estimates than that for the Euler equations.  On the other hand, under the stability condition \ef{23.8} for the isentropic Euler-Poisson equations, we still need the same order of derivatives as that for the Euler equations to close the higher-order estimates.

A critical part of the proof of the  main results of this paper is the establishment of  the higher-order  a priori estimates for the free boundary problem of  compressible fluids equations with variable entropy, to illustrate the role of variations of entropy to the free surface motions and deal with intricate coupling of fluid variables with the electric potential,  for general bounded domains with free boundaries not restricting to
graphs, by adopting the geometric approach developed in \cite{CL00} for the  incompressible Euler equations and used in \cite{lindluo} for isentropic compressible Euler equations.
Our estimates in the present work are based on the following equations: the momentum equations for velocity $v$,  the wave equation for pressure $p$ with the sound speed as the wave speed  (it should be noted the divergence of $v$ also satisfies the same type of wave equation),  the nonlinear elliptic equation for electric potential $\phi$ (it should be noticed that the elliptic equations for gravitational fields or warm plasmas discussed in, for instance \cite{GL, hadjang1, luocy, LN},  are linear for potential),  the transport equations for vorticity ${\rm curl} v$ and entropy $s$. We will deal with  the interaction of the surface wave,
sound wave, entropy wave and electric wave in the paper.

We will construct  higher-order energy functionals and estimate their time derivatives, motivated by \cite{CL00}.  Compared with that in \cite{lindluo,lcy}, the construction of these functionals are quite different even for the fluids variables of density, pressure and velocity, in addition to the entropy and electric potential.  In \cite{lindluo,lcy}, the higher-order energy functionals involve the space-time mixed derivatives of velocity and pressure, while in our construction, only the
space derivatives are involved for the velocity filed, and the mixed  space-time derivatives of pressure contain at most one space derivative.  Moreover, we have terms $\sum_{r=0}^{5-\kappa}\int_{\mathscr{D}_t}  |\pl^r s|^2dx$ and $ \sum_{r=0}^{4-\ka} \int_{\mathscr{D}_t}   |\pl^{r} {\rm div} v|^2 dx $  for the non-isentropic flow, which is not needed in \cite{lindluo,lcy}, since the derivatives of divergence can be neatly controlled for the isentropic flow by the space-time mixed derivatives of the pressure and the higher-order functionals defined  in \cite{lindluo,lcy}. But this is not the case for the non-isentropic flow, in particular, the presence of the electric potential which satisfies a nonlinear elliptic equation makes the estimates of derivatives of ${\rm div} v$  in terms of the other terms in our higher-order energy
functionals extremely difficult. Therefore, we include the term $ \sum_{r=0}^{4-\ka} \int_{\mathscr{D}_t}   |\pl^{r} {\rm div} v|^2 dx $   in our
higher-order energy functional constructions.

The analysis becomes much more involved when the electric filed is taken into consideration. It should be noted that the equation for
the electric potential $\phi$ in \ef{1-d} is nonlinear, which makes the coupling of density, velocity, pressure, entropy, electric potential and the evolving geometry extremely intricate and highly nonlinear. Handling this coupling by various elliptic  estimates is one of the   main concerns of this paper.
Indeed, a big challenge for  free boundary problems of inviscid fluids is that the regularity of the boundary enters to the highest-order estimates, which is
in particularly so for the problem when $\kappa=1$ studied in this paper due to the elaborate coupling mentioned above.

\subsection{Some remarks}
As mentioned in subsection \ref{sec1.1}, we identify the stability condition \ef{tls-2}
for the free boundary problem of non-isentropic  compressible Euler-Poisson equations, that is, \ef{1} with $\kappa=1$ and \ef{1-d}; and condition   \ef{tls-1}
for that  of  non-isentropic compressible Euler equations, that is,  \ef{1} with $\kappa=0$.
Indeed, the Taylor sign condition of the pressure, $-\pl_N p>0$ on $\pl \mathscr{D}_t$,
plays an important role to the stability in the study of  the free surface problems of inviscid fluids, excluding the Rayleigh-Taylor type instability, without which problems may become ill-posed. (See \cite{Ebin} for the problem of incompressible Euler equations.)
For the problem of non-isentropic compressible Euler-Poisson equations, we find that only the Taylor sign condition for the pressure may not be adequate.
In fact, it can be seen from \ef{1-b} that
$$ D_t v\cdot N=-\rho^{-1} {\pl_N p}  +\pl_N \phi \  \ {\rm on} \  \ \pl\mathscr{D}_t. $$
This means that the acceleration of the free surface $\pl\mathscr{D}_t$ is due to two parts, $-\rho^{-1} {\pl_N p}$ and $\pl_N \phi$. Therefore, we propose the stability condition \ef{tls-2} so that $ D_t v\cdot N>0$ on   {\rm $\pl\mathscr{D}_t$.

The stability condition \ef{tls-2} for the general variable entropy case can be replaced by a neat and simple one
in the isentropic case. For  the isentropic problem \ef{22.1}  and  \ef{1-d},
the stability condition we propose in the paper  is \ef{23.8}.  This is motivated by the equivalent form of \ef{22.1-b}:
\be\label{8.8}
D_t v+\pl \mathcal{P}=0, \ \  {\rm where} \ \  \mathcal{P}=h(\rho)-\phi.
\ee
Indeed, $\mathcal{P}$ is a constant and $\pl \mathcal{P} =N \pl_N \mathcal{P} $ on $ \pl\mathscr{D}_t $,  thus on $ \pl\mathscr{D}_t $.
$$
  (-\pl_N \mathcal{P})^{-1} D_t\pl^r \mathcal{P}-N^m \pl^r v_m = (-\pl_N \mathcal{P})^{-1} \pl^rD_t\mathcal{P} +(-\pl_N \mathcal{P})^{-1} \mathcal{R}_r   , $$
where $\mathcal{R}_r=[D_t, \pl^r]\mathcal{P}+(\pl_m \mathcal{P})\pl^r v^m$.
This serves as boundary conditions when one performs the higher-order energy estimates.

Motivated by \ef{tls-1},
one may attempt to try  the following stability condition for the non-isentropic flow,
$$ \rho^{-1} {\pl_N p(\rho, s)}  -\pl_N \phi<0 \ \  {\rm on} \ \  \pl\mathscr{D}_t.  $$
However this does not work since one cannot  write $\rho^{-1} {\pl_N p(\rho, s)}  -\pl_N \phi$ as $\pl_N \mathcal{Q}$
for some scaler function $\mathcal{Q}$ as in the isentropic case. This is a big
difference between the isentropic and non-isentropic flows.
Therefore, we believe the  condition \ef{tls-2} is suitable for the problem of non-isentropic compressible Euler-Poisson equations given by \ef{1} with $\kappa=1$ and \ef{1-d}.

Motivated by \cite{CL00}, we define higher-order energy functionals consisting of
 boundary parts and  interior parts. To define the boundary integrals, we need to
project the equations to the tangent space of the boundary.
\begin{defn}\label{defn1.1} The orthogonal projection $\Pi$ to the tangent space of the boundary of a $(0,r)$ tensor $\al$ is defined to be the projection of each component along the normal:
$$
(\Pi \al)_{i_1\cdots i_r} =  \Pi_{i_1}^{j_1}\cdots \Pi_{i_r}^{j_r} \al_{j_1\cdots j_r}, \ \ {\rm where} \ \  \Pi_{i}^{j}=\da_i^j- {N}_i {N}^j .
$$
The tangential derivative of the boundary is defined by $\overline{\pl}_i=\Pi_i^j\pl_j$, and the second fundamental form of the boundary is defined by $\ta_{ij}=\overline{\pl}_i  { N} _j$.
\end{defn}

The energy functionals we choose in this article contain the $H^{5-\ka}(\mathscr{D}_t)$ norm of the variables
for the non-isentropic compressible problem \ef{1}, which ensures us to estimate the $L^\iy$-bound for the second fundamental form $\theta$ on the boundary.
Recall that the projection formula
$$\theta= (\pl_N p)^{-1}\Pi \pl^2 p  \  \ {\rm on} \ \  \pl \mathscr{D}_t $$
was used to estimate the $L^\iy$-bound for $\theta$
in \cite{CL00}
for the study of the free boundary problem for  incompressible
Euler equations.
The reason why this can work in \cite{CL00}
is because one may obtain the $L^\iy$-bound for $\pl^2 p$ on $\pl \mathscr{D}_t$ independent of
that for $\theta$, which, together with the lower bound for $-\pl_N p$ due to the Taylor sign
condition,  gives the $L^\iy$-bound for $\theta$.
However, we can only obtain, for problem \ef{1} that
\begin{align}\label{23.4}
\|\pl^2 p\|_{L^\iy(\pl\mathscr{D}_t)}
\le C \|\theta\|_{L^\iy(\pl\mathscr{D}_t)}
\|\pl^2 D_t^2 \rho \|_{L^2(\pl\mathscr{D}_t)}
+\textrm{other terms},
 \end{align}
from which it is clear that the projection formula used in \cite{CL00} to give the $L^\iy$-bound for $\theta$ cannot work directly for our problems.
Here and thereafter, ``other terms" means
terms that do not affect those we single out to discuss, and $C$ denotes a certain constant independent of the $L^\iy$-bound for $\theta$.
Indeed, \ef{23.4} follows from Sobolev lemmas and the following estimates:
\begin{subequations}\label{23.6}\begin{align}
&\|\pl^4 p\|_{L^2(\pl\mathscr{D}_t)}
\le C\|\Pi \pl^4 p\|_{L^2(\pl\mathscr{D}_t)}
+C \sum_{0\le r\le 3}\| \pl^r \Da p\|_{L^2(\mathscr{D}_t)},\label{23.6-a}\\
&\| \pl^3 \Da p\|_{L^2(\mathscr{D}_t)}
\le \| \pl^3 D_t^2 \rho \|_{L^2(\mathscr{D}_t)}
+ \textrm{other terms},\label{23.6-b}\\
&  \| \pl^3 D_t^2 \rho \|_{L^2(\mathscr{D}_t)} \le C \|\Pi \pl^3 D_t^2 \rho \|_{L^2(\pl\mathscr{D}_t)}+\textrm{other terms},
\label{23.6-c}\\
&
\|\Pi \pl^3 D_t^2 \rho \|_{L^2(\pl\mathscr{D}_t)} \le 3 \|\theta \|_{L^\iy(\pl\mathscr{D}_t)}\| \pl^2 D_t^2 \rho \|_{L^2(\pl\mathscr{D}_t)}
+\textrm{other terms}. \label{23.6-d}
 \end{align}\end{subequations}
Here \ef{23.6-a} and \ef{23.6-c} follow from elliptic estimates, \ef{23.6-b} from the equation
$\Da p= D_t^2 \rho + \textrm{other terms} $,
and \ef{23.6-d} from the projection formula.

In the study of the non-isentropic compressible Euler equations \ef{1} with $\ka=0$,
we use the same Taylor sign condition \ef{tls-1} but  choose a new energy functional which contains the
$H^{5}(\mathscr{D}_t)$ norm of $p$, instead of
the energy functional chosen in \cite{CL00} containing the $H^{4}(\mathscr{D}_t)$ norm of $p$, since it follows from Sobolev lemmas that
$L^\iy$-bound for $\pl^2 p$ on $\pl \mathscr{D}_t$ can be bounded directly by  the
$H^{5}(\mathscr{D}_t)$ norm of $p$.

For the problem of the non-isentropic compressible Euler-Poisson equations \ef{1} and \ef{1-d} with $\ka=1$ under the stability condition \ef{tls-2}, we use the stability
condition $\pl_N \phi>0$ on $\pl\mathscr{D}_t$ in  \ef{tls-2}, instead of the Taylor sign condition
$-\pl_N p(\rho,s)>0$ on  $\pl \mathscr{D}_t$,
and the projection formula
$$\theta= (\pl_N \phi)^{-1}\Pi \pl^2 \phi$$
to obtain the $L^\iy$-bound for $\theta$.
Indeed,
one can obtain the $L^\iy$-bound for $\pl^2 \phi$ on $\pl \mathscr{D}_t$ independent of
that for $\theta$, which, together with the lower bound for $\pl_N \phi$ due to  \ef{tls-2}, gives the $L^\iy$-bound for $\theta$. It should be noted that the energy functional chosen in this case
contains the $H^{4}(\mathscr{D}_t)$ norm of the variables. However, for the problem of the isentropic compressible Euler-Poisson equations \ef{22.1} and \ef{1-d} under the stability condition \ef{23.8},
we can only obtain, noting $\Da (  h-\phi)=\rho^{-1} D_t^2 \rho +\textrm{other terms}$ and using the same derivation as that of \ef{23.4},  that
$$
\|\pl^2 (h-\phi)\|_{L^\iy(\pl\mathscr{D}_t)}
\le C \|\theta\|_{L^\iy(\pl\mathscr{D}_t)}
\|\pl^2 D_t^2 \rho \|_{L^2(\pl\mathscr{D}_t)}
+\textrm{other terms}.
$$
So, the projection formula used in \cite{CL00} to give the $L^\iy$-bound for $\theta$ cannot work directly for our problem under  condition \ef{23.8}. Other than   the energy functional chosen in \cite{CL00} containing the $H^{4}(\mathscr{D}_t)$ norm of the variables, we choose a new one containing the $H^{5}(\mathscr{D}_t)$ norm, which gives the
$L^\iy$-bound for $\pl^2 (h-\phi)$ on $\pl \mathscr{D}_t$  directly by Sobolev lemmas.

\subsection{Related works}
In many important physical situations, fluid free boundary  problems arise  naturally.  Such problems can be used to model a wide range of phenomena such as water waves, shape of stars, liquid drops,  vortex sheets, etc, depending on the particular hypotheses on fluids.   Much attention has been given to the  case of  homogeneous,
incompressible, and usually inviscid fluids with applications in oceanography
through the water wave problem (cf. \cite{AZ,ABZ,AM,CL00, Coutand,Ebin,L1, L2, Lindblad2, SZ, wu1, wu2,zhang}). To solve those problems, important analytic and geometric techniques are developed.  More recently, the methods developed for
this situation have been brought to bear on models of more complicated fluids.

For incompressible inviscid flows,   the local-in-time well-posedness in Sobolev spaces was first proved in \cite{wu1, wu2} for the irrotational case,  and then in \cite{AZ,ABZ,AM,CL00,Coutand,Ebin,Lindblad2, PN,MW,  SZ,zhang} for prominent progresses including the cases without irrotational assumptions, finite depth water waves, lower regularities, uniform estimates with respect to surface tension, in domains with corner, and etc; the global or almost global-in-time existence for water waves was achieved first  in \cite{GMS,wu4,wu3}, and then in \cite{AD1,BMSW,dengyu,IT1, IT2, IT3, IP1, IP2, IP3, IP4} for recent  developments on this topic including life-span estimates and  the other long-time  well-posedness for related problems; and the singularity formation was proved in \cite{CCFG,CS,wu5}. One may refer to the  survey  \cite{L2} for more references.  The major tools in the study of the above problems , except \cite{ CL00,Coutand,Ebin,Lindblad2,SZ},  rely  on  Fourier analysis,  pseudo-differential operators and analysis on singular integral operators. In this paper,  we adopt a more elementary geometrical approach of dealing with the coupling of the boundary geometry and interior solutions in the spirit of those developed in \cite{CL00, HaoLuo, HaoWang, Lindblad2, lindluo, lcy, LZ, SZ, SZ1, SZ2}, in particular,  \cite{CL00}.
For compressible inviscid  flows, the local-in-time well-posedness of smooth solutions was established for liquids  in  \cite{22,35} (see also \cite{CHS} for  zero surface tension limits); while for gases  with  physical vacuum singularity, the related results  can be found in  \cite{10', 16', LXZ} for the local-in-time theories, and in \cite{hadjang1,hadjang2,LZ1, HZeng, HZeng1} for the global or almost global-in-time ones.

For  the physical vacuum free boundary problem with density tending to zero  in the  rate of $\rho^{\gamma-1}(x, t)\sim {\rm dist} (x, \Gamma(t))$ (where $\gamma>1$ is the adiabatic exponent) near the vacuum boundary $\Gamma(t)$ of the isentropic Euler-Poisson equations describing self-gravitating  fluids modeling gaseous star in 3-D,  the local-in-time well-posedness is obtained in \cite{GL}, the unconditional uniqueness for $\gamma\in(1,2)$ for general 3-D solutions  and the local-in-time well-posedness for the radial solutions without the compatibility conditions of the derivative of the density at the center of symmetry was proved  in \cite{LXZ}, the global expanding solutions are constructed in \cite{hadjang1}, and  the instability of the stationary solution and continued gravitational collapse for $\gamma<4/3$ was was proved in \cite{Jang} and \cite{GHJ}, respectively. It should be noted that those results are for isentropic fluids and the density vanishes on the boundary, while the problem considered in this paper is for the non-isentropic flow and the density is strictly positive on the boundary, so, the stability mechanisms are different. Moreover, the gravitation potential discussed in the above works has an integral representation in terms of density satisfying a linear elliptic equation of the form $\Delta \phi=k\rho$ for a positive constant $k$ with the positive sign of $\Delta \phi$, which models one species of gas. Therefore, for the local-in-time well-posedness,  the gravitational field appears as a lower-order term (see also \cite{luocy} in the study of the local-in-time existence for the case that density is positive  on the boundary).  For the problem studied in this paper,   the electric potential satisfies the Dirichlet problem of a nonlinear elliptic equation, which cannot be treated as the lower-order term. Moreover, $\Delta \phi$ does not have a definite sign which models two species charged particles of ions and electrons. All these  make the coupling of interior solutions and the geometry of the free surfaces much stronger.

\subsection{Organization of the paper and notations}
The rest of the paper are organized as follows. In Section \ref{sec2}, we prepare the necessary materials  including the geometry and regularity of the boundary, Hodge type inequalities,  Sobolev lemmas, interpolation inequalities, elliptic estimates,  estimates for the boundary, commutator estimates  and derivations of higher-order equations for later use. In Section \ref{sec3},
we study the free boundary problem for the non-isentropic
 compressible Euler-Poisson   equations, \ef{1} with $\kappa=1$ and \ef{1-d} under stability condition \ef{tls-2}. The main results in this section are given in Theorem \ref{mainthm}.
The free boundary problem \ef{1} with $\kappa=0$ for the non-isentropic compressible Euler equations is studied  in Section \ref{sec4} with the main results  given in Theorem \ref{thm-2}.
In Section \ref{sec5}, we investigate the free boundary problem \ef{22.1}  and  \ef{1-d} of the isentropic Euler-Poisson equations under the stability condition \ef{23.8} with the main results given in Theorem \ref{thm-3}. The main results in Sections \ref{sec3}-\ref{sec5} give the a priori estimates of  the Sobolev norms of the fluid  and electric variables  and bounds for geometric quantities such as the
second fundamental form and the injectivity radius of the normal exponential map,  of free surface for the corresponding problems. The main general strategy of the proofs of these results is:  identify suitable higher-order functionals,  make appropriate a priori assumptions,  estimate the necessary norms in terms of higher-order energy functionals, perform higher-order energy estimates, prove the main theorems by closing the a priori assumptions.

Throughout the rest of paper, $C$ will denote a universal constant unless stated otherwise, which can change from one inequality to another.
We will employ the
notation
$a\lesssim b$ to denote $a\le Cb$, where $C$ is the universal constant as defined above.
We will use $C(\beta)$ and $C_k(\beta)$ to denote
certain positive constants
depending continuously on quantity $\beta$, which can change from one inequality to another. We will adopt the notation ${\rm Vol} \mathscr{D}_t=\int_{\mathscr{D}_t} dx$.

\section{Preliminaries}\label{sec2}
In this section, we prepare the necessary materials for later use. These materials include the geometry and regularity of the boundary, Hodge type inequalities,  Sobolev lemmas, interpolation inequalities, elliptic estimates,  estimates for the boundary, commutator estimates  and derivations of higher-order equations.

First, we give the following definitions.

\begin{defn}\label{defn1.2}
Let $\iota_0=\iota_0(t)$ be the injectivity radius of the normal exponential map of $\pl \mathscr{D}_t$, that is, the largest number such that the map
$$
\pl \mathscr{D}_t \times (-\iota_0,  \iota_0)   \  \to \ \{x\in \mathbb{R}^n: \  {\rm dist} (x,\pl \mathscr{D}_t)< \iota_0 \} :  \
(\bar x, \iota)   \   \mapsto \ x=\bar x + \iota {N}(\bar x)
$$
is an injection.
\end{defn}

\begin{defn}\label{defn1.3}
Let $d_0$ be a fixed number such that $\iota_0/16 \le d_0 \le \iota_0/2$, and $\eta$ be a smooth cutoff function on $[0, \iy)$ satisfying  $0\le \eta(s)\le 1$, $\eta(s)=1$ when $s\le  d_0/4$, $\eta(s)=0$ when $ s\ge d_0/2$, and  $|\eta'(s)|\le 8/d_0$.
Set
$$ d(t, x)={\rm dist}(x,\pl \mathscr{D}_t),  \ {N}_i (t,x) =- \pl_i d(t,x),
\
   {N}^j(t,x)=\da^{ij}  {N}_i(t,x),     $$
and define
\begin{align*}
&\zeta_{ij}(t,x)=\da_{ij}-\eta^2(d(t,x)) {N}_i(t,x)  {N}_j(t,x) \ \   {\rm in} \ \ \mathscr{D}_t,\\
&\zeta^{ij}(t,x)=\da^{ij}-\eta^2(d(t,x)) {N}^i(t,x)  {N}^j(t,x) \ \   {\rm in} \ \ \mathscr{D}_t.
\end{align*}
In particular, $\zeta$ gives the the induced metric on the tangential space to the boundary:
$$\zeta_{ij}=\da_{ij}- {N}_i {N}_j,  \ \ \zeta^{ij}=\da^{ij}- {N}^i {N}^j  \  {\rm on}  \ \pl \mathscr{D}_t.$$
\end{defn}

\begin{defn}\label{defn2.4} For the multi-indices $I=(i_1, \cdots, i_r)$ and $J=(j_1, \cdots, j_r)$, set
$\da^{IJ}=\da^{i_1j_1} \cdots \da^{i_r j_r}$ and   $\zeta^{IJ}=\zeta^{i_1j_1} \cdots \zeta^{i_r j_r} $.
If $\alpha$ is a $(0,r)$ tensor, define
$|\alpha|^2= \da^{IJ} \alpha_I \alpha_J $.
Then for the projection $(\Pi \alpha)_I=\zeta_I^J \alpha_J$,
$|\Pi\alpha|^2=\zeta^{IJ} \alpha_I \alpha_J$
on $\pl \mathscr{D}_t$.
The $L^p$-norms of a $(0,r)$-tensor $\alpha$ on $\mathscr{D}_t$ and $\pl\mathscr{D}_t$ are denoted, respectively, by $\| \alpha   \| _{L^p}$ and  $ | \alpha    |_{L^p}$:
\begin{align*}
& \| \alpha  \| _{L^p}= \lt(\int_{\mathscr{D}_t} |\alpha|^p dx\rt)^{1/p}   \ {\rm for}  \  1\le p < \iy, \ \ \| \alpha  \| _{L^\iy}={ {\rm ess} \ {\rm sup} }_{\mathscr{D}_t}|\alpha|,\\
&|  \alpha  |_{L^p}= \lt(\int_{\pl \mathscr{D}_t} |\alpha|^p d s\rt)^{1/p}  \ {\rm for}  \ 1\le p  < \iy , \ \  | \alpha  |_{L^\iy}={ {\rm ess} \ {\rm sup} }_{\pl \mathscr{D}_t} |\alpha| .
\end{align*}
\end{defn}
We need another geometric quantity $\iota_1$on $\pl \mathscr{D}_t$ of which it is easier to control the time evolution than $\iota_0$:
\begin{defn}\label{defn2.3}
Let $0<\epsilon_1\le 1/2$ be a fixed number, and let $\iota_1=\iota_1(t)$ depending on $\epsilon_1$ be the largest number such that
$$
| {N}(t,\bar x_1)- {N}(t,\bar x_2)| \le \epsilon_1 \  {\rm whenever} \
|\bar x_1 - \bar x_2|\le \iota_1, \      \bar x_1, \bar x_2 \in \pl\mathscr{D}_t.
$$
\end{defn}
The following Lemma shows that $\iota_1$ is equivalent to $\iota_0$ in conjunction with a bound of the second fundamental form.

\begin{lem}
Suppose that $|\theta|\le K$, and let $\iota_0$ and $\iota_1$ be as in  Definitions \ref{defn1.2} and \ref{defn2.3}. Then
\begin{equation}\label{lemkk1}
\iota_0\ge \min\{\iota_1/2, \ 1/ {K}\} \ \ {\rm and} \ \  \iota_1 \ge \min\{2\iota_0, \ \epsilon_1/ {K}\}.
\end{equation}
\end{lem}
This is Lemma 3.6 in \cite{CL00}.

\begin{lem}   With the notations in Definitions \ref{defn1.2}-\ref{defn2.4}, we have
\begin{align}
\| \pl \za\|_{L^\iy} \le 512\lt(|  \theta |_{L^\iy}+ 1/\iota_0\rt) \ \
{\rm and} \ \  \|D_t \zeta\|_{L^\iy} \le 256 \|\pl v \|_{L^\iy}. \label{CL-3.28}
\end{align}
\end{lem}
The proof of this lemma follows from that of Lemma 3.11 in \cite{CL00}.
Indeed, it holds that
\be\label{dtn}
\pl N=\theta \ \ {\rm and}  \ \
D_t N_i=- N^k (\pl_i -  N_i  \pl_N )v_k=-N^k \overline{\pl}_i v_k \ \  {\rm on} \ \  \pl\mathscr{D}_t
.\ee

The following Lemma gives the Hodge-type inequality.

\begin{lem}\label{lem5.5}
Let $w$ be a $(0,1)$ tensor and define  a scalar ${\rm div} w=\delta^{ij}\pl_i w_j$ and a $(0, 2)$ tensor ${\rm curl} w_{ij}=\pl_i w_j -\pl_j w_i$. If $|\theta|+1/\iota_0\le K$ and $\iota_1\ge 1/K_1$,  then for any nonnegative integer $r$,
\begin{subequations}\label{23.2.19}\begin{align}
&|\pl^{r+1} w|^2 \le C\left(\delta^{ij} \za^{kl} \za^{IJ} (\pl_k \pl^r_I w_i) \pl_l \pl^r_J w_j + |\pl^ r {\rm div} w|^2 + |\pl^r {\rm curl} w|^2 \right), \label{CL-5.16}
\\
&| \pl^r w |_{L^2}^2 \le C \left(\|\pl^{r+1}w\|_{L^2}  + K_1 \|\pl^r w\|_{L^2} \right)\|\pl^r w\|_{L^2}.\label{CL-5.19'}
\end{align}\end{subequations}
\end{lem}
The proof of this lemma follows from that of
Lemmas 5.5 and 5.7 in \cite{CL00}.

Some elliptic estimates are given in the following Lemma.
\begin{lem}\label{prop5.8}
Let $q=q_b$ on $\pl \mathscr{D}_t$ with $q_b$ being a constant. If $|\theta|+1/\iota_0\le K$ and $\iota_1\ge 1/K_1$,  we have  for any $r\ge 2$,
\begin{subequations}\begin{align}
&\|q-q_b\|_{L^2}   \le C({\rm Vol}\mathscr{D}_t)  \|\pl q\|_{L^2}  , \   \|\pl q\|_{L^2}  \le C( {\rm Vol}\mathscr{D}_t) \|\Delta q\|_{L^2}  ,  \label{CL-A.17} \\
& \|\pl^2  q\|_{L^2}  \le C(K, {\rm Vol}\mathscr{D}_t) \|\Delta q\|_{L^2}, \label{LZ-A5a} \\
&\|\pl^r q\|_{L^2}  + | \pl^r q |_{L^2} \le C|  \Pi \pl^r q |_{L^2} + C(K_1, {\rm Vol}\mathscr{D}_t) \sum_{0\le s\le r-1} \|\pl^s \Delta q \|_{L^2} , \label{CL-5.28'}\\
& \|\pl^r q\|_{L^2}  + | \pl^{r-1} q |_{L^2}  \le C |    \Pi \pl^r q |_{L^2} + C( K, {\rm Vol}\mathscr{D}_t) \sum_{0\le s \le r-2} \|\pl^s \Delta q \|_{L^2}. \label{CL-5.29}
\end{align}\end{subequations}
\end{lem}
The proof of \ef{CL-A.17} follows from that of Lemma A.5 in \cite{CL00}, \ef{LZ-A5a} from Corollary A-4 in \cite{LZ}, \ef{CL-5.28'} and \ef{CL-5.29}
from  Proposition 5.8 in \cite{CL00}.

The following Lemma gives some boundary estimates.
\begin{lem}\label{2ndfun} Let $q=q_b$ on $\pl \mathscr{D}_t$ with $q_b$ being a constant and $\iota_1\ge 1/K_1$, then
\begin{subequations} \label{hb1}\begin{align}
& | \Pi \pl^r q  |_{L^2} \le  2 | \pl_N q |_{L^\iy} | \overline{\pl}^{r-2}\theta |_{L^2} + C \sum_{1\le k \le r-1} |  \theta |_{L^\iy}^k |  \pl^{r-k}  q |_{L^2}
 \notag \\
&\quad +C \sum_{1\le k\le r-3} |     \theta   |_{L^\iy} |    \pl_N q   |_{L^\iy} |   \overline{\pl}^{k} \theta  |_{L^2} ,
\ \ r=2,3,4,
\\
& | \Pi \pl^5 q  |_{L^2} \le  2 | \pl_N q |_{L^\iy} | \overline{\pl}^{3}\theta |_{L^2}    + C \sum_{1\le k \le 4} |  \theta |_{L^\iy}^k | \pl^{5-k}  q |_{L^2}  \notag \\
& \quad +C(K_1, | \theta |_{L^\iy})\lt(| \overline{\pl}^{2}\theta |_{L^2}+ | \overline{\pl} \theta |_{L^2}\rt)\sum_{1\le k \le 4}  |  \pl^{k}  q |_{L^2}
.
\end{align}\end{subequations}
If, in addition, $|\pl_N q|\ge \epsilon$ and $|\pl_N q| \ge 2 \epsilon | \pl_N q |_{L^\iy}$ on $\pl \mathscr{D}_t$ for a certain positive constant $\epsilon$, then
\begin{subequations}\label{hb2}\begin{align}
& |  \overline{\pl}^{r-2}\theta  |_{L^2} \le \epsilon^{-2}    | \Pi \pl^r q |_{L^2} + C \epsilon^{-3} \sum_{1\le k \le r-1} | \theta |_{L^\iy}^k |  \pl^{r-k}  q |_{L^2}   \notag \\
&\quad  +C  \epsilon^{-2} \sum_{1\le k\le r-3} |     \theta   |_{L^\iy} |   \pl_N q   |_{L^\iy} |   \overline{\pl}^{k} \theta   |_{L^2} ,     \ \ r=2,3,4, \\
& |  \overline{\pl}^{3}\theta |_{L^2} \le  \epsilon^{-2}    | \Pi \pl^5 q |_{L^2} +  C  \epsilon^{-3}\sum_{1\le k \le 4} |  \theta |_{L^\iy}^k |  \pl^{5-k}  q |_{L^2} \notag \\
&\quad +C(K_1, |  \theta |_{L^\iy})\epsilon^{-3} \lt(| \overline{\pl}^{2}\theta |_{L^2}+ | \overline{\pl} \theta |_{L^2}\rt)\sum_{1\le k\le 4}  |  \pl^{k}  q |_{L^2}.
\end{align}\end{subequations}
\end{lem}
The proof of this lemma follows from that of Proposition 5.9 of \cite{CL00} or Lemmas 2.8 and 2.9 of \cite{LZ}. It should be noticed that on $\pl\mathscr{D}_t$,
\begin{subequations}\label{Pidt2}\begin{align}
&\Pi \pl^2 q - (\pl_N q)\theta=0, \label{22.5.7-1}\\
& |\Pi \pl^3 q - (\pl_N q)\overline{\pl} \theta |
\le 3|\theta||\pl^2 q|
+2 |\theta|^2  |\pl q| . \label{22.5.7-2}
\end{align}\end{subequations}
Indeed, $\epsilon$ appearing on the right-hand side of \eqref{hb2} can be chosen as
\begin{align}
\epsilon= |  (\pl_N q)^{-1}   |_{L^\iy}^{-1} \min\left\{1,  \   2^{-1} |   \pl_N q    |_{L^\iy}^{-1}\right\}. \label{rt1}
\end{align}

\begin{lem}  Let $q=q_b$ on $\pl \mathscr{D}_t$ with $q_b$ being a constant. If $|\theta|+1/\iota_0\le K$ and  $\iota_1\ge 1/K_1$, then
\begin{equation}\label{CL-5.34}
 | \pl  q |_{L^\iy} \le
 \|\pl^2  \Da q \|_{L^2} +  C  ( K,K_1,| \overline{\pl}\theta |_{L^2}, {\rm Vol}\mathscr{D}_t )  ( \|\pl  \Da q \|_{L^2}  +  \|  \Da q \|_{L^2} ).
\end{equation}
\end{lem}
The proof of this lemma follows from that of  Proposition 5.10 of \cite{CL00} or Lemma 2.11 of \cite{LZ}.

The following lemma gives Sobolev inequalities.

\begin{lem}   If $\al$ is a $(0,r)$ tensor, then
\be\label{CL-A.4}
 |  \overline{\pl} \al  |_{L^4}^{2} \le C  |   \al  |_{L^\iy}
|  \overline{\pl}^2 \al |_{L^2} .
\ee
If $\al$ is a $(0,r)$ tensor and $\iota_1\ge 1/K_1$, then
\begin{subequations}\label{22.6.6}\begin{align}
& |  \al  |_{L^{\iy}} \le   C(  K_1 )
( |\pl^2 \al|_{L^2}+|\pl  \al|_{L^2}+| \al|_{L^2} )      , \label{CL-A.8} \\
&\| \al \|_{L^6} \le C(K_1)  (\|\pl \al \|_{L^2}+ \|\al \|_{L^2}),   \label{CL-A.14} \\
&\| \al \|_{L^{\iy}} \le C(K_1) (\|\pl^2 \al \|_{L^2}+   \|\pl \al \|_{L^2}+ \|\al \|_{L^2}). \label{CL-A.15}
\end{align}\end{subequations}
If $\al$ is a $(0,r)$ tensor and $|\theta|+1/\iota_0\le K$, then
\be\label{CL-A.20}
|\al|_{L^{4}} \le C(r, K, {\rm Vol}\mathscr{D}_t)(\|\pl \al \|_{L^2}+\|\al \|_{L^2}).
\ee
\end{lem}
The proof of \ef{CL-A.4} follows from that of Lemma A.1 in \cite{CL00}, \ef{22.6.6}  from Lemmas A.2 and A.4 in \cite{CL00}, \ef{CL-A.20} from Lemma A.7 in \cite{CL00}.

We also need the following    transport formula.
Notice a fact that for any function $f$,
\begin{subequations}\label{3.2-1}\begin{align}
& \frac{d}{dt}  \int_{\mathscr{D}_t } f dx=\int_{\mathscr{D}_t }\lt(D_t f +f {\rm div} v\rt) dx, \label{3.2-1-a} \\
 & \frac{d}{dt}  \int_{\pl \mathscr{D}_t } f ds =\int_{\pl \mathscr{D}_t }\lt(D_t f +f {\rm div} v-fN^i \pl_N v_i \rt)ds
\notag\\&\quad
 = \int_{\pl \mathscr{D}_t }\lt(D_t f -fN^i \pl_N v_i \rt)ds ,\label{3.2-1-b}
\end{align}\end{subequations}
where we have used ${\rm div} v=-\rho^{-1} D_t\rho =-\rho^{-1}\rho_ p D_t p=0$ on $\pl \mathscr{D}_t$, due to $p=\bar p$ on $\pl \mathscr{D}_t$.
This
 implies, with the aid of $D_t \rho + \rho {\rm div} v=0$, that
\be\label{3.2-2}
 \frac{d}{dt}  \int_{\mathscr{D}_t }\rho f dx=\int_{\mathscr{D}_t }\rho D_t f dx.
\ee

The following commutator estimates are useful. Simple calculations show that
for any function $f$, positive integer $r$
and non-negative integer $m$,
 \begin{subequations}\label{7.7}
\begin{align}
& \sum_{1\le k\le r} \lt|\pl^{r -k }[D_t, \pl^k ] f \rt|   \le C(r)
\sum_{1\le k \le r}|\pl^k v| |\pl^{r+1-k} f| ,  \label{3-7-3}\\
&|D_t^m [D_t, \pl^r ] f | \le C(m,r) \sum_{0\le i \le m}
\sum_{1\le k\le r}|D_t^i \pl^k v||D_t^{m-i}\pl^{r+1-k}f|,   \label{22.3.17}\\
& |D_t^m [D_t, \Delta ] f|  \le C(m)  \sum_{0\le i \le m}
(|D_t^{i} \Da v||D_t^{m-i} \pl f|
+ |D_t^{i} \pl v||D_t^{m-i} \pl^2 f|).
\label{3-8-3}
\end{align}\end{subequations}

We also need the following useful inequality. It follows from the H$\ddot{o}$lder inequality   that  for any functions $f$ and $g$,
\be\label{22.5.16}
\|fg\|_{L^2}\le \|f\|_{L^4}\|g\|_{L^4}\le ({\rm Vol}\mathscr{D}_t)^{1/6} \|f\|_{L^6}\|g\|_{L^6}.
 \ee

Finally, we list here some higher-order derived equations from Euler-Poisson equations or Euler equations. It follows from \ef{1-b} that
\be\label{7.5-1}
D_t v_i+ \rho^{-1} \pl_i p= \ka \pl_i \phi,
\ee
which implies, due to $[D_t, \pl_i]=-(\pl_i v^k)\pl_k$, that
\begin{subequations}\label{7.5-2}\begin{align}
&D_t\pl_j v_i +\rho^{-1}  \pl_{ij} p= \rho^{-2}{(\pl_j \rho)\pl_i p} -(\pl_jv^k)\pl_k v_i+\ka \pl_{ij}\phi ,\label{7.5-2-a}
\\
&  D_t {\rm div} v +\rho^{-1}\Delta p= \rho^{-2}{(\pl \rho)\cdot \pl p}- (\pl_jv^k)\pl_k v^j+\ka  \Da \phi, \label{7.5-2-b}\\
 &  D_t ({\rm curl} v)_{ij}=\rho^{-2} \lt((\pl_i \rho)\pl_j p-(\pl_j \rho)\pl_i p
  \rt)+(\pl_j v^k)\pl_k v_i -  (\pl_iv^k) \pl_k v_j.\label{7.5-2-c}
 \end{align}\end{subequations}
In view of \ef{1-a} and \ef{7.5-2-b}, we see that
\be\label{7.5-3}
 D_t^2\rho-\Delta p=\rho \mathcal{H}_0-\kappa \rho \Da \phi,
\ee
where
$
\mathcal{H}_0=({\rm div} v)^2-\rho^{-2}{(\pl \rho)\cdot \pl p}+ (\pl_jv^k)\pl_k v^j
$. This means for $r\ge 1$,
\be\label{22.4.2}
\Da D_t^r p= D_t^{r+2} \rho  - \sum_{0\le k\le r-1}D_t^k[D_t, \Da] D_t^{r-1-k}p- D_t^r(\rho \mathfrak{H}_0)+\kappa  D_t^r(\rho \Da \phi).
\ee
For the potential $\phi$, it follows from \ef{1-d} that
\begin{align}\label{22Mar9-1}
\Da D_t^r \phi - e^{-\phi} D_t^r \phi=  D_t^r \rho-\mathfrak{G}_r, \ \  1\le r\le 4,
 \end{align}
where $\mathfrak{G}_1=[D_t, \Da]\phi$ and for $2\le r\le 4$,
$$\mathfrak{G}_r=\sum_{0\le k \le  r-1}D_t^k [D_t, \Da]D_t^{r-1-k}\phi
-\sum_{0 \le k \le  r-2} D_t^k \lt(e^{-\phi} (D_t \phi) D_t^{r-1-k}\phi \rt).
$$

\section{The non-isentropic  compressible Euler-Poisson   equations}\label{sec3}
In this section, we study the free boundary problem for the non-isentropic
 compressible Euler-Poisson   equations, \ef{1} with $\kappa=1$ and \ef{1-d}, under the stability condition \ef{tls-2}. The main results are given in Theorem \ref{mainthm}.
We define the higher-order  energy functionals as follows:
\begin{align}
& \mathscr{E}_I(t)= \int_{\mathscr{D}_t} \rho |v|^2 dx
+\sum_{1\le r \le 4} \int_{\pl \mathscr{D}_t} |\Pi \pl^r p|^2 (-\pl_N p)^{-1}  ds
\notag\\
&
\quad +\sum_{1\le r \le 4} \int_{\mathscr{D}_t}
\lt(\rho \delta^{mn} \zeta^{IJ} (\pl^r_I v_m)\pl^r_J v_n+ |\pl^{r-1}{\rm curl} v|^2+|\pl^{r-1} {\rm div} v|^2\rt)dx\notag\\
&
\quad + \sum_{0\le r \le 4} \int_{\mathscr{D}_t}  (|\pl^r \rho |^2 + |\pl^r p|^2+|\pl^r s|^2  +|D_t^{r+1} \rho|^2+  \rho_p |\pl D_t^{r}p|^2  )dx,
\notag\\
& \mathscr{E}_{II}(t)=
\sum_{0\le r \le 4} \int_{\mathscr{D}_t}   |\pl^r \phi |^2 dx+\sum_{1\le r \le 4} \int_{\pl \mathscr{D}_t} \rho |\Pi \pl^r \phi|^2 (\pl_N \phi)^{-1}  ds,
\notag\\
& \mathscr{E}_{EP}(t)=\mathscr{E}_I(t)
+\mathscr{E}_{II}(t). \label{23.3}
\end{align}
In order to state the main result , we set
\begin{subequations}\begin{align}
&\underline{\varrho}=\min_{x\in \mathscr{D}_0} \rho_0(x), \  \  \overline{\varrho}=\max_{x\in \mathscr{D}_0}\rho_0(x),   \ \,
\overline{s}=\max_{x\in \mathscr{D}_0} |s_0(x)|,
\label{in1}\\
& \varepsilon_1= \min_{x\in \pl\mathscr{D}_0}(-\pl_N p)( 0,x) , \ \   \varepsilon_2= \min_{x\in \pl\mathscr{D}_0} \pl_N \phi( 0,x) , \label{in2}\\
&K_0= \max_{x\in \pl\mathscr{D}_0} |\theta(0,x)|
+|{\iota_0}^{-1}(0)|,
\label{in3}
\end{align}\end{subequations}
where $p(0,x)=p(\rho_0(x), s_0(x))$, and $\phi(0,x)$ is determined by the Dirichlet problem \ef{1-d}. With these notations, the main results of this section are stated as follows:

\begin{thm}\label{mainthm}
Let $\ka=1$ in \ef{1-b}, and \eqref{conP} hold.
Suppose that
$$0<  \underline{\varrho}, \overline{\varrho}, \overline{s},
\varepsilon_1, \varepsilon_2,   K_0, \mathscr{E}_{EP}(0), {\rm Vol} \mathscr{D}_0<\iy.$$
Then there exists a continuous function
$\mathscr{T}\lt(
\underline{\varrho}^{-1},\overline{\varrho},
\overline{s}, \varepsilon_1^{-1}, \varepsilon_2^{-1}, K_0, \mathscr{E}_{EP}(0),
{\rm Vol} \mathscr{D}_0\rt)>0$
such that  any smooth solution of the free boundary problem \eqref{1}-\eqref{1-d} for $0\le t\le T$ with $T\le \mathscr{T}$
satisfies the  following estimates:  for $0\le t\le T$,
\begin{subequations}\begin{align}
&\mathscr{E}_{EP}(t)\le 2 \mathscr{E}_{EP}(0),  \ \
 2^{-1}{\rm Vol} \mathscr{D}_0\le  {\rm Vol} \mathscr{D}_t\le 2 {\rm Vol} \mathscr{D}_0 , \label{}\\
&2^{-1} \underline{\varrho} \le \min_{x\in \mathscr{D}_t} \rho(t, x), \  \  \max_{x\in \mathscr{D}_t}\rho(t, x)\le 2\overline{\varrho},   \  \ \max_{x\in \mathscr{D}_t}|s(t,x)|\le \overline{s}, \label{29ac}\\
&2^{-1}\varepsilon_1\le \min_{x\in \pl\mathscr{D}_t}(-\pl_N p)( t,x) , \ \
  2^{-1} \varepsilon_2\le \min_{x\in \pl\mathscr{D}_t} \pl_N \phi( t,x),  \\
&\max_{x\in \pl\mathscr{D}_t}|\theta(t,x)|+|\iota_0^{-1}(t)|\le C(  \underline{\varrho}^{-1}, \overline{\varrho}, \varepsilon_2^{-1}, K_0, \mathscr{E}_{EP}(0),
 {\rm Vol} \mathscr{D}_0).
\end{align}\end{subequations}
\end{thm}

\begin{rmk}
It follows from \ef{CL-5.16} and \ef{29ac} that
$\|v(t, \cdot)\|^2_{H^4(\mathscr{D}_t)}\le C \underline{\varrho}^{-1} \mathscr{E}_{I}(t).$
\end{rmk}

To prove the theorem, We make the following a priori assumptions: for $t\in [0,T]$,
\begin{subequations}\label{7.4-2}\begin{align}
&2^{-1} {\rm Vol} \mathscr{D}_0\le  {\rm Vol} \mathscr{D}_t \le 2  {\rm Vol} \mathscr{D}_0, & \label{7.4-2-a}  \\
&2^{-1}\underline{\varrho} \le \rho( t,x)\le 2\overline{\varrho}, \ \  |s(t,x)|\le \overline{s}  & {\rm in}  \ \ \mathscr{D}_t , \label{7.4-2-c} \\
&|\pl (v, p, s)(t,x)|      \le M     \ \ {\rm and} \ \  |\pl \phi(t,x)|\le \bar M  & {\rm in} \ \ \mathscr{D}_t, \label{7.4-2-e}
\\
&|\theta(t,x)|+  \iota_0^{-1}(t)\le K \ \ {\rm and} \ \   \iota_1^{-1}(t)\le K_1 &   {\rm on} \ \  \pl \mathscr{D}_t, \label{7.4-2-g}\\
&-\pl_N p(t,x) \ge \ea_{b} \ \ {\rm and} \ \    |\pl_N D_t p(t,x)| \le  L &   {\rm on} \ \  \pl \mathscr{D}_t,\label{7.4-2-f}\\
& \pl_N \phi(t,x) \ge  \bar\ea_{b}  \ \ {\rm and} \ \    |\pl_N D_t \phi(t,x)| \le \bar L  &   {\rm on} \ \  \pl \mathscr{D}_t,\label{7.4-2-h}
\end{align}\end{subequations}
where $M$, $\bar M$, $K$, $K_1$,  $\ea_{b}$, $\bar\ea_{b}$, $L$ and $\bar L$ are positive constants.
It follows from the maximal principle that for $t\in [0,T]$ and $x\in \mathscr{D}_t$,
\begin{subequations}
\label{7.8.1}\begin{align}
& \phi(t,x )\le \max\lt\{0, \ -\ln \min_{x\in \mathscr{D}_t} \rho(t, x)\rt\},  \label{7.8.1-a}\\
&  \phi(t,x )\ge \min\lt\{0, \ -\ln \max_{x\in \mathscr{D}_t} \rho(t, x)\rt\}. \label{7.8.1-b}
\end{align}
\end{subequations}
Indeed, if there exists $x_0\in \mathscr{D}_t$ such that $\phi( t,x_0)=\max_{x\in \bar{\mathscr{D}}_t }\phi( t,x)$, then
$\Delta \phi(t, x_0)\le 0$ and
$e^{-\phi(t, x_0)} \ge \rho (t,x_0)$, due to \ef{1-d}. That is,
$\phi(t, x_0)\le -\ln \rho (t, x_0)\le -\ln \min_{x\in \mathscr{D}_t} \rho(t, x)$.
By the boundary condition $\phi|_{\pl \mathscr{D}_t}=0$, we have \ef{7.8.1-a}
for $  t \le T$  and $x\in \mathscr{D}_t$.
By the same argument, \ef{7.8.1-b} holds  for $  t \le T$  and $x\in \mathscr{D}_t$.

For the simplicity of the presentation, we may assume without loss of generality that
$$ {\rm Vol} \mathscr{D}_0= {4\pi}/{3},   \  \
\underline{\varrho} =2^{-1},  \  \ \overline{\varrho}=2, \  \ \overline{s} =1,
$$
which implies, due to \ef{7.4-2-a}, \ef{7.4-2-c} and \ef{7.8.1}, that for $t\in [0,T]$,
\begin{subequations}\label{7.4-4}\begin{align}
& {2\pi}/3  \le  {\rm Vol} \mathscr{D}_t \le {8\pi}/3 ,  \label{7.4-4-a}\\
&4^{-1}\le \rho(t,x)\le 4, \ \
|s( t,x )|\le 1,  \  \  |\phi(t,x )|\le \ln 4 \ \  {\rm in} \ \ \mathscr{D}_t.\label{7.4-4-b}
\end{align}\end{subequations}
In view of \ef{1-a} and \ef{1-c}, we see that
\be\label{7.4-5}
D_t \rho=-\rho {\rm div}  v, \ \
 D_t p=p_\rho D_t \rho  =-\rho p_\rho  {\rm div}  v,
 \ee
which, together with  $\pl\rho=\rho_p \pl p + \rho_s \pl s $,
\ef{conP}, \ef{7.4-2-e} and \ef{7.4-4-b}, means that for
$t\in [0,T]$,
\be\label{7.4-6}
 |D_t \rho(t,x)| +|D_t p(t,x)|+|\pl \rho(t,x)|     \les    M \ \ {\rm in} \ \ \mathscr{D}_t.
\ee

\subsection{Regularity estimates}

\begin{prop} Let $\kappa=1$ in \ef{1-b}, then
it holds that
\begin{align}
& \sum_{0\le r \le 4}\big(\|\pl^r v\|_{L^2}^2+
\|D_t \pl^r(\rho, p, s, \phi) \|_{L^2}^2 \big)+  \| \pl^4 {\rm div} v\|_{L^2}^2 +\|D_t v\|_{L^2}^2  \notag\\
& +\sum_{0\le r \le 3} \big(\|D_t \pl^r({\rm curl} v , \ {\rm div}  v )\|_{L^2}^2 +\|D_t \pl D_t^r p \|_{L^2}^2\big)
 +\| D_t^6 \rho - \Da D_t^4 p  \|_{L^2}^2  \notag\\
&
+\|\pl D_t^5 \rho - \rho_p D_t \pl D_t^4 p \|_{L^2}^2
 +\sum_{1\le r\le 4}(|\pl^r(p, \phi)|_{L^2}^2
+|\Pi\pl^r D_t(p, \phi)|_{L^2}^2)
\notag\\
&\le
 C(M,\bar M, K,K_1,\ea_b^{-1},{\bar\ea_b}^{-1},L, \bar L) \sum_{1\le i \le 4}  \mathscr{E}_{EP}^i.
 \notag
\end{align}
\end{prop}

The proof consists of the following four lemmas, Lemmas \ref{22lem-1}-\ref{22lem-4}.
Since some estimates also hold for the non-isentropic compressible Euler equations, $\kappa=0$ in \ef{1-b}, we establish the
estimates depending on $\kappa$ for $\kappa=1$ or $\kappa=0$, which will be used for both cases.

\begin{lem}
\label{22lem-1}
Let $\kappa=0,1$ in \ef{1-b}, then
it holds that
\begin{subequations}\label{lem-1}
\begin{align}
&\|\pl^i v \|^2_{L^2}+  \| D_t^i p \|_{L^2}^2    \le C \mathscr{E}_I ,  \ \ i=1, 2, 3, 4, \label{7.6.1}\\
& \sum_{i=1,2}\|\pl^i ( v,\rho, p, s) \|_{L^\infty}^2
+\sum_{i+j=0,1}\| \pl^i D_t \pl^j (\rho, p) \|_{L^\infty}^2+ \|D_t \pl s \|_{L^\infty}^2
\notag\\
&\ \ +\|\pl^3 ( v,\rho, p, s) \|_{L^6}^2
+ \sum_{i=2,3}\|  D_t^{i}   (\rho,p)  \|_{L^6}^2
+  \sum_{i+j=2} \|  \pl^i D_t \pl^j (\rho,p)  \|_{L^6}^2
\notag\\ & \ \
  +\|\pl  D_t^2  \rho  \|_{L^2}^2
   +\sum_{i=1,2,3} |\pl^i ( v, p) |_{L^2}^2 +|\pl D_t p |^2_{L^2}
\le C(M, K_1)\mathscr{E}_I, \label{7.7.1} \\
& \lt\| (D_t^5 p, \ D_t \pl^3 {\rm curl} v ,\ \pl D_t^3 \rho) \rt\|_{L^2 }^2
+\sum_{i+j=3} \|\pl^i D_t \pl^j (\rho, p)   \|_{L^2}^2
\notag\\ & \ \ +\sum_{1\le i \le 4}\lt\|D_t \pl^i  s\rt\|_{L^2}^2 +|\pl^2 D_t p |^2_{L^2}
\le C(M,K_1)(\mathscr{E}_I^2+\mathscr{E}_I),
\label{7.7.2}\\
& |\overline\pl^{i} \theta|_{L^2}^2 \le C( M, K,K_1,\ea_b^{-1})\mathscr{E}_I, \ \  i=0,1,2. \label{7.6.3}
\end{align}\end{subequations}
\end{lem}

{\em Proof}. \ef{7.6.1} follows from \ef{CL-5.16} and \ef{CL-A.17}. It follows  from
 \ef{CL-5.19'}  and \ef{7.6.1} that
\be\label{7.6.2}
 |\pl^r ( v, p) |_{L^2}^2
   \le C(K_1)\mathscr{E}_I ,  \ \ r=1,2,3 ,
\ee
 which proves  \ef{7.6.3}, using
  \ef{hb2}, \ef{rt1}  and
$|\Pi \pl^4 p|_{L^2}^2 \le |\pl p|_{L^\iy} \mathscr{E}_I$.
In view of \ef{CL-A.14}, \ef{7.6.1} and $|D_t^r \rho| \les \sum_{k=1}^r M^{r-k}|D_t^k p|$ for $r=1, 2,3$, we see that for $r=1, 2,3$,
\be\label{7.2-1}
\|D_t^r p \|^2_{L^6}\le C(K_1)\mathscr{E}_I, \ \
 \|D_t^r \rho \|^2_{L^6}\le C(M, K_1)\mathscr{E}_I,
\ee
which, together with $|D_t^5 p|\les \sum_{j=1}^5 M^{5-j}|D_t^j \rho|+|D_t^3 \rho||D_t^2 \rho|
+M|D_t^2 \rho|^2$ and \ef{22.5.16}, implies
\be\label{7.2-3}
\|D_t^5 p \|^2_{L^2}\le C(M, K_1)(\mathscr{E}_I +\mathscr{E}_I^2).
\ee
It follows from \ef{CL-A.14}, \ef{CL-A.15}, \ef{7.6.1} and \ef{7.4-5}
 that
\begin{align}\label{7.6.4}
\sum_{r=1,2}\|\pl^r ( v,\rho, p, s) \|_{L^\infty}^2
+\|  D_t (\rho, p) \|_{L^\infty}^2 +\|\pl^3 ( v,\rho, p, s) \|_{L^6}^2
\le C( K_1)\mathscr{E}_I .
\end{align}
Notice that
 \begin{align*}
&|\pl^3  D_t {\rm curl} v |\les
M\lt(|\pl^4 (v,\rho,p)|+|\pl^2 \rho ||\pl^2(\rho, p)|\rt)+M^2|\pl^3 (\rho, p)|
\\
&\quad+M^3|\pl^2 (\rho, p)|+ M^4    |\pl p| + |\pl^2 v||\pl^3 v|  + |\pl^2 p||\pl^3 \rho|  + |\pl^2 \rho||\pl^3 p| ,
\end{align*}
which is due to  \ef{7.5-2-c}, then, we use
\ef{7.6.1}, \ef{7.6.4}
and \ef{3-7-3} to get
\be\label{7.2-2}
\lt\| D_t \pl^3 {\rm curl} v \rt\|_{L^2 }^2 \le C( M, K_1) \lt( \mathscr{E}_I  +   \mathscr{E}_I^2 \rt).
\ee
In view of  \ef{1-c}, \ef{3-7-3}, \ef{7.6.1} and \ef{7.6.4}, we see  that for  $1 \le  r \le 4$,
\begin{align}\label{7.7.4}
 \lt\|D_t \pl^r  s\rt\|_{L^2}^2 \le  C( M, K_1) \lt( \mathscr{E}_I  +   \mathscr{E}_I^2 \rt), \  \ \|D_t \pl s\|_{L^\iy}^2 \le  C( M, K_1)\mathscr{E}_I.
\end{align}
It follows from \ef{1-a}  that
\begin{align*}
&|\pl D_t\rho|\les |\pl^2 v|+M|\pl \rho|,\ \
|\pl^2 D_t\rho|\les |\pl^3 v|+M|\pl^2 (v, \rho)| ,\\
&|\pl^3 D_t\rho|\les
|\pl^4 v|+M|\pl^3(v, \rho)|
+|\pl^2 \rho| |\pl^2 v|,
\end{align*}
which gives, using   \ef{7.6.1}, \ef{7.6.4}  and \ef{3-7-3},   that
\begin{subequations}\label{7.7.5}\begin{align}
& \|  (\pl D_t, D_t \pl)  \rho \|_{L^\infty}^2
 +  \sum_{i+j=2} \|  \pl^i D_t \pl^j \rho  \|_{L^6}^2  \le C( M, K_1)\mathscr{E}_I ,
 \label{}\\
 & \sum_{i+j=3} \|  \pl^i D_t \pl^j \rho  \|_{L^2}^2 \le C( M, K_1)\lt(\mathscr{E}_I+\mathscr{E}_I^2\rt) .
 \label{}
\end{align}\end{subequations}
Due to \ef{7.4-5}, we have
$
\pl^r D_t p= - \rho p_\rho \pl^r  {\rm div}  v
+ \mathfrak{F}_r$ for  $r=1,2,3$,
where
\begin{align*}
&|\mathfrak{F}_1|\les M |\pl v|,\ \
|\mathfrak{F}_2|\les M|\pl^2(v,\rho ,s)| +  M^2    |\pl v|,\\
&|\mathfrak{F}_3|\les |\pl^2 v||\pl^2(\rho,s)|+ M|\pl^3(v,\rho,s)| +  M^2 |\pl^2(v,\rho,s)| +M^3  |\pl v| .
\end{align*}
Thus, we use
 \ef{7.6.1}, \ef{7.6.4},   \ef{3-7-3},  \ef{CL-5.19'} and \ef{CL-A.14}
to obtain
\begin{subequations}
\label{7.7.6}
\begin{align}
 & \|(D_t \pl  , \pl D_t)  p  \|_{L^\infty}^2  + \sum_{i+j=2} \|\pl^i D_t \pl^j p  \|_{L^6}^2+|\pl  D_t p |^2_{L^2}  \le C( M, K_1)\mathscr{E}_I,
\label{7.7.6-a}\\
& \sum_{i+j=3} \|\pl^i D_t \pl^j p   \|_{L^2}^2 + |\pl^2 D_t p |^2_{L^2}  \le C( M,   K_1)\lt(\mathscr{E}_I+\mathscr{E}_I^2\rt).
 \end{align}
\end{subequations}
Note that
\begin{align*}
& |\pl D_t^2 \rho|\les |\pl D_t^2 p|+ M|(D_t^2,\pl D_t) p|+M^2|D_t p|,\\
&|\pl D_t^3 \rho|\les \sum_{0\le i\le 1,\ 1\le j\le 3} M^{4-i-j}|\pl^i D_t^j p|+|\pl D_t  p||D_t^2 p|,
\end{align*}
then it follows from  \ef{7.6.1} and \ef{7.7.6-a} that
$$ \|\pl  D_t^2  \rho  \|_{L^2}^2 \le C( M )\mathscr{E}_I ,\ \
  \|\pl D_t^3 \rho \|_{L^2}^2  \le C( M, K_1)\lt(\mathscr{E}_I+\mathscr{E}_I^2\rt) .$$
This, together with \ef{7.6.2}-\ef{7.7.6}, proves \ef{7.7.1} and \ef{7.7.2}.
\hfill $\Box$

\begin{lem}\label{22lem-2} Let $\kappa=0,1$ in \ef{1-b}, then
it holds that
\begin{subequations}\label{lem-2}
\begin{align}
& \sum_{i=1,2}\|\pl^i \phi\|_{L^\iy}^2 +\|\pl^3 \phi\|_{L^6}^2+\sum_{i=1,2,3}|\pl^i \phi|_{L^2}^2
 \le C(K_1) \mathscr{E}_{II}, \label{7.7.8}\\
&\|  D_t^2 (\rho, p) \|_{L^\infty}^2+\sum_{i+j=0,1}\|\pl^i D_t \pl^j v \|_{L^\infty}^2
+ \sum_{i+j=2} \big( \|\pl^i D_t \pl^j  v  \|_{L^6}^2 \notag \\ & \quad+ \|D_t^i \pl D_t^j (  \rho, p) \|_{L^6}^2\big)
\le C(M, K_1)(\mathscr{E}_I + \kappa\mathscr{E}_{II}), \label{7.8.5}\\
& \| \pl D_t^4 \rho\|^2_{L^2 }
+\sum_{i+j=3} \|\pl^i D_t \pl^j v  \|_{L^2}^2
+ \sum_{i+j+k=2}\|\pl^i D_t \pl^j  D_t \pl^k  (\rho, p)  \|_{L^2}^2
\notag \\ & \quad
 +  |\pl D_t^2 p |^2_{L^2} \
  \le C( M, K_1)(\mathscr{E}_I+1)\lt(\mathscr{E}_I+\ka  \mathscr{E}_{II} \rt) ,\label{7.8.4}\\
&  \|\pl D_t^5 \rho - \rho_p D_t \pl D_t^4 p \|_{L^2}^2
 \le  C( M, K_1) (\mathscr{E}_I^2+ \mathscr{E}_I+ 1 ) \lt(\mathscr{E}_I+\ka  \mathscr{E}_{II} \rt).\label{7.8.3}
\end{align}\end{subequations}
\end{lem}

{\em Proof}.   Due to   \ef{CL-5.19'}, \ef{CL-A.14} and \ef{CL-A.15}, one has \ef{7.7.8}.
It follows from \ef{7.5-1} that
$
\pl^r D_t v =-\pl^r (\rho^{-1}\pl p) + \ka \pl^{r+1} \phi$ for $  0\le r\le 3$,
which leads to
\begin{align*}
& |D_t v- \ka \pl \phi| \les |\pl p|, \ \
 |\pl D_t v -\ka \pl^2 \phi|\les |\pl^2 p|+ M|\pl p|
 ,\\
&|\pl^2 D_t v-\ka \pl^3 \phi |\les |\pl^3 p|+M |\pl^2(\rho,p)| +M^2 |\pl p|  ,\\
&|\pl^3 D_t v -\ka \pl^4 \phi|\les \sum_{1\le i \le 4}M^{4-i} |\pl^i (\rho,p )| +|\pl^2 p||\pl^2 \rho|
.
\end{align*}
This, together with  \ef{7.6.1}, \ef{7.7.1},
\ef{7.7.8} and  \ef{3-7-3},  implies
\begin{subequations}\label{7.7.9}\begin{align}
 &\sum_{i+j=0,1}\|\pl^i D_t \pl^j v \|_{L^\infty}^2 +  \sum_{i+j=2} \|\pl^i D_t \pl^j  v  \|_{L^6}^2  \le C( M, K_1)(\mathscr{E}_I
   + \ka \mathscr{E}_{II}) ,
 \label{7.7.9-a}\\
& \sum_{i+j=3} \|\pl^i D_t \pl^j v  \|_{L^2}^2 \le C( M, K_1)\lt(\mathscr{E}_I+\mathscr{E}_{I}^2+ \ka \mathscr{E}_{II}\rt)
 .
 \label{7.7.9-b}
\end{align}\end{subequations}
In view of \ef{7.5-3}, we see that
\begin{align*}
&|\pl D_t^2 \rho| \les |\pl^3(p,\ka\phi)|
+M|\pl^2(v,\rho,p,\ka\phi)|+M^2|\pl(v,p)|,\\
& |\pl D_t^2 p|\les   |\pl D_t^2 \rho|+ M|(D_t^2,\pl D_t) \rho |+M^2|D_t \rho|,
\end{align*}
which gives, using \ef{7.7.1}, \ef{7.7.8}, \ef{7.7.9-a}, \ef{3-7-3} and \ef{22.3.17}, that
\be\label{7.8.2}
\sum_{i+j=2}
  \|D_t^i \pl D_t^j (  \rho, p) \|_{L^6}^2\le C( M, K_1)(\mathscr{E}_I
   + \ka \mathscr{E}_{II}).
\ee
Due to \ef{7.5-3}, $D_t^2 p= p_\rho D_t^2 \rho + p_{\rho \rho }(D_t \rho)^2$, \ef{7.7.1} and \ef{7.7.8}, one obtains
\be\label{7.6.5}
\|  D_t^2 (\rho, p) \|_{L^\infty}^2
\le C(M, K_1)(\mathscr{E}_I + \kappa\mathscr{E}_{II}) ,
\ee
which means, with the aid of \ef{7.6.1} and \ef{7.7.1}, that
\be\label{7.6.8}
\| \pl D_t^4 \rho\|^2_{L^2 }
\le C(M, K_1)(1+\mathscr{E}_I)\lt(\mathscr{E}_I+\ka  \mathscr{E}_{II} \rt)  ,
\ee
since
\begin{align*}
&|\pl D_t^4 \rho|\les
  \sum_{\substack{0\le i \le 1,\  1\le j\le 4}}M^{5-i-j} |\pl^i D_t^j p|
 +  \sum_{\substack{0\le i \le 1, \  1\le j\le 4 \\ 2\le i+j, \ 2\le l \le 4-j  }} M^{5-i-j-l} |\pl^i D_t^j p||D_t^l p|
.
\end{align*}
As a consequence of \ef{7.7.9-a}, \ef{7.8.2} and
 \ef{7.6.5}, we prove \ef{7.8.5}.

It follows from  \ef{7.5-3} that
\begin{align*}
& |\pl^2 D_t^2 \rho|\les
|\pl^4 (p,\ka \phi) | + M|\pl^3(v,\rho, p, \ka\phi)|+ M^2|\pl^2(v,\rho, p)|
\\
& \quad +M^3|\pl p|+|\pl^2 v|^2
 +|\pl^2 \rho||\pl^2 (p, \ka \phi) |,\\
& |\pl^2 D_t^2 p|\les
\sum_{\substack{0\le i, j \le 2 , \ 1\le i+j}}M^{4-i-j} |\pl^i D_t^j(\rho ,s)|+|\pl^2(\rho, s)||D_t^2 \rho|
 +|\pl D_t \rho |^2,
\end{align*}
which, together with  \ef{7.7.1} and \ef{7.7.8}, implies
$$
\|\pl^2 D_t^2 (\rho, p) \|_{L^2}^2
\le C( M, K_1)(1+\mathscr{E}_I)\lt(\mathscr{E}_I+\ka  \mathscr{E}_{II} \rt).
$$
Thus, we use \ef{CL-5.19'},  \ef{3-7-3},
 \ef{22.3.17},   \ef{7.7.1}    and \ef{7.8.5}  to get
$$
\sum_{i+j+k=2}\|\pl^i D_t \pl^j  D_t \pl^k  (\rho, p)  \|_{L^2}^2
 +  |\pl D_t^2 p |^2_{L^2} \
 \le C( M, K_1)(1+\mathscr{E}_I)\lt(\mathscr{E}_I+\ka  \mathscr{E}_{II} \rt) .
$$
Hence, we prove \ef{7.8.4} by virtue of
\ef{7.7.9-b} and
\ef{7.6.8}.

Finally, \ef{7.8.3} follows from  \ef{22.5.16}, \ef{7.6.1}, \ef{7.7.1} and \ef{7.8.5} that
\begin{align*}
&|\pl D_t^5 \rho - \rho_p D_t \pl D_t^4 p|=|\pl D_t^4( \rho_p D_t p) - \rho_p D_t \pl D_t^4 p|\\
&\le |\pl D_t^4( \rho_p D_t p) - \rho_p \pl  D_t^5 p|+\rho_p|\pl  D_t^5 p-D_t \pl D_t^4 p |\\
& \les
\sum_{\substack{0\le i \le 1, \ 1\le j\le 5,\ i+j\le 5}}M^{6-i-j}|\pl^i D_t^j p|
+ |D_t^2 p|^2 |\pl D_t p|\\
&\quad +  \sum_{\substack{0\le i \le 1, \ 1\le j\le 4 \\  2\le i+j , \  2\le l \le 5-j  }} M^{6-i-j-l} |\pl^i D_t^j p||D_t^l p|.
\end{align*}
This finishes the proof of the lemma. \hfill $\Box$

\begin{lem}\label{22lem-3} Let $\kappa=0,1$ in \ef{1-b}, then
it holds that
\begin{subequations}\label{22.lem4}
\begin{align}
&
\sum_{i+j+k=2}\|\pl^i D_t \pl^j D_t \pl^k  v \|_{L^2}^2
 +\sum_{i+j=3}\|D_t^i\pl D_t^j (v,\rho,p) \|_{L^2}^2
 \notag\\
&\quad
+\sum_{i+j=4}\|D_t^i \pl D_t^j ( \rho,p)\|_{L^2}^2 +\sum_{i+j=2}\|D_t^i  \pl D_t^j  v\|_{L^6}^2
  \le C(M,K_1) \mathfrak{C}_1, \label{7.10.1}\\
& \sum_{i+j+k=3}\|D_t^i\pl D_t^j \pl D_t^k (\rho, p)\|_{L^2}^2
+\sum_{i+j=4}\|D_t^i \pl D_t^j  v \|_{L^2}^2
\le C(M,K,K_1)\mathfrak{C}_2,
\label{7.10.2}\\
&\sum_{i+j=4}\|\pl^i D_t \pl^j (\rho,p) \|_{L^2}^2+ |\Pi \pl^3 D_t p|_{L^2}^2 + |\pl^3 D_t p|_{L^2}^2 + |\Pi \pl^4 D_t p|_{L^2}^2
\notag\\
&\quad
+
\|\pl^4 {\rm div} v\|_{L^2}^2
\le C(M,K,K_1,\ea_b^{-1},L) \mathfrak{C}_3,
\label{7.10.3}\\
&\sum_{i+j+k=3}\|\pl^i D_t \pl^j D_t \pl^k ( \rho, p) \|_{L^2}^2
+\| D_t^3 \Da v \|_{L^2}^2
+\| D_t \pl^3 {\rm div} v \|_{L^2}^2
 \notag\\
&\quad
\le C( M, K, K_1, \ea_b^{-1})\mathfrak{C}_4
, \label{lem5-a}\\
&|\pl^4 p|_{L^2}^2
 \le C( M, K, K_1, \ea_b^{-1})\lt(\mathfrak{C}_4
 +\ka \|\pl^3(\rho \Da \phi)\|_{L^2}^2\rt)
, \label{lem5-b}\\
&\| D_t^6 \rho - \Da D_t^4 p  \|_{L^2}^2\le
C( M, K, K_1, \ea_b^{-1})
\lt(\mathfrak{C}_1^2+\mathfrak{C}_2+\mathfrak{C}_4
+\ka\|D_t^4(\rho \Da \phi)\|_{L^2}^2\rt).\label{lem5-c}
\end{align}
\end{subequations}
where
\begin{align*}
&\mathfrak{C}_1=(\mathscr{E}_I+1+\ka  \mathscr{E}_{II} )\lt(\mathscr{E}_I+\ka  \mathscr{E}_{II} +\ka    \|  D_t \pl^2\phi\|_{L^2}^2 \rt)+\ka \|  (D_t^2,  \pl D_t) \pl^2\phi\|_{L^2}^2,\\
&\mathfrak{C}_2=\mathfrak{C}_1+ \ka \|D_t^3 (\rho \Da \phi)\|_{L^2}^2
+ \ka \|  D_t^3 \pl^2\phi\|_{L^2}^2,\ \
 \mathfrak{C}_3=\mathfrak{C}_2+ \ka \|\pl^2 D_t (\rho \Da \phi)\|_{L^2}^2,\\
& \mathfrak{C}_4=
\mathscr{E}_I^3+\mathfrak{C}_1+\ka \| D_t^2  \pl \Da \phi\|_{L^2}^2
+ \ka \| \pl D_t^2(\rho \Da \phi)\|_{L^2}^2.
\end{align*}
\end{lem}

{\em Proof}. It follows from \ef{7.5-2-a} that
\begin{align*}
&|D_t^2 \pl v-\ka D_t\pl^2\phi|\les \sum_{0\le i, j\le 1}M^{2-i-j}|D_t^i \pl^{j+1}p|  + M  |D_t \pl(v,\rho)|,\\
& |\pl D_t^2 \pl v- \ka \pl D_t \pl^2 \phi|\les \sum_{0\le i,j,k \le 1}M^{3-i-j-k} |\pl^i D_t^j \pl^{k+1}p|\\
&\quad+ \sum_{\substack{0\le i, j, k\le 1,\ 2\le i+j+k}}M^{4-i-j-k}|\pl^i D_t^j \pl^k \rho| + |D_t\pl p||\pl^2 \rho| \\
&\quad +|(D_t\pl, \pl D_t)\rho||\pl^2 p| +\sum_{0\le i, j \le 1}|\pl^i D_t^j \pl v||\pl^{1-i}D_t^{1-j}\pl v| ,\\
& |D_t^3 \pl v-\ka D_t^2\pl^2 \phi | \les \sum_{0\le i \le 2, \ 1\le j\le 2}M^{4-i-j}|D_t^i \pl^j p|+|D_t \pl p||D_t \pl \rho|\\
&\quad+|D_t^2 \rho||\pl^2 p| + \sum_{1\le i \le 2, \ 0\le j\le 1}M^{4-i-j}|D_t^i \pl^j \rho|
+\sum_{0\le i\le 2}|D_t^i \pl v||D_t^{2-i} \pl v|
,
\end{align*}
which, together with
\ef{7.7.1}, \ef{7.7.2}, \ef{7.8.5}, \ef{7.8.4}, \ef{3-7-3}, \ef{22.3.17} and \ef{CL-A.14}, gives
\begin{subequations}\label{7.9.1}\begin{align}
& \sum_{i+j=2}\|D_t^i  \pl D_t^j  v\|_{L^2}^2 \les C(M, K_1)(\mathscr{E}_I + \kappa\mathscr{E}_{II})+\ka \|D_t\pl^2\phi\|_{L^2}^2, \label{lem4-a} \\
& \sum_{i+j+k=2}\|\pl^i D_t \pl^j D_t \pl^k  v \|_{L^2}^2
+  \sum_{i+j=2}\|D_t^i  \pl D_t^j  v\|_{L^6}^2
\les  \mathfrak{C} +\ka \|\pl D_t\pl^2\phi\|_{L^2}^2, \label{lem4-b}\\
& \sum_{i+j=3}\|D_t^i\pl D_t^j v \|_{L^2}^2 \les \mathfrak{C}+ \ka \|  (D_t, D_t^2) \pl^2\phi\|_{L^2}^2,
 \label{lem4-c}\\
 & \sum_{i+j=3}\|D_t^i\pl D_t^j (\rho,p)\|_{L^2}^2 \les  \mathfrak{C} +\ka \|D_t\pl^2\phi\|_{L^2}^2,
\end{align}
\end{subequations}
where $\mathfrak{C}=C( M, K_1)(1+\mathscr{E}_I+\ka  \mathscr{E}_{II} )\lt(\mathscr{E}_I+\ka  \mathscr{E}_{II} \rt)$. We then have, using \ef{3-7-3}, \ef{22.3.17}, \ef{7.7.1},  \ef{7.8.5} and \ef{7.8.4}, that
$$
\sum_{i+j=4}\|D_t^i \pl D_t^j ( \rho,p)\|_{L^2}^2
  \les \mathfrak{C}+ \ka (1+\mathscr{E}_I)
  \|  D_t \pl^2\phi\|_{L^2}^2+\ka \|  D_t^2 \pl^2\phi\|_{L^2}^2,
$$
which proves \ef{7.10.1}, by virtue of \ef{7.9.1}.

In view of \ef{22.4.2}, \ef{7.5-2-a},  \ef{3-7-3} and  \ef{22.3.17}, we see that
\begin{align*}
&| \Da D_t^3 p-D_t^5 \rho -\ka  D_t^3(\rho \Da \phi)|  \les \mathcal{J}_1\\
&
  \quad +\sum_{1\le i\le 2} |\pl^i v||\pl^{3-i }D_t^2 p |+\sum_{0\le i, j \le 1}
| D_t^i \pl^{j+1} v||  D_t^{1-i}\pl^{2-j} D_t p|\\
& \quad
+\sum_{\substack{0\le i  \le 2\\ 0\le j\le 1}}
| D_t^i \pl^{j+1} v| |  D_t^{2-i}\pl^{2-j}   p|
+\sum_{\substack{0\le i,j  \le 3\\   i+j\le 3  }} |  D_t^i \rho||  D_t^j \pl v || D_t^{3-i-j} \pl v|,
\\
& |D_t^4 \pl v +\rho^{-1} D_t^3 \pl^2 p-\ka D_t^3 \pl^2 \phi| \les  \mathcal{J}_1+ \sum_{\substack{0\le i  \le 3 }} |  D_t^i \pl v || D_t^{3-i} \pl v|
\\
& \quad +\sum_{0\le i\le 2} M^{3-i} |D_t^i \pl^2 p| +\sum_{\substack{0\le i\le 1 ,\
   2\le j \le 3-i}} M^{3-i-j} |D_t^i \pl^2 p| |D_t^j \rho| ,
\\
 &
 |\pl^2 D_t^3 \rho|
\les \sum_{\substack{0\le i\le 2, \ 0\le j\le 3,\
1\le i+j}} M^{5-i-j}  |\pl^i D_t^j(p,s)|
\\
& \quad + \sum_{\substack{0\le i\le 2, \ 0\le j\le 3, \ 0\le k\le 2-i \\ 0\le l\le 3-j,\
2\le i+j,\ k+l}} M^{5-i-j-k-l}|\pl^i D_t^j(p,s)||\pl^k D_t^l p|   ,
\end{align*}
where
\begin{align*}
& \mathcal{J}_1=\sum_{0\le i\le 3}M^{4-i}\big(|D_t^i\pl p|+\sum_{\max\{0,\ 2-i\}\le j\le 1} M^{1-j}|D_t^i \pl^j \rho|\big)
\\ & \quad
+\sum_{\substack{1\le i\le 2, \ 2\le j\le 3-i \\ \max\{0, \ 2-j\}\le k\le 1}} M^{4-i-j-k} |D_t^i \pl p||D_t^j \pl^k \rho|+ M|D_t^2 \rho||D_t \pl \rho|.
\end{align*}
This proves \ef{7.10.2}, with the aid of \ef{LZ-A5a}, \ef{3-7-3}, \ef{22.3.17}, \ef{7.6.1}  \ef{7.7.1}, \ef{7.8.5}, \ef{7.8.4} and \ef{7.9.1}.

It follows from \ef{22.4.2}  and \ef{3-7-3} that
\begin{align*}
&|\pl^2 \Da D_t p-\ka\pl^2D_t(\rho \Da \phi)|
\les |\pl^2 D_t^3 \rho|+M|\pl^2 \rho||(D_t \pl, \pl D_t)\rho|\\
&+ \sum_{\substack{0\le i, k \le 2, \ 0\le j,l\le 1\\
i+k\le 2, \ j+l\le 1  }} |\pl^i D_t^j \rho||\pl^k D_t^l \pl v ||\pl^{2-i-k} D_t^{1-j-l} \pl v|+\sum_{1\le i\le 4} |\pl^i v||\pl^{5-i }p |   \\
&+\sum_{\substack{0\le i \le 2 ,\ 0\le j\le 1}}M^{4-i-j}\big(|\pl^i D_t^j \pl p|
+\sum_{\max\{0,\ 2-i-j\}\le k\le 1 }M^{1-k}|\pl^i D_t^j \pl^k \rho|\big)\\
&+\sum_{\substack{0\le i \le 2, \ 0\le j,m\le 1, \
1\le i+j  \\
0\le k\le 2-i,\ 0\le l \le 1-j, \   2 \le k+l+m}}M^{4-i-j-k-l-m} |\pl^i D_t^j \pl p||\pl^k D_t^l \pl^m \rho |,
\end{align*}
which gives, with the aid of \ef{7.6.1}, \ef{7.7.1}, \ef{7.7.2}, \ef{7.8.5}, \ef{7.8.4} and \ef{7.10.2},  that
\be\label{22.5.1}
\| \pl^2 \Da D_t p\|_{L^2}^2
\le
C(M, K ,K_1)  \mathfrak{C}_3.
\ee
With
 \ef{lem-1}, \ef{lem-2} and \ef{22.5.1} at hand, we
use \ef{hb1} to obtain the bound for $|\Pi \pl^3 D_t p|_{L^2}$ in \ef{7.10.3}, \ef{CL-5.28'} for $|\pl^3 D_t p|_{L^2}$, \ef{hb1}  for $|\Pi \pl^4 D_t p|_{L^2}$, \ef{CL-5.29} for $\|\pl^4 D_t p\|_{L^2}$, and
 \ef{3-7-3} for $\sum_{i+j=4}\|\pl^i D_t \pl^j p\|_{L^2}$, step by step.
Due to
\ef{1-a} and \ef{1-c}, one has
\begin{align*}
&|\pl^4 ({\rm div} v, D_t \rho)|
= |\pl^4( (-\rho^{-1}, 1) \rho_p D_t p))|
\les M|\pl^2(p,s)|^2 \\
&\quad + \sum_{0\le i\le 4} M^{4-i} |\pl^i D_t p|
+\sum_{0\le i\le 2, \ 2\le j\le 4-i}M^{4-i-j}|\pl^i D_t p||\pl^j(p,s)|,
\end{align*}
which leads to the bounds for
$\|\pl^4 {\rm div} v\|_{L^2}$ and $\sum_{i+j=4}\|\pl^i D_t \pl^j \rho \|_{L^2}$ in  \ef{7.10.3}, using \ef{3-7-3}, \ef{7.6.1},
\ef{7.7.1}, \ef{7.7.2} and the bound just obtained for $\|\pl^4 D_t p \|_{L^2}$.

It follows from   \ef{22.4.2},
 \ef{3-7-3} and \ef{22.3.17}, that
\begin{align*}
&| \pl \Da D_t^2 p- \pl D_t^4 \rho -\ka \pl D_t^2(\rho \Da \phi)|
\les \sum_{1\le i\le 3} |\pl^i v||\pl^{4-i }D_t p |
\\ &
\ +\sum_{\substack{0\le i,k  \le 1, \ 0\le j, l \le 2\\
i+k\le 1,\ j+l\le 2  }} |\pl^i D_t^j \rho||\pl^k D_t^l \pl v ||\pl^{1-i-k} D_t^{2-j-l} \pl v|
+M|D_t^2 \rho||\pl^2 \rho|\\
& \
+\sum_{0\le i, j,  k\le 1}
|\pl^i D_t^j \pl^{k+1} v||\pl^{1-i} D_t^{1-j}\pl^{2-k} p|+M|\pl D_t \rho||D_t \pl \rho|
\\
&\  +\sum_{\substack{0\le i  \le 1, \ 0\le j\le 2  }}M^{4-i-j}\big(|\pl^i D_t^j \pl p|
+\sum_{\max\{0,\ 2-i-j\}\le k\le 1 }M^{1-k}|\pl^i D_t^j \pl^k \rho|\big)\\
&\ +\sum_{\substack{0\le i,m  \le 1, \ 0\le j\le 2, \
1\le i+j  \\
0\le k\le 1-i, \ 0\le l \le 2-j, \ 2\le k+l+m}}M^{4-i-j-k-l-m} |\pl^i D_t^j \pl p||\pl^k D_t^l \pl^m \rho |,
\end{align*}
which, together with \ef{lem-1}, \ef{lem-2}, \ef{7.10.1} and \ef{lem4-a},
means
\begin{subequations}\label{}
\begin{align}
& \|\pl \Da D_t^2 p\|_{L^2}^2
\les C( M,  K_1 ) \mathfrak{C}_1
+ \ka \| \pl D_t^2(\rho \Da \phi)\|_{L^2}^2
.\label{22.5.4}
\end{align}\end{subequations}
In view of \ef{Pidt2}, \ef{CL-A.4} and \ef{CL-A.20}, we see that
\begin{subequations}\label{}
\begin{align}
&|\Pi \pl^2 D_t^2 p|_{L^2}
\le |\theta|_{L^\iy}|\pl D_t^2 p|_{L^2}\le K |\pl D_t^2 p|_{L^2},\label{22.5.15-1}\\
& |\Pi \pl^3 D_t^2 p|_{L^2}
\le |\Pi \pl^3 D_t^2 p- (\pl_N D_t^2 p)\overline{\pl} \theta|_{L^2}
+|(\pl_N D_t^2 p)\overline{\pl} \theta|_{L^2}\notag\\
&\quad \le 3K|\pl^2 D_t^2 p|_{L^2}+2K^2 |\pl D_t^2 p|_{L^2}
+|\pl D_t^2 p|_{L^4}|\overline{\pl} \theta |_{L^4}\notag\\
&\quad \le C(K) \sum_{i=1,2}(|\pl^i D_t^2 p|_{L^2} + \|\pl^i D_t^2 p\|_{L^2} |\theta|_{L^\iy}^{1/2}|\overline{\pl}^2 \theta |_{L^2}^{1/2}).
\label{22.5.15-2}
\end{align}\end{subequations}
We use
\ef{CL-5.28'}, \ef{7.8.4},
\ef{22.5.4} and \ef{22.5.15-1}  to get
$$
| \pl^2 D_t^2 p|_{L^2}^2
\le C( M, K, K_1)
\lt(\mathfrak{C}_1
+ \ka \| \pl D_t^2(\rho \Da \phi)\|_{L^2}^2\rt),
$$
and use \ef{CL-5.29}, \ef{7.6.3}, \ef{7.8.4},
\ef{22.5.4} and \ef{22.5.15-2} to obtain the bound for
$\|\pl^3 D_t^2 p  \|_{L^2}$ in \ef{lem5-a}.
Then,
the bound for  $\|\pl^3 D_t^2 \rho  \|_{L^2}$  follows from     \ef{7.7.1}, \ef{7.7.2}, \ef{7.8.5}, \ef{7.8.4} and that
\begin{align*}
& |\pl^3 D_t^2 \rho|\les
\sum_{\substack{0\le i\le 3, \ 0\le j\le 2, \
1\le i+j}} M^{5-i-j}  |\pl^i D_t^j(p,s)|
\\
& \quad + \sum_{\substack{0\le i\le 3, \ 0\le j\le 2, \ 0\le k\le 3-i \\ 0\le l\le 2-j,\
2\le i+j,  k+l}} M^{5-i-j-k-l}|\pl^i D_t^j(p,s)||\pl^k D_t^l p|   .
\end{align*}
The bound for $\sum_{i+j+k=3}\|\pl^i D_t \pl^j D_t \pl^k (p, \rho) \|_{L^2}$ in \ef{lem5-a} follows from \ef{3-7-3}, \ef{22.3.17}, \ef{lem-1}, \ef{lem-2}  and the bound just obtained for
$\|\pl^3 D_t^2 (\rho,p)  \|_{L^2}$.
It follows from  \ef{7.5-2-a} that
\begin{align}
& D_t \Da v_i+\rho^{-1} \pl_i \Da p=\da^{lj}\lt\{\rho^{-2} (\pl_l \rho) \pl_{ij} p
-(\pl_l v^k)\pl_{kj} v_i\rt.
 \notag \\&\quad  \lt.+ \pl_l \lt( \rho^{-2}{(\pl_j \rho)\pl_i p} -(\pl_jv^k)\pl_k v_i \rt) \rt\}+\ka \pl_i \Da \phi, \label{23.1.15}
\end{align}
which implies
\begin{align*}
&|D_t^3 \Da v +\rho^{-1}D_t^2 \pl \Da p -\ka D_t^2 \pl \Da \phi |
 \les
\sum_{0\le i\le 2}|D_t^i\pl v||D_t^{2-i}\pl^2 v|
\\ &
 \ +\sum_{\substack{0\le i,j\le 2\\ i+j\le 3 }} M^{4-i-j} |D_t^i\pl^{j+1}p|
+\sum_{\substack{0\le i,j\le 2\\ 2\le  i+j }} M^{5-i-j} |D_t^i\pl^{j}\rho|+ M |D_t^2 \rho||\pl^2 \rho|
\\
&\  +\sum_{\substack{0\le i,j\le 2, \ 1\le i+j \le  3\\
0\le k\le 2-i, \ 0 \le l \le 2-j , \ 2 \le k+l}} M^{4-i-j-k-l} |D_t^i\pl^{j+1}p|  |D_t^k\pl^{l}\rho|.
\end{align*}
Then, the bound for $\|D_t^3 \Da v\|_{L^2}$  in \ef{lem5-a} follows from
\ef{lem-1}, \ef{lem-2}, \ef{7.10.1}, \ef{lem4-a} and the bound just obtained for
$\|D_t^2 \pl^3 p \|_{L^2}$.
The bound for $\|D_t \pl^3 {\rm div} v\|_{L^2}$ follows from
\ef{lem-1}, \ef{lem-2}, the bound just obtained for  $\|D_t \pl^3 D_t \rho \|_{L^2}$, and that
\begin{align*}
& |D_t \pl^3 {\rm div} v |=|D_t \pl^3 (\rho^{-1}D_t\rho)|
\les
\sum_{\substack{0\le i, k \le 1, \ 0\le j\le 3\\
1\le i+j+k}}M^{5-i-j-k}|D_t^i\pl^j D_t^k \rho|\\
& \quad
+\sum_{\substack{0\le i, k \le 1, \ 0\le j\le 3, \ 0\le l\le 3- j\\ 0\le m \le 1-k ,\
2\le i+j+k, \ l+m}}M^{5-i-j-k-l-m}|D_t^i\pl^j D_t^k \rho||\pl^l D_t^m \rho|.
\end{align*}
In view of \ef{7.5-3}, we see that
\begin{align*}
&|\pl^3 \Da p-\ka \pl^3(\rho \Da \phi)| \les |\pl^3 D_t^2 \rho|
+\sum_{1\le i \le 4} M^{5-i}|\pl^i (v,\rho, p)|
+|\pl^3 v||\pl^2 v|\\& \quad  +|\pl^3 \rho||\pl^2 p|
+|\pl^3 p||\pl^2 \rho|
+M (|\pl^2 \rho||\pl^2(v,\rho, p)|+|\pl^2 v|^2)
,
\end{align*}
which proves \ef{lem5-b}, using \ef{CL-5.28'},
\ef{lem-1} and the bound just obtained for  $\|\pl^3 D_t^2 \rho \|_{L^2}$.

It follows from \ef{22.4.2}  and  \ef{7.7} that
\begin{align*}
&|  \Da D_t^4 p-D_t^6 \rho - \ka D_t^4(\rho \Da \phi) |
\\
& \les \sum_{\substack{0\le i,j  \le 4, \   i+j\le 4  }} |  D_t^i \rho||  D_t^j \pl v || D_t^{4-i-j} \pl v| +\sum_{1\le i\le 2} |\pl^i v||\pl^{3-i }D_t^3 p |\\
& \quad
+\sum_{0\le i, j \le 1}
| D_t^i \pl^{j+1} v||  D_t^{1-i}\pl^{2-j} D_t^2 p|+M|D_t^3\Da v|+|D_t^3 \pl v||\pl^2 p|
\\ & \quad
+\sum_{0\le i  \le 2, \ 0\le j\le 1}
| D_t^i \pl^{j+1} v| ( |  D_t^{2-i}\pl^{2-j} D_t p|+|D_t^{3-i}\pl^{2-j}  p|)+\mathcal{J}_2,
\end{align*}
where
\begin{align}
& \mathcal{J}_2=\sum_{0\le i\le 4}M^{5-i}\big(|D_t^i\pl p|+\sum_{\max\{0,\ 2-i\}\le j\le 1} M^{1-j}|D_t^i \pl^j \rho|\big)
\notag\\ & \ \
+\sum_{\substack{1\le i\le 3, \ 1\le j\le 4-i \\ \max\{0, \ 2-j\}\le k\le 1}} M^{5-i-j-k} |D_t^i \pl p||D_t^j \pl^k \rho|+|D_t \pl p||D_t^2 \rho||D_t\pl \rho|
\notag\\ & \ \
+\sum_{\substack{1\le i\le 3, \ 0\le j\le 1, \ 1\le k\le 4-i \\ 0\le l\le 1-j,  \  2\le i+j, k+l }} M^{6-i-j-k-l} |D_t^i \pl^j \rho||D_t^k \pl^l \rho|. \label{23-1-21}
\end{align}
This proves \ef{lem5-c}, using \ef{22.5.16}, \ef{lem-1},
\ef{lem-2},   \ef{7.10.1}, \ef{7.10.2} and
\ef{lem5-a}. \hfill$\Box$

\begin{lem} \label{22lem-4}  Let $\ka=1$ in \ef{1-b}, then it holds that
\begin{subequations}\label{22.lem3}\begin{align}
&|\pl^4 \phi|_{L^2}^2+|\pl^2 \phi|_{L^\iy}^2 +\|\pl^3(\rho \Da \phi)\|_{L^2}^2\le C(M, \bar M, K_1)   \mathscr{E}_{EP}, \label{22.6.26}\\
&\|  D_t  \phi\|_{L^\iy}^2 +
\sum_{i+j=1,2}\|\pl^i D_t \pl^j  \phi\|_{L^2}^2  + |\pl D_t  \phi|_{L^2}^2   \le C(M, \bar M, K, K_1)\mathscr{E}_{EP}, \label{22Mar9-4}\\
&\sum_{i+j=3,4} \| \pl^i D_t  \pl^j \phi \|_{L^2}^2
  + \sum_{i+j=1,2} \|  \pl^i D_t \pl^j \phi  \|_{L^\iy}^2 + \sum_{i=2,3}|\pl^i D_t  \phi|_{L^2}^2
  \notag\\
 &  \quad  + \sum_{i=2,3,4} |\Pi \pl^i D_t  \phi|_{L^2}^2
   \le C( M, \bar M, K, K_1, \ea_b^{-1}, \bar L)(\mathscr{E}_{EP}^2+\mathscr{E}_{EP}), \label{22Mar9-5}
\\
& \sum_{i+j+k=2}\|\pl^i D_t \pl^j D_t \pl^k  \phi\|_{L^2}^2 +\sum_{i+j=2} \| D_t^i \pl D_t^j  \phi \|_{L^6}^2  +\|  D_t^2  \phi\|_{L^\iy}^2
+\|D_t^2 \pl    \Da \phi\|_{L^2}^2 \notag\\
 &\quad +\|(\pl^2 D_t, \pl D_t^2)(\rho \Da \phi)\|_{L^2}^2
\le C( M,\bar M, K, K_1)(\mathscr{E}_{EP}^2+\mathscr{E}_{EP}), \label{8.1-1}\\
  &\|D_t^3 \phi\|_{L^\iy}^2
+\sum_{i+j=3}\|D_t^i\pl D_t^j \phi\|_{L^2}^2
+\sum_{i+j+k=3}\|D_t^i\pl D_t^j \pl D_t^k \phi\|_{L^2}^2 \notag\\
&\quad + \|   D_t^3 (\rho \Da \phi)\|_{L^2}^2 \le C( M, \bar M, K, K_1, \ea_b^{-1}, \bar L) \sum_{1\le i\le 3}
\mathscr{E}_{EP}^i,\label{8.1-5}\\
&\|  D_t^4\phi \|_{L^\iy}^2+\sum_{i+j=4}\|D_t^i\pl D_t^j \phi\|_{L^2}^2
+ \|\pl^2 D_t^4 \phi\|_{L^2}^2
+\|D_t^4(\rho \Da \phi)\|_{L^2}^2\notag\\
& \quad
\le C( M, \bar M, K, K_1, \ea_b^{-1}, \bar L)(\mathscr{E}_{EP}^4+\mathscr{E}_{EP}^3
+\mathscr{E}_{EP}^2 +\mathscr{E}_{EP}
). \label{8.1-2}
\end{align}\end{subequations}
\end{lem}

{\em Proof}. We assume $\ka=1$ in the proof. The bound for $|\pl^4 \phi|_{L^2}$
in \ef{22.6.26} follows from
$|\Pi \pl^4 \phi |_{L^2}^2 \le |\rho^{-1}\pl \phi|_{L^\iy} \mathscr{E}_{II}$, \ef{CL-5.28'}
and
$$ |\pl^3 \Da \phi|\le|\pl^3 e^{-\phi}|+|\pl^3 \rho| \les |\pl^3\rho |+\sum_{1\le i\le 3} \bar M^{3-i}|\pl^i \phi|;$$
that for $|\pl^2 \phi|_{L^\iy}$
from \ef{CL-A.8}, \ef{7.7.8} and the bound just obtained for $|\pl^4 \phi|_{L^2}$;
that for $\|\pl^3(\rho \Da \phi)\|_{L^2}$ from
$$|\pl^3(\rho \Da \phi)|
\le |\pl^3(\rho e^{-\phi})|+|\pl^3\rho^2|\les \sum_{0\le i\le 2}(M^i+\bar M^i)|\pl^{3-i}(\rho, \phi)|. $$

We multiply \ef{22Mar9-1} by $D_t^r \phi$ to get
$$|\pl D_t^r \phi|^2+ e^{-\phi}|D_t^r \phi|^2
= {\rm div} \lt(D_t^r \phi\pl D_t^r \phi \rt)
-(D_t^r \rho-\mathfrak{G}_r)D_t^r \phi,$$
which implies, due to   $D_t^r \phi=0$ on
$\pl \mathscr{D}_t$, that for $1\le r\le 4$,
\be\label{22Mar9-2}
\|\pl D_t^r \phi\|_{L^2}^2+ \|D_t^r \phi\|_{L^2}^2
\les
\|D_t^r \rho\|_{L^2}^2 + \|\mathfrak{G}_r\|_{L^2}^2 .
\ee

In view of \ef{3-7-3}, we see that
\begin{align*}
&|\mathfrak{G}_1|\les M |\pl^2 \phi |+\bar M |\pl^2  v|,\ \ |\pl \mathfrak{G}_1|\les M |\pl^3 \phi |+\bar M |\pl^3 v|+|\pl^2 v||\pl^2 \phi|,\\
&|\pl^2 \mathfrak{G}_1|\les
M |\pl^4 \phi |+\bar M |\pl^4  v|+|\pl^3 v||\pl^2 \phi|+|\pl^2 v||\pl^3 \phi|,
\end{align*}
which, together with \ef{7.6.1}, \ef{7.7.1} and \ef{7.7.8}, means
\be\label{22-3-9-1}
\|\mathfrak{G}_1 \|_{L^2}^2 \le C(M,\bar M)\mathscr{E}_{EP}, \ \
\|\pl \mathfrak{G}_1 \|_{L^2}^2 +\|\pl^2 \mathfrak{G}_1 \|_{L^2}^2\le C(M,\bar M,K_1)(\mathscr{E}_{EP} +\mathscr{E}_{EP}^2).
\ee
Thus, we use \ef{22Mar9-2}, \ef{3-7-3} and  \ef{22Mar9-1} to obtain
$$
\|D_t  \phi\|_{L^2}^2 + \|(\pl D_t  , D_t \pl) \phi ) \|_{L^2}^2+\|\Da D_t  \phi\|_{L^2}^2  \le C(M,\bar M)\mathscr{E}_{EP}.
$$
Hence, \ef{22Mar9-4} follows from \ef{LZ-A5a}, \ef{3-7-3}, \ef{7.6.1}, \ef{CL-A.15}  and \ef{CL-5.19'}.
Note that
\begin{align*}
&|\pl( e^{-\phi} D_t  \phi)|\les |\pl D_t  \phi| + \bar M |D_t    \phi|,   \\
& |\pl^2 ( e^{-\phi} D_t   \phi)|\les
|\pl^2 D_t  \phi| +  \bar M |\pl D_t  \phi|+(|\pl^2 \phi| + \bar M^2) |D_t    \phi|,
\end{align*}
then it yields from \ef{22Mar9-1}, \ef{22Mar9-4},
\ef{7.7.1} and \ef{22-3-9-1} that
\be\label{22-3-9-3}
\|\pl \Da D_t  \phi\|_{L^2}^2+ \|\pl^2 \Da D_t  \phi\|_{L^2}^2
\le  C( M, \bar M, K, K_1)\lt(\mathscr{E}_{EP}+\mathscr{E}_{EP}^2\rt).
\ee
With \ef{7.6.3}, \ef{22Mar9-4} and \ef{22-3-9-3} at hand, we prove \ef{22Mar9-5} as follows.  We use \ef{hb1} and \ef{22Mar9-4} to obtain the bound for $ |\Pi \pl^2 D_t  \phi|_{L^2}$,  \ef{CL-5.28'}  for $ | \pl^2 D_t  \phi|_{L^2}$,  \ef{hb1} for $ |\Pi \pl^3 D_t  \phi|_{L^2}$,
\ef{CL-5.28'}  for $ \| \pl^3 D_t  \phi\|_{L^2}$ and $| \pl^3 D_t  \phi|_{L^2}$, \ef{hb1} for $ |\Pi \pl^4 D_t  \phi|_{L^2}$,
\ef{CL-5.29}  for $ \| \pl^4 D_t  \phi\|_{L^2}$,
\ef{CL-A.15} for  $ \| \pl D_t  \phi\|_{L^\iy}$ and $ \| \pl^2 D_t  \phi\|_{L^\iy}$,
and \ef{3-7-3}, \ef{7.6.1}, \ef{7.7.1} and \ef{7.7.8}
for the remaining terms,
step by step.

The bound for $\|\pl^2 D_t(\rho \Da \phi)\|_{L^2}$ in \ef{8.1-1} follows  from
\ef{7.7.1}, \ef{22Mar9-4}, and that
\begin{align*}
& |\pl^2 D_t(\rho \Da \phi)| \les |\pl^2 D_t(\rho, \phi)|+(M+\bar M)|\pl D_t(\rho, \phi)|+(M \\
&\quad +|D_t \phi|)|\pl^2(\rho, \phi)|
+(M^2+M\bar M +{\bar M}^2)|D_t(\rho, \phi)|,
\end{align*}
which is due to \ef{1-d}.
Next, we prove
\begin{align}
& \sum_{i+j=2} \| D_t^i \pl D_t^j  \phi \|_{L^2}^2 + \sum_{i+j+k=2}\|\pl^i D_t \pl^j D_t \pl^k  \phi\|_{L^2}^2 +\sum_{i+j=2} \| D_t^i \pl D_t^j  \phi \|_{L^6}^2 \notag\\
 &\quad +\|  D_t^2  \phi\|_{L^\iy}^2
\le C( M,\bar M, K, K_1)(\mathscr{E}_{EP}+\mathscr{E}_{EP}^2). \label{22Mar17-2}
\end{align}
It follows from \ef{3-7-3} and \ef{22.3.17}  that
$$
|\mathfrak{G}_2|\les
\sum_{0\le r, j\le 1, \ 0\le i\le r} |D_t^i \pl^{j+1} v|
|D_t^{r-i}\pl^{2-j}D_t^{1-r}\phi| + |D_t\phi|^2  ,
$$
which implies, due to \ef{7.7.1}, \ef{7.8.5},  \ef{22Mar9-4} and \ef{22Mar9-2}, that
$$
\|\mathfrak{G}_2\|_{L^2}^2 +  \| \pl D_t^2  \phi\|_{L^2}^2
+  \|  D_t^2  \phi\|_{L^2}^2 \le C(M,\bar M, K,K_1)(\mathscr{E}_{EP} +\mathscr{E}_{EP}^2)
.$$
Then, we use \ef{LZ-A5a} and \ef{22Mar9-1} to obtain
the bound for $\|\pl^2 D_t^2\phi\|_{L^2}$ in \ef{22Mar17-2}; \ef{CL-A.15} and \ef{CL-A.14}
for $\|D_t^2   \phi\|_{L^\iy}$ and $\|\pl D_t^2 \phi\|_{L^6}$;
and \ef{3-7-3},  \ef{22.3.17}, \ef{7.7.1}, \ef{7.8.5} and  \ef{22Mar9-4}  for the remaining terms, step by step. Due to \ef{1-d}, one has
\begin{align*}
& |D_t^2 \pl    \Da \phi |\les |  D_t^2\pl(\rho, \phi)|+
\bar M(|D_t^2 \phi|+|D_t \phi|^2) + |D_t \phi| |D_t \pl \phi|
,\\
& |\pl D_t^2 (\rho \Da \phi)|\les |\pl  D_t^2(\rho, \phi)|+(M +|D_t \phi|)|\pl D_t(\rho, \phi)|+(M \\
&\quad +\bar M)(|D_t^2(\rho, \phi)|+|D_t \phi|^2
+M|D_t(\rho, \phi)|),
\end{align*}
which give the bounds for $\|D_t^2 \pl    \Da \phi\|_{L^2}$ and $\| \pl D_t^2 (\rho \Da \phi)\|_{L^2}$ in \ef{8.1-1}, by virtue of \ef{7.7.1}, \ef{7.8.5},
\ef{22Mar9-4} and \ef{22Mar17-2}.

In a  similar way to deriving   \ef{22Mar9-4}, we can prove
\begin{align}
  &\|D_t^3 \phi\|_{L^\iy}^2
+\sum_{i+j=3}\|D_t^i\pl D_t^j \phi\|_{L^2}^2
+\sum_{i+j+k=3}\|D_t^i\pl D_t^j \pl D_t^k \phi\|_{L^2}^2 \notag\\
&\quad \le C( M, \bar M, K, K_1, \ea_b^{-1}, \bar L)(\mathscr{E}_{EP}+\mathscr{E}_{EP}^2
+\mathscr{E}_{EP}^3),\label{7.29-1}\\
&\|  D_t^4\phi \|_{L^\iy}^2+\sum_{i+j=4}\|D_t^i\pl D_t^j \phi\|_{L^2}^2
+ \|\pl^2 D_t^4 \phi\|_{L^2}^2\notag\\
&\quad
\le
C( M, \bar M, K, K_1, \ea_b^{-1}, \bar L)
\sum_{1\le i\le 4}
\mathscr{E}_{EP}^i ,\notag
\end{align}
step by step,
noting that
\begin{align*}
& |\mathfrak{G}_3|
\les
  \sum_{0\le r\le 2, \ 0\le i\le r,\ 0\le j\le 1}|D_t^i \pl^{j+1} v ||D_t^{r-i}\pl^{2-j}D_t^{2-r} \phi|
\\
& \ \ +|D_t^2 \phi| |D_t\phi|+|D_t\phi|^3,\\
&|\mathfrak{G}_4|\les \sum_{0\le r\le 2, \ 0\le i\le r,\ 0\le j\le 1}|D_t^i \pl^{j+1} v ||D_t^{r-i}\pl^{2-j}D_t^{3-r} \phi|
\\
&\ \ +\sum_{0\le i\le 3}(|D_t^i \Da v||D_t^{3-i}\pl \phi| + |D_t^i \pl v||D_t^{3-i }\pl^2 \phi|)\\
&\ \
+|D_t^3 \phi||D_t \phi|+|D_t^2 \phi|^2+ |D_t^2 \phi||D_t\phi|^2
+|D_t \phi|^4,
\end{align*}
which are due to
\ef{7.7}.
The bound for $\| D_t^3 (\rho \Da \phi)\|_{L^2}$ in \ef{8.1-5} follows from  \ef{22Mar9-4},
\ef{22Mar17-2}, \ef{7.29-1} and that
$$ | D_t^3 (\rho \Da \phi)|\les |  D_t^3(\rho, \phi)|+(M +|D_t \phi|)(|  D_t^2(\rho, \phi)|+|D_t \phi|^2)+M^2 |D_t\rho |,$$
which is due to \ef{1-d}.
The bound for $\| D_t^4 (\rho \Da \phi)\|_{L^2}$ in \ef{8.1-2}
can be obtained similarly by noticing that
\begin{align*}
&| D_t^4 (\rho \Da \phi)|\les |  D_t^4 (\rho, \phi)|+(M +|D_t \phi|)(|  D_t^3 (\rho, \phi)|+|D_t \phi||D_t^2\phi|\\
&\quad +|D_t \phi|^3)+|D_t^2(\rho, \phi)|^2 + (M^2+|D_t \phi|^2)|D_t^2 \rho|
+M^3 |D_t\rho |.
\end{align*}
This finishes the proof of the lemma.
\hfill $\Box$

\subsection{Energy estimates}
\begin{prop}  Let $\ka=1$ in \ef{1-b}, then it holds that
\be\label{3.2-3}
\frac{d}{dt} \mathscr{E}_{EP}\le C(M,\bar M, K,K_1,\ea_b^{-1},{\bar\ea_b}^{-1},L, \bar L) \sum_{1\le i \le 4}  \mathscr{E}_{EP}^i  .
\ee
\end{prop}

The proof consists of the following two lemmas,  Lemmas \ref{lem3.9} and \ref{lem3.10}.

\begin{lem}\label{lem3.9}
Let $\ka=0,1$ in \ef{1-b}, then it holds that for $1\le r\le 4$,
 \begin{align}
 \frac{d}{dt}E_r\le C(M,K,K_1,\ea_b^{-1},L) \sum_{1\le i \le 3} \lt(  \mathscr{E}_I^i
+\ka C(\bar M, {\bar\ea_b}^{-1}, \bar L) \mathscr{E}_{EP}^i \rt),\label{3.2-4}
\end{align}
where
\begin{align*}
& E_r=  \int_{\mathscr{D}_t}
\rho \delta^{mn} \zeta^{IJ} (\pl^r_I v_m)\pl^r_J v_ndx +  \int_{\pl \mathscr{D}_t} |\Pi \pl^r p|^2 (-\pl_N p)^{-1}  ds\\
&\quad  +\ka \int_{\pl \mathscr{D}_t}  \rho |\Pi \pl^r \phi|^2(\pl_N \phi)^{-1}  ds.
\end{align*}
\end{lem}

{\em Proof}. It follows from \ef{7.5-1} that for $r\ge 1$,
$$
D_t \pl^r v + \rho^{-1} \pl^r \pl  p -\ka \pl^r \pl \phi=[D_t, \pl^r] v  +  \rho^{-1} \pl^r \pl  p-  \pl^r (\rho^{-1}  \pl p),
$$
which implies that
\begin{align}
&2^{-1}\rho D_t(\da^{mn} \zeta^{IJ} (\pl^r_I v_m) \pl^r_J v_n)\notag \\
=&\da^{mn} \zeta^{IJ} (\rho D_t \pl^r_I v_m) \pl^r_J v_n +2^{-1}\rho \da^{mn} (D_t \zeta^{IJ}) (\pl^r_I v_m) \pl^r_J v_n\notag \\
=&\mathcal{H}_r-\ka\overline{\mathcal{H}}_r   -{\rm div} \lt(\zeta^{IJ}(\pl^r_I p-\ka  \rho \pl^r_I \phi)\pl^r_J v\rt) ,
\label{3.2-5}
\end{align}
where
\begin{align}
& \mathcal{H}_r= \zeta^{IJ}(\pl^r_I p ) \pl^r_J {\rm div} v
 +\da^{mn}\{2^{-1}\rho (D_t \zeta^{IJ}) \pl^r_I v_m  + (\pl_m \zeta^{IJ}) \pl^r_I p \notag\\
 & \ \    + \zeta^{IJ}\lt( (\pl^r_I \pl_m p- \rho \pl^r_I (\rho^{-1}  \pl_m p)) +\rho [D_t, \pl^r_I] v_m  \rt)   \}\pl^r_J v_n, \label{8.6.1}\\
 & \overline{\mathcal{H}}_r= \{\rho \zeta^{IJ} \pl^r_J {\rm div} v
 +\da^{mn} ( \rho  \pl_m \zeta^{IJ}    + \zeta^{IJ}\pl_m \rho     )\pl^r_J v_n\}  \pl^r_I \phi. \notag
\end{align}
Due to $\pl_m p= (\overline{\pl}_m +  N_m \pl_N ) p= N_m \pl_N p$ on $\pl \mathscr{D}_t$, we have on $\pl \mathscr{D}_t$,
\begin{align}
&2^{-1}D_t\lt((\pl_N p)^{-1}  \zeta^{IJ}(\pl^r_I p )\pl^r_J p \rt)\notag\\
=& (\pl_N p)^{-1} \zeta^{IJ} (\pl^r_I p )D_t \pl^r_J p+ 2^{-1}\lt(D_t((\pl_N p)^{-1} \zeta^{IJ})\rt) (\pl^r_I p )\pl^r_J p
\notag\\
=& \mathcal{L}_r-(\pl_N p)^{-1}
\zeta^{IJ} (\pl^r_I p )
(\pl^r_J v^m) \pl_m p
=  \mathcal{L}_r-   N_m  \zeta^{IJ} (\pl^r_I p )
\pl^r_J v^m,\label{3.2-6}
\end{align}
where
\begin{align}
&\mathcal{L}_r=(\pl_N p)^{-1} \zeta^{IJ}(\pl^r_I p ) (D_t \pl^r_J p- \pl^r_J D_t p + (\pl^r_J v^m) \pl_m p  )
\notag\\
& \quad +(\pl_N p)^{-1}  \zeta^{IJ}(\pl^r_I p ) \pl^r_J D_t p +
2^{-1}\lt(D_t((\pl_N p)^{-1}\zeta^{IJ})\rt) (\pl^r_I p )\pl^r_J p. \label{8.6.2}
\end{align}
Similarly, we have on $\pl \mathscr{D}_t$,
\begin{align}
2^{-1}D_t\lt(\rho (\pl_N \phi)^{-1} \zeta^{IJ} (\pl^r_I \phi )\pl^r_J  \phi \rt)
=& \mathcal{R}_r- \rho N_m \zeta^{IJ} (\pl^r_I \phi )\pl^r_J v^m, \label{22.6.11-2}
\end{align}
where
\begin{align*}
&\mathcal{R}_r=\rho (\pl_N \phi )^{-1}  \zeta^{IJ} (\pl^r_I \phi ) (D_t \pl^r_J \phi- \pl^r_J D_t\phi + (\pl^r_J v^m) \pl_m \phi  )
\\
& \quad +\rho (\pl_N \phi)^{-1}  \zeta^{IJ} (\pl^r_I \phi ) \pl^r_J D_t \phi +
2^{-1}\lt(D_t(\rho (\pl_N \phi)^{-1} \zeta^{IJ})\rt) (\pl^r_I \phi )\pl^r_J \phi.
\end{align*}
In view of $\int_{\mathscr{D}_t}   {\rm div} \lt(\zeta^{IJ}(\pl^r_I p-\ka  \rho \pl^r_I \phi)\pl^r_J v\rt)  dx
= \int_{\pl \mathscr{D}_t}   N_m
\zeta^{IJ}(\pl^r_I p-\ka  \rho \pl^r_I \phi )\pl^r_J v^m
ds$, \ef{3.2-1-b}, \ef{3.2-2}, \ef{3.2-5}, \ef{3.2-6} and \ef{22.6.11-2},
we see that for $r\ge 1$,
\begin{align}
&\frac{d}{dt}E_r=
\int_{\pl \mathscr{D}_t}D_t (|\Pi \pl^r p|^2 (-\pl_N p)^{-1}  + \ka \rho |\Pi \pl^r \phi|^2(\pl_N \phi)^{-1} )ds
\notag\\
&-\int_{\pl \mathscr{D}_t} (|\Pi \pl^r p|^2 (-\pl_N p)^{-1}  + \ka \rho |\Pi \pl^r \phi|^2(\pl_N \phi)^{-1} ) N^i \pl_N  v_ids
\notag\\
&+\int_{\mathscr{D}_t}  \rho  D_t(\da^{mn} \zeta^{IJ} (\pl^r_I v_m) \pl^r_J v_n) dx
=\mathcal{G}_r, \label{22.6.11-3}
\end{align}
where
\begin{align*}
&\mathcal{G}_r=2 \int_{\mathscr{D}_t} (\mathcal{H}_r- \ka \overline{\mathcal{H}}_r) dx
+ 2  \int_{\pl \mathscr{D}_t}
(\ka \mathcal{R}_r - \mathcal{L}_r) ds
\\
&-\int_{\pl \mathscr{D}_t} (|\Pi \pl^r p|^2 (-\pl_N p)^{-1}  + \ka \rho |\Pi \pl^r \phi|^2(\pl_N \phi)^{-1} ) N^i \pl_N v_i  ds.\end{align*}

Due to \ef{CL-3.28} and \ef{3-7-3}, one has that for $1\le r\le 4$,
\begin{align*}
&|\mathcal{H}_r|\les |\pl^r p ||\pl^r{\rm div} v|+(M|\pl^r v | + K|\pl^r p|+H_r )|\pl^r v|,\\
&|\overline{\mathcal{H}}_r|\les
(|\pl^r{\rm div} v|+ (K+M)|\pl^r v|)|\pl^r \phi|,
\end{align*}
where
\begin{align*}
& H_1 =
 M|\pl (v,p)|,
 \ \  H_2 = M|\pl^2 (v,\rho,p)|+M^2|\pl p|,  \\
& H_3 = M |\pl^3 (v,\rho,p)|+ M^2|\pl^2 (\rho,p)|+M^3|\pl p|+|\pl^2 v|^2+|\pl^2 \rho||\pl^2 p|,  \\
& H_4 = M |\pl^4 (v,\rho,p)|+ M^2|\pl^3 (\rho,p)|+ M^3|\pl^2 (\rho,p)|+M^4|\pl p|\\
&\quad +|\pl^2 v| |\pl^3 v|+|\pl^2 \rho||\pl^3 p|+|\pl^3 \rho||\pl^2 p|+M|\pl^2 \rho||\pl^2 (\rho,p)|.
\end{align*}
This, together with \ef{lem-1}, \ef{7.10.3} and \ef{22.lem3}, gives that for $1\le r\le 4$,
\begin{subequations}\label{3.2-8}\begin{align}
&\| \mathcal{H}_r\|_{L^1}\le C(M,K,K_1,\ea_b^{-1},L)\lt(\mathscr{E}_I+\mathscr{E}_I^2
+\ka C(\bar M,  \bar L)\sum_{1\le i \le 3}\mathscr{E}_{EP}^i \rt),\\
& \|\overline{\mathcal{H}}_r\|_{L^1}\le
C(M,\bar M,K,K_1,\ea_b^{-1},L,\bar L)\lt(\mathscr{E}_{EP}+\mathscr{E}_{EP}^2+\mathscr{E}_{EP}^3\rt).
\end{align}\end{subequations}
Note that
\begin{align*}
&|D_t \pl^r  p- \pl^r  D_t p + (\pl^r  v^m) \pl_m p|
\les \sum_{1\le i \le r-1}|\pl^i v||\pl^{r+1-i} p|
,\\
&\lt||\pl^2 v|(|\pl^2 p|+|\pl^3 p|)\rt|_{L^2}
+\lt||\pl^3 v||\pl^2 p|\rt|_{L^2}\\
&\quad \les |\pl^2 v|_{L^4}(|\pl^2 p|_{L^4}+|\pl^3 p|_{L^4})
+|\pl^3 v|_{L^4}|\pl^2 p|_{L^4}\\
& \quad \le C(K) \sum_{2\le i \le 4} \|\pl^i v\|_{L^2}\sum_{2\le j \le 4} \|\pl^j p\|_{L^2}\le C(K) \mathscr{E}_I,
\end{align*}
which is due to \ef{CL-A.20} and \ef{7.6.1}.  Then, we use \ef{7.7.1}, \ef{lem5-b} and \ef{22.lem3} to get that for $1\le r \le 4$,
\begin{align}\label{8.1-3}
& |D_t \pl^r  p- \pl^r  D_t p + (\pl^r  v^m) \pl_m p|_{L^2}^2
\notag\\
& \le C(M,K,K_1,\ea_b^{-1})  \lt(\sum_{1\le i \le 3}  \mathscr{E}_I^i
+\ka C(\bar M,  \bar L) (\mathscr{E}_{EP}+\mathscr{E}_{EP}^2) \rt).
\end{align}
In view of \ef{dtn} and \ef{3-7-3}, we see that $|D_t \zeta^{IJ}|\les |\pl v| \le M $ and
\begin{align*}
|D_t(\pl_N p)^{-1}| & \le | (\pl_N p)^{-2}|(|D_t N^i||\pl_i p|
+|\pl_N D_t p|+|N^i[D_t, \pl_i]p|)\notag\\
& \les \ea_b^{-2}(M^2+L) \ \  {\rm on} \ \  \pl \mathscr{D}_t,
\end{align*}
which implies, using \ef{7.7.1}, \ef{7.10.3}, \ef{lem5-b}, \ef{22.lem3} and \ef{8.1-3},  that for $1\le r\le 4$,
\begin{align}
& |\mathcal{L}_r|_{L^1 }  \le C(  \ea_b^{-1},M, L)(|\pl^r p|_{L^2}^2+ |\Pi \pl^r D_t p|_{L^2}^2 \notag\\
&\quad +|D_t \pl^r  p- \pl^r  D_t p + (\pl^r  v^m) \pl_m p|_{L^2}^2)\notag \\
&\quad \le C(M,K,K_1,\ea_b^{-1},L) \sum_{1\le i \le 3} \lt(  \mathscr{E}_I^i
+\ka C(\bar M,  \bar L) \mathscr{E}_{EP}^i \rt) .\label{8.1-4}
\end{align}
Similarly, we have $|D_t(\pl_N \phi)^{-1}|   \les \bar{\ea}_b^{-2}(M \bar M+\bar{L})$ and $D_t \rho=\rho_p D_t p=0$ on $\pl \mathscr{D}_t$,
and use  \ef{7.7.8}  and \ef{22.lem3} to get that for $1\le r\le 4$,
\begin{align}
& |\mathcal{R}_r|_{L^1 }  \le C(  \bar\ea_b^{-1},M, \bar M, \bar L)(|\pl^r \phi|_{L^2}^2+ |\Pi \pl^r D_t \phi|_{L^2}^2 \notag\\
&\quad +|D_t \pl^r  \phi- \pl^r  D_t \phi + (\pl^r  v^m) \pl_m \phi|_{L^2}^2)\notag \\
&\quad \le C( M, \bar M, K, K_1, \ea_b^{-1}, \bar\ea_b^{-1}, \bar L)
(\mathscr{E}_{EP}+ \mathscr{E}_{EP}^2 ).\notag
\end{align}
This proves \ef{3.2-4}, with the help of
\ef{22.6.11-3}, \ef{3.2-8} and \ef{8.1-4}.
\hfill $\Box$

\begin{lem}\label{lem3.10}
Let $\ka=0,1$ in \ef{1-b}, then it holds that
for $0\le r \le 4$,
\begin{align}\label{8.2}
\frac{d}{dt}(P_r + W_r) \le
C(M,K,K_1,\ea_b^{-1},L)\sum_{1\le i\le 4}\lt(  \mathscr{E}_I^i
+\ka C(\bar M,  \bar L) \mathscr{E}_{EP}^i \rt) ,
\end{align}
where
\begin{align*}
&P_r=\int_{\mathscr{D}_t}  \lt( |D_t^{r+1} \rho|^2+  \rho_p |\pl D_t^{r}p|^2\rt) dx, \ r\ge 0, \\
&W_0= \int_{\mathscr{D}_t} \lt(\rho |v|^2+ \rho^2
+p^2+ s^2 +\ka \phi^2  \rt)dx ,\\
&W_r= \int_{\mathscr{D}_t}
\lt(|\pl^{r-1}({\rm curl} v, {\rm div} v)|^2+ |\pl^r (\rho, p,s)|^2 +\ka |\pl^r \phi |^2   \rt)dx , \  r\ge 1.
\end{align*}

\end{lem}

{\em Proof}. We use \ef{3.2-1-a} and
$|D_t\pl D_t^{r}p-\pl D_t^{r+1}p|
\les M |\pl D_t^{r}p|$
to obtain that for $0\le r\le 3$,
\begin{align}
&\frac{d}{dt} P_r \le \int_{\mathscr{D}_t} \{
2   |D_t^{r+1} \rho| |D_t^{r+2} \rho|
+2 \rho_p |\pl D_t^{r}p||D_t \pl D_t^{r}p| +
(D_t \rho_p) |\pl D_t^{r}p|^2
 \}dx\notag\\
&+M P_r \le  \int_{\mathscr{D}_t}(
     |D_t^{r+2} \rho|^2
+  \rho_p |D_t \pl D_t^{r}p|^2 )dx
+C(M) P_r \notag \\
& \le 2 P_{r+1} +C(M) P_r\le C(M)\mathscr{E}_I.\label{8.2-1}
\end{align}
For $P_4$, we
notice that
\begin{align*}
& (\Da D_t^4 p )  D_t^5 \rho
={\rm div}\lt( (D_t^5 \rho) \pl D_t^4 p \rt)
- \pl D_t^5 \rho \cdot \pl D_t^4 p\\
&  ={\rm div}\lt( (D_t^5 \rho) \pl D_t^4 p \rt)- (\pl D_t^5 \rho - \rho_p D_t \pl D_t^4 p )   \cdot \pl D_t^4 p\\
&\quad - 2^{-1} D_t(\rho_p |\pl D_t^4 p|^2)
+2^{-1} (D_t \rho_p) |\pl D_t^4 p|^2
.
\end{align*}
Then, it follows from $D_t^5 p=0$ on $\pl \mathscr{D}_t $,  \ef{3.2-1-a}, \ef{7.8.3}, \ef{lem5-c} and \ef{22.lem3} that
\begin{align}
&\frac{d}{dt} P_4
 \le \int_{\mathscr{D}_t} \lt(  2 (D_t^{6} \rho)D_t^5 \rho +  D_t (\rho_p |\pl D_t^{4}p|^2)\rt) dx
+MP_4
\notag\\
\le & 2 \int_{\mathscr{D}_t } ( (D_t^6 \rho - \Da D_t^4 p) D_t^5 \rho + {\rm div}\lt( (D_t^5 \rho ) \pl D_t^4 p \rt)-(\pl D_t^5 \rho
\notag\\
&- \rho_p D_t \pl D_t^4 p )   \cdot \pl D_t^4 p+ 2^{-1} (D_t \rho_p) |\pl D_t^4 p|^2 )dx
 +MP_4
 \notag\\
\les & \int_{\mathscr{D}_t }  (|D_t^6 \rho - \Da D_t^4 p|^2 +
|\pl D_t^5 \rho - \rho_p D_t \pl D_t^4 p |^2 )dx
+ (M+1)P_4 \notag \\
 \le & C(M,K,K_1,\ea_b^{-1})\sum_{1\le i\le 4}\lt(  \mathscr{E}_I^i
+\ka C(\bar M,  \bar L) \mathscr{E}_{EP}^i \rt). \label{8.2-2}
\end{align}

In view of  \ef{3.2-1-a} and \ef{3.2-2}, we see  that
\begin{align}
&\frac{d}{dt}W_0\le 2\int_{\mathscr{D}_t} \lt(\rho |v||D_t v|+|\rho D_t \rho|
+|p D_t p|  +\ka |\phi D_t \phi|  \rt)dx+M W_0\notag\\
\le & \int_{\mathscr{D}_t} \lt(\rho |D_t v|^2+|D_t \rho |^2
+|D_t p |^2  +\ka |  D_t \phi|^2  \rt)dx+(M+1) W_0\notag\\
\le & C(M,K_1)\lt(\mathscr{E}_I+\ka C(\bar M, K) \mathscr{E}_{EP} \rt), \label{8.2-3}
\end{align}
where we have used  \ef{1-c}  to derive the first inequality;  \ef{7.6.1}, \ef{7.8.5} and \ef{22Mar9-4} to the last inequality.
It follows from \ef{3.2-1-a} that for $1\le r\le 4$,
\begin{align}
&\frac{d}{dt} W_r
 \le  2 \int_{\mathscr{D}_t}
\lt(|\pl^{r-1}({\rm curl} v, {\rm div} v)||D_t\pl^{r-1}({\rm curl} v, {\rm div} v)| \rt. \notag\\
&\lt. + |\pl^r (\rho, p,s)| |D_t \pl^r (\rho, p,s)| +\ka |\pl^r \phi ||D_t \pl^r \phi |  \rt)dx
 +M W_r\notag \\
& \le   \int_{\mathscr{D}_t}
\lt(|D_t\pl^{r-1}({\rm curl} v, {\rm div} v)|^2   +   |D_t \pl^r (\rho,p,s)|^2 +\ka  |D_t \pl^r \phi |^2  \rt)dx
 \notag\\
 &+(M+1) W_r
 \le C(M,K,K_1,\ea_b^{-1},L) \sum_{1\le i \le 3} \lt(  \mathscr{E}_I^i
+\ka C(\bar M,   \bar L) \mathscr{E}_{EP}^i \rt),
\label{8.2-4}
\end{align}
where we have used \ef{7.7.1}, \ef{7.7.2}, \ef{7.8.5},
 \ef{7.8.4}, \ef{7.10.3},  \ef{lem5-a} and \ef{22.lem3} to derive the last inequality.
As a consequence of \ef{8.2-1}-\ef{8.2-4}, we prove \ef{8.2} and finish the proof of the lemma. \hfill $\Box$

\subsection{Proof of Theorem \ref{mainthm}}
The proof  follows from Proposition \ref{prop23.2.20}, which is stated as follows.

\begin{prop}\label{prop23.2.20}
Let $\ka=1$ in \ef{1-b}, then there exists a continuous function
$\overline{T}>0$ such that
\begin{align*}
& 2^{-1}{\rm Vol} \mathscr{D}_0\le  {\rm Vol} \mathscr{D}_t\le 2 {\rm Vol} \mathscr{D}_0 , \label{}\\
 & 2^{-1}\underline{\varrho} \le \rho( t,x)\le 2\overline{\varrho}, \ \  |s(t,x)|\le \overline{s}, \ \  x\in \mathscr{D}_t,
 \label{}\\
 &
-\pl_N p(t,x)\ge   2^{-1}\varepsilon_1, \  \
\pl_N \phi(t,x)\ge   2^{-1}\varepsilon_2,  \ \ x\in  \pl \mathscr{D}_t,  \\
& \iota_1^{-1}(t)\le   16K_0,\ \
 \mathscr{E}_{EP}(t)\le 2 \mathscr{E}_{EP}(0),\\
& \|   \pl  (v, p, s )(t,\cdot) \|_{L^\infty} \le 2 \|   \pl  ( v, p, s)(0,\cdot) \|_{L^\infty}, \\
&
\|   \pl   \phi(t,\cdot) \|_{L^\infty} \le 2 \|   \pl    \phi(0,\cdot) \|_{L^\infty},\\
&|\theta(t,\cdot)|_{L^\iy}+\iota_0^{-1}(t)
  \le C( \underline{\varrho}^{-1},\overline{\varrho},
   \varepsilon_2^{-1}, K_0, \mathscr{E}_{EP}(0),
{\rm Vol} \mathscr{D}_0), \\
&
 |\pl_N D_t \phi(t,\cdot)|_{L^\iy} + |\pl_N D_t p(t,\cdot)|_{L^\iy} \notag \\
& \le C\lt(
\underline{\varrho}^{-1},\overline{\varrho},
\overline{s}, \varepsilon_1^{-1}, \varepsilon_2^{-1}, K_0, \mathscr{E}_{EP}(0),
{\rm Vol} \mathscr{D}_0\rt)
\end{align*}
for $t\le \overline{T}( \underline{\varrho}^{-1},\overline{\varrho},
\overline{s}, \varepsilon_1^{-1}, \varepsilon_2^{-1}, K_0, \mathscr{E}_{EP}(0),
{\rm Vol} \mathscr{D}_0)$.
\end{prop}

{\em Proof}.  We use the Lagrangian map:  let  $x= x(t,y) $ be the change of variables given by
$$
\pl_t x(t,y) =   v \left(t,  x(t,y)\right)  \ \  {\rm and} \ \ x(0,y) =y, \ \  y \in \mathscr{D}_0.
$$
 For each $t$ we will then have a change of coordinates $x: \mathscr{D}_0
\to \mathscr{D}_t $, taking $y \to x(t,y)$.
The proof consists of three steps.

{\em Step 1}. We prove that for $ t \le M^{-1} \ln 2$,
\begin{subequations}
\label{22.6.28}\begin{align}
& 2^{-1}{\rm Vol} \mathscr{D}_0\le  {\rm Vol} \mathscr{D}_t\le 2 {\rm Vol} \mathscr{D}_0 , \label{23.2.18}\\
& |s(t,x)|\le \overline{s},
 \ \  x\in \mathscr{D}_t,\label{22.6.20-1}\\
 & 2^{-1} \underline{\varrho} \le \rho(t, x) \le 2\overline{\varrho}, \ \  x\in \mathscr{D}_t.
 \label{22.6.20-5}
\end{align}
\end{subequations}
It follows from \ef{3.2-1} that
$$\lt|\frac{d}{dt}{\rm Vol} \mathscr{D}_t\rt| = \lt|\frac{d}{dt}  \int_{\mathscr{D}_t } 1 dx\rt|=\lt|\int_{\mathscr{D}_t }  {\rm div} v  dx\rt|\le M {\rm Vol} \mathscr{D}_t,  $$
which implies that
$$
e^{-Mt}{\rm Vol} \mathscr{D}_0\le  {\rm Vol} \mathscr{D}_t\le e^{Mt}{\rm Vol} \mathscr{D}_0 .  $$
Thus, we have \ef{23.2.18}  for $  t \le M^{-1} \ln 2$.
Since $D_t s=0$, we have
$
s(t,x(t, y))=s(0,x(0, y))=s_0(y)$ and then \ef{22.6.20-1} holds for $  t \le T$.
In view of \ef{1-a}, we see that
$$ \rho(t, x(t,y))=\rho_0(y)e^{-\int_0^t {\rm div } v( \tau, x( \tau,y))d\tau}, $$
which implies
$$e^{-Mt}\rho_0(y) \le \rho(t, x(t,y)) \le e^{Mt}\rho_0(y) . $$
Then,  \ef{22.6.20-5} holds for $ t \le M^{-1} \ln 2$. Hence, we prove \ef{22.6.28}.

{\em Step 2}. We prove that there exists a continuous function
$\widetilde{T}>0$ such that
\begin{subequations}\label{22.6.30}\begin{align}
&
-\pl_N p(t,x)\ge   2^{-1}\varepsilon_1, \  \
\pl_N \phi(t,x)\ge   2^{-1}\varepsilon_2,  \ \ x\in  \pl \mathscr{D}_t,  \\
& \iota_1^{-1}(t)\le   16K_0,\ \
 \mathscr{E}_{EP}(t)\le 2 \mathscr{E}_{EP}(0),\\
& \|   \pl  (v, p, s )(t,\cdot) \|_{L^\infty} \le 2 \|   \pl  ( v, p, s)(0,\cdot) \|_{L^\infty}, \\
&
\|   \pl   \phi(t,\cdot) \|_{L^\infty} \le 2 \|   \pl    \phi(0,\cdot) \|_{L^\infty},
\end{align}\end{subequations}
for $t\le \widetilde{T}(  K,   L, \bar L, \underline{\varrho}^{-1}, \overline{\varrho}, \overline{s}, \varepsilon_1^{-1}, \varepsilon_2^{-1}, K_0,\mathscr{E}_{EP}(0),  {\rm Vol} \mathscr{D}_0) $.
Moreover, \ef{22.6.28} also holds for $t\le \widetilde{T}$.

Due to \ef{dtn}, we have on $\pl \mathscr{D}_t$,
$$
|D_t \pl_N p| \le |D_t N^i||\pl_i p|+|\pl_N D_t p|+|N^i[D_t, \pl_i] p|
\le 3M^2+L,$$
which gives, using \ef{in2}, that   for $  t \le 2^{-1}(3M^2+L)^{-1}\varepsilon_1 $  and $x\in \pl\mathscr{D}_t$,
\begin{align}
&-\pl_N p(t,x(t,y))=- \pl_N p (0,x(0,y))
-\int_0^t  D_\tau( \pl_N p (\tau,x(\tau,y))) d\tau
\notag\\
&\qquad \ge \varepsilon_1 - (3M^2+L)t \ge 2^{-1}\varepsilon_1.
\label{22.6.30-1}
\end{align}
Similarly, one has  that  for $  t \le 2^{-1}(3\bar M^2+\bar L)^{-1}\varepsilon_2 $  and $x\in \pl\mathscr{D}_t$,
\be\label{22.6.30-2}
\pl_N \phi(t,x)\ge   2^{-1}\varepsilon_2.
\ee

Let $\epsilon_1\in (0,1/2]$ be a fixed constant (for example, $\epsilon_1=1/4$), $\underline{\iota_1}$ the largest number such that
\begin{equation}\label{CL-7.55}\begin{split}
&| {N}(0, x(0,y_1))-{N}(0, x(0,y_2))|\le 2^{-1}\epsilon_1 ,
\\
&  \ \  \ \ \ \ {\rm whenever} \ \
|x(0,y_1) -x(0,y_2)| \le  2\underline{\iota_1},  \ \
 y_1, y_2 \in  \pl\mathscr{D}_0.
\end{split}\end{equation}
Then we have from   \eqref{lemkk1}  and \eqref{in3} that
\begin{equation}\label{iota10}
\underline{\iota_1} \ge  \epsilon_1/(4K_0) \ \ {\rm and} \ \ \iota_1(0) \ge \ea_1/K_0
.
\end{equation}
It follows from \ef{7.5-1} and \ef{22.6.20-5} that
$|D_t v|\le 2 \underline{\varrho}^{-1}M+ \bar M$
and
\begin{align*}
|v(t,x(t,y))| \le
|v_0(y)|+ (2 \underline{\varrho}^{-1}M+ \bar M) t\le
 2 |v_0|_{L^\iy}
\end{align*}
for $t\le \min\{M^{-1}\ln 2, \ (2 \underline{\varrho}^{-1}M+ \bar M)^{-1}|v_0|_{L^\iy} \}$ and $y\in \pl \mathscr{D}_0$. This implies
\begin{align}\label{6.30-1}
|x(t,y)-x(0,y)|\le \int_0^t |v(\tau, x(\tau, y))|d \tau \le 2 |v_0|_{L^\iy} t \le 2^{-1}\underline{\iota_1}
\end{align}
for $t\le T_1=\min\{M^{-1}\ln 2, \ (2 \underline{\varrho}^{-1}M+ \bar M)^{-1}|v_0|_{L^\iy}, \  (4|v_0|_{L^\iy})^{-1} \underline{\iota_1} \}$ and $y\in \pl \mathscr{D}_0$.
In view of
\ef{dtn}, we see that $|D_t N|\le 2|\pl v|\le 2M$ and
$$
|N(t,x(t,y))-N(0, x(0,y))|\le   2 M t \le 4^{-1} \ea_1
$$
for $t\le  (8M)^{-1}\ea_1 $ and $y\in \pl \mathscr{D}_0$.
This, together with \ef{CL-7.55} and \ef{6.30-1}, gives
\begin{equation*}\label{CL-7.56}\begin{split}
&|{N}(t, x(t,y_1))- {N}(t, x(t,y_2))|\le \epsilon_1  ,
\\
& \ \  \ \  \ \  {\rm whenever} \ \
|x(t,y_1) -x(t,y_2)| \le   \underline{\iota_1},  \ \
 y_1, y_2 \in \pl\mathscr{D}_0
\end{split}\end{equation*}
for $t\le T_2=\min\{T_1,\ (8M)^{-1}\ea_1\}$.
So, it yields from \ef{iota10} that for $t\le T_2$,
\be\label{22.6.23}
\iota_1(t)\ge \underline{\iota_1}\ge \epsilon_1/(4K_0).
\ee

It follows from \ef{3.2-3} and \ef{22.6.28} that
there exists a continuous function $T_3\in (0, M^{-1} \ln 2]$ such that
\begin{align}
\mathscr{E}_{EP}(t)\le 2 \mathscr{E}_{EP}(0) \label{22.6.27-1}
\end{align}
for $t \le T_3( M, \bar M, K, K_1, \ea_b^{-1}, \bar\ea_b^{-1},L, \bar L, \underline{\varrho}^{-1}, \overline{\varrho}, \overline{s},{\rm Vol} \mathscr{D}_0)$.
This means, with the aid of
\ef{7.7.1}, \ef{7.8.5} and \ef{22Mar9-5}, that for $t\le T_3$,
\begin{align*}
\| D_t \pl  (v,  p, s,  \phi)(t,\cdot) \|_{L^\infty}
\le C( M, \bar M, K, K_1, \ea_b^{-1}, \bar L,\mathscr{E}_{EP}(0)).
\end{align*}
Thus, there exists a continuous function $T_4\in (0, T_3]$ such that
\begin{subequations}\label{22.6.27-2}\begin{align}
&\|   \pl  (  v, p, s)(t,\cdot) \|_{L^\infty} \le 2
\|   \pl  (  v, p, s)(0,\cdot) \|_{L^\infty}, \\
& \|   \pl   \phi(t,\cdot) \|_{L^\infty} \le 2
\|   \pl    \phi(0,\cdot) \|_{L^\infty},
\end{align}\end{subequations}
for $t\le T_4( M, \bar M, K, K_1, \ea_b^{-1}, \bar\ea_b^{-1},L, \bar L,  \mathscr{E}_{EP}(0),\|   \pl  ( v, p, s, \phi)(0,\cdot) \|_{L^\infty},\underline{\varrho}^{-1}, \overline{\varrho}, \overline{s},{\rm Vol} \mathscr{D}_0)$.

In view of   \ef{CL-A.8}, \ef{CL-A.15},
\ef{23.2.19}, \ef{in1}, \ef{22.6.20-5}
 and \eqref{iota10}, we see that
\be\label{7.1-1}
|v_0|_{L^\iy}^2+ \|   \pl  (v, p, s, \phi)(0,\cdot) \|_{L^\infty}^2 \le C(K_0/\ea_1, \underline{\varrho}^{-1})  \mathscr{E}_{EP}(0).
\ee
This, together with  \ef{22.6.30-1}, \ef{22.6.30-2} and
\ef{22.6.23}-\ef{22.6.27-2}, implies that
\ef{22.6.30} holds for $t\le \widetilde{T}$ for some continuous function $\widetilde{T}(  K,   L, \bar L, K_0,\mathscr{E}_{EP}(0), \varepsilon_1^{-1}, \varepsilon_2^{-1},  \underline{\varrho}^{-1}, \overline{\varrho}, \overline{s},{\rm Vol} \mathscr{D}_0)>0$,
 by choosing
\begin{subequations}\label{7.1-2}\begin{align}
& \epsilon_b= \varepsilon_1/4, \  \ \bar\epsilon_b= \varepsilon_2/4, \ \ M=4 \|   \pl  (v, p, s) (0,\cdot) \|_{L^\infty}, \\
& \bar M=4 \|   \pl  \phi(0,\cdot) \|_{L^\infty} ,  \ \  \ea_1=1/4, \ \  K_1=32K_0.
\end{align}\end{subequations}
Clearly, \ef{22.6.28} holds for $t\le \widetilde{T}$.

{\em Step 3}. We prove that there exists a continuous function
$\overline{T}>0$ such that
\begin{subequations}\label{22.7.1}\begin{align}
&
|\theta(t,\cdot)|_{L^\iy}+\iota_0^{-1}(t)
  \le C( \underline{\varrho}^{-1},\overline{\varrho},
   \varepsilon_2^{-1}, K_0, \mathscr{E}_{EP}(0),
{\rm Vol} \mathscr{D}_0),\\
&
|\pl_N D_t \phi(t,\cdot)|_{L^\iy} + |\pl_N D_t p(t,\cdot)|_{L^\iy} \notag \\
 &\quad \le C( \underline{\varrho}^{-1},\overline{\varrho},
\overline{s}, \varepsilon_1^{-1}, \varepsilon_2^{-1}, K_0, \mathscr{E}_{EP}(0),
{\rm Vol} \mathscr{D}_0)
\end{align}\end{subequations}
for $t\le \overline{T}( \underline{\varrho}^{-1},\overline{\varrho},
\overline{s}, \varepsilon_1^{-1}, \varepsilon_2^{-1}, K_0, \mathscr{E}_{EP}(0),
{\rm Vol} \mathscr{D}_0)$.
Moreover, \ef{22.6.28} and \ef{22.6.30} also hold for $t\le \overline{T}$.

It follows from \ef{22.5.7-1}, \ef{22.6.26}, \ef{22.6.30}, \ef{7.1-1} and \ef{7.1-2}  that
for $t\le \widetilde{T}$,
\begin{align}
&|\theta(t,\cdot)|_{L^\iy}=|(\pl_N \phi)^{-1}(t,\cdot)|_{L^\iy}|\Pi \pl^2 \phi(t,\cdot)|_{L^\iy}\notag\\
& \le |(\pl_N \phi)^{-1}(t,\cdot)|_{L^\iy}| \pl^2 \phi(t,\cdot)|_{L^\iy}
\le C(\bar\epsilon_b^{-1}, M, \bar M, K_1,\mathscr{E}_{EP})\notag\\
& \le C_1( \underline{\varrho}^{-1},\overline{\varrho},
   \varepsilon_2^{-1}, K_0, \mathscr{E}_{EP}(0),
{\rm Vol} \mathscr{D}_0),
\label{7.1-3}
\end{align}
which implies, with the help of \ef{lemkk1} and  \ef{22.6.30}, that for $t\le \widetilde{T}$,
\be\label{7.1-4}
\iota_0^{-1}(t) \le \max\{32K_0,\ \
C_1( \underline{\varrho}^{-1},\overline{\varrho},
   \varepsilon_2^{-1}, K_0, \mathscr{E}_{EP}(0),
{\rm Vol} \mathscr{D}_0)\}.
\ee
So, we may choose
\be\label{7.1-5}
K=64K_0+4C_1( \underline{\varrho}^{-1},\overline{\varrho},
   \varepsilon_2^{-1}, K_0, \mathscr{E}_{EP}(0),
{\rm Vol} \mathscr{D}_0).
\ee
With \ef{CL-5.34}, \ef{7.6.3}, \ef{22.6.28}, \ef{22.6.30}, \ef{7.1-1}, \ef{7.1-2} and \ef{7.1-5} at hand, we use   \ef{22Mar9-4} and
\ef{22-3-9-3} to get that for $t\le \widetilde{T}$,
\begin{align}\label{7.1-6}
&|\pl_N D_t \phi(t,\cdot)|_{L^\iy}\le C  ( M, \bar M, K,K_1,\ea_b^{-1}, {\rm Vol}\mathscr{D}_t, \mathscr{E}_{EP})\notag\\
&\le C_2( \underline{\varrho}^{-1},\overline{\varrho},
 \varepsilon_1^{-1},  \varepsilon_2^{-1}, K_0, \mathscr{E}_{EP}(0),
{\rm Vol} \mathscr{D}_0).
\end{align}
Thus, we can choose
\be\label{7.1-7}
\bar L= 2 C_2( \underline{\varrho}^{-1},\overline{\varrho},
 \varepsilon_1^{-1},  \varepsilon_2^{-1}, K_0, \mathscr{E}_{EP}(0),
{\rm Vol} \mathscr{D}_0).
\ee
Similarly, we use \ef{7.7.1}, \ef{7.7.2}, \ef{22.5.1}, \ef{22.lem3} and \ef{7.1-7} to obtain that for $t\le \widetilde{T}$,
\begin{align}
& |\pl_N D_t p(t,\cdot)|_{L^\iy}\le
C  ( M, \bar M, K,K_1,\ea_b^{-1}, \bar L, {\rm Vol}\mathscr{D}_t, \mathscr{E}_{EP})\notag\\
&\le C_3( \underline{\varrho}^{-1},\overline{\varrho},
\overline{s}, \varepsilon_1^{-1}, \varepsilon_2^{-1}, K_0, \mathscr{E}_{EP}(0),
{\rm Vol} \mathscr{D}_0)
\label{7.1-8}
\end{align}
and then choose
\be
L=2 C_3( \underline{\varrho}^{-1},\overline{\varrho},
\overline{s}, \varepsilon_1^{-1}, \varepsilon_2^{-1}, K_0, \mathscr{E}_{EP}(0),
{\rm Vol} \mathscr{D}_0).
\label{7.1-9}
\ee

It yields from \ef{7.1-3}, \ef{7.1-4}, \ef{7.1-6} and \ef{7.1-8} that \ef{22.7.1} holds for $t\le \overline{T}$
for some continuous function $\overline{T}( \underline{\varrho}^{-1},\overline{\varrho},
\overline{s}, \varepsilon_1^{-1}, \varepsilon_2^{-1}, K_0, \mathscr{E}_{EP}(0),
{\rm Vol} \mathscr{D}_0)>0$ by choosing $K, \bar L$ and $L$ as in \ef{7.1-5}, \ef{7.1-7} and \ef{7.1-9}, respectively. Clearly,
\ef{22.6.28} and \ef{22.6.30} also hold for $t\le \overline{T}$.
This finishes the proof of the proposition. \hfill $\Box$

\section{The non-isentropic compressible Euler    equations}
\label{sec4}
In this section,
we study the free boundary problem for the non-isentropic compressible Euler equations,  \ef{1} with $\kappa=0$, under the stability condition \ef{tls-1}. The main results are given in Theorem \ref{thm-2}. For this, we define the higher-order energy functionals as follows:
\begin{align}
& \mathscr{E}_{E}(t)= \int_{\mathscr{D}_t} \rho |v|^2 dx
+\sum_{1\le r \le 5} \int_{\pl \mathscr{D}_t} |\Pi \pl^r p|^2 (-\pl_N p)^{-1}  ds
\notag\\
&
\quad +\sum_{1\le r \le 5} \int_{\mathscr{D}_t}
\lt(\rho \delta^{mn} \zeta^{IJ} (\pl^r_I v_m)\pl^r_J v_n+ |\pl^{r-1}{\rm curl} v|^2+|\pl^{r-1} {\rm div} v|^2\rt)dx\notag\\
&
\quad + \sum_{0\le r \le 5} \int_{\mathscr{D}_t}  (|\pl^r \rho|^2 + |\pl^r p|^2+ |\pl^r s|^2  +|D_t^{r+1} \rho|^2+  \rho_p |\pl D_t^{r}p|^2  )dx. \notag
\end{align}
In order to state the main results, we set
\begin{align*}
&\underline{\varrho}=\min_{x\in \mathscr{D}_0} \rho_0(x), \  \  \overline{\varrho}=\max_{x\in \mathscr{D}_0}\rho_0(x),   \ \,
\overline{s}=\max_{x\in \mathscr{D}_0} |s_0(x)|,
\label{}\\
& \varepsilon_1= \min_{x\in \pl\mathscr{D}_0}(-\pl_N p)( 0,x) , \ \
K_0= \max_{x\in \pl\mathscr{D}_0} |\theta(0,x)|
+|{\iota_0}^{-1}(0)|,
\label{}
\end{align*}
where $p(0,x)=p(\rho_0(x), s_0(x))$. With these notations, the main results of this section are stated as follows:

\begin{thm}\label{thm-2}
Let $\ka=0$ in \ef{1-b}, and \eqref{conP} hold.
Suppose that
$$0<\underline{\varrho},  \overline{\varrho}, \overline{s},
\varepsilon_1,  K_0,  {\rm Vol} \mathscr{D}_0, \mathscr{E}_{E}(0)<\iy.$$
Then there exists a continuous function
$\mathscr{T}\lt(  \underline{\varrho}^{-1},\overline{\varrho},
\overline{s}, \varepsilon_1^{-1}, K_0, \mathscr{E}_{E}(0),
{\rm Vol} \mathscr{D}_0\rt)>0$
such that  any smooth solution of the free boundary problem \eqref{1}  for $0\le t\le T$ with $T\le \mathscr{T}$
satisfies the  following estimates:  for $0\le t\le T$,
\begin{subequations}\label{23.thm-2}\begin{align}
&\mathscr{E}_{E}(t)\le 2 \mathscr{E}_{E}(0),  \ \  2^{-1}{\rm Vol} \mathscr{D}_0\le  {\rm Vol} \mathscr{D}_t\le 2 {\rm Vol} \mathscr{D}_0 , \label{}\\
&2^{-1} \underline{\varrho} \le \min_{x\in \mathscr{D}_t} \rho(t, x), \  \  \max_{x\in \mathscr{D}_t}\rho(t, x)\le 2\overline{\varrho},  \label{23.thm-2-d}\\
  & \max_{x\in \mathscr{D}_t}|s(t,x)|\le \overline{s},
 \ \
 2^{-1}\varepsilon_1\le \min_{x\in \pl\mathscr{D}_t}(-\pl_N p)( t,x) ,   \label{}\\
&\max_{x\in \pl\mathscr{D}_t}|\theta(t,x)|+|\iota_0^{-1}(t)|\le C\lt(\varepsilon_1^{-1}, K_0,  \mathscr{E}_{E}(0)\rt).  \label{}
\end{align}\end{subequations}
\end{thm}

\begin{rmk}
It follows from \ef{CL-5.16} and \ef{23.thm-2-d} that
$\|v(t, \cdot)\|^2_{H^5(\mathscr{D}_t)}\le C \underline{\varrho}^{-1} \mathscr{E}_{E}(t).$
\end{rmk}

We make the following a priori assumptions: for $t\in [0,T]$,
\begin{align*}
&2^{-1} {\rm Vol} \mathscr{D}_0\le  {\rm Vol} \mathscr{D}_t \le 2  {\rm Vol} \mathscr{D}_0, & \\
&2^{-1}\underline{\varrho}  \le \rho( t,x)\le 2 \overline{\varrho}, \ \  |s( t,x)|\le \overline{s}   & {\rm in}  \ \ \mathscr{D}_t , \\
&|\pl (v, p, s)(t,x)|      \le M   & {\rm in} \ \ \mathscr{D}_t,
\\
&|\theta(t,x)|+  \iota_0^{-1}(t)\le K \ \ {\rm and} \ \   \iota_1^{-1}(t)\le K_1 &   {\rm on} \ \  \pl \mathscr{D}_t, \\
&-\pl_N p(t,x) \ge \ea_{b} , \     |\pl_N D_t p(t,x)| \le  L ,   \
|\pl_N D_t^2 p(t,x)| \le  \widetilde{L}  &   {\rm on} \ \  \pl \mathscr{D}_t,
\end{align*}
where $M$,  $K$, $K_1$,  $\ea_{b}$, $L$ and $\widetilde{L}$  are positive constants.
As in Section \ref{sec3}, we may assume without loss of generality that
$$ {\rm Vol} \mathscr{D}_0= {4\pi}/{3},    \
 \underline{\varrho}=2^{-1},  \  \overline{\varrho}=2,
   \
\overline{s}=1, $$
which implies that for $t\in [0,T]$,
 $$ {2\pi}/3  \le  {\rm Vol} \mathscr{D}_t \le {8\pi}/3 ,\
4^{-1}\le \rho(t,x)\le 4 \ {\rm and} \
|s( t,x )|\le 1  \  {\rm in}  \ \mathscr{D}_t. $$
Similarly, it holds that for $t\in [0,T]$,
$$
 |D_t \rho(t,x)| +|D_t p(t,x)|+|\pl \rho(t,x)|     \les    M \   \  {\rm in}   \ \ \mathscr{D}_t.
$$

\subsection{Regularity estimates}
\begin{prop}  Let $\kappa=0$ in \ef{1-b}, then
it holds that
\begin{align}
& \sum_{0\le r \le 5}\big(\|\pl^r v\|_{L^2}^2+
\|D_t \pl^r(\rho, p, s) \|_{L^2}^2 \big)+  \| \pl^5 {\rm div} v\|_{L^2}^2 +\|D_t v\|_{L^2}^2  \notag\\
& +\sum_{0\le r \le 4} \big(\|D_t \pl^r({\rm curl} v , \ {\rm div}  v )\|_{L^2}^2 +\|D_t \pl D_t^r p \|_{L^2}^2\big)
 +\| D_t^7 \rho - \Da D_t^5 p  \|_{L^2}^2  \notag\\
&
+\|\pl D_t^6 \rho - \rho_p D_t \pl D_t^5 p \|_{L^2}^2
 +\sum_{1\le r\le 5}(|\pl^r p|_{L^2}^2
+|\Pi\pl^r D_t p|_{L^2}^2)
\notag\\
&\le C(M, K,K_1,\ea_b^{-1}, L,  \widetilde{L}) \sum_{1\le i \le 4}  \mathscr{E}_{E}^i.
 \notag
\end{align}
\end{prop}

The proof consists of the following two lemmas, Lemmas \ref{23lem-1} and \ref{23lem-2}.

\begin{lem}\label{23lem-1} Let $\ka=0$ in \ef{1-b}, then it holds that
\begin{subequations}\label{23-01}
\begin{align}
& \sum_{i=1,2,3}  \|\pl^i ( v,\rho, p, s) \|_{L^\infty}^2+ \sum_{i+j=0,1,2} \|  \pl^i D_t \pl^j (v,\rho,p)  \|_{L^\infty}^2 +\|D_t \pl s \|_{L^\infty}^2  \notag \\
&\  \
+\sum_{i=2,3}
\|  D_t^{i} (\rho, p) \|_{L^\infty}^2
 +\sum_{i+j=2}
  \|D_t^i \pl D_t^j (v,  \rho, p) \|_{L^\iy}^2
 +|\pl D_t p |^2_{L^2}  \notag \\
 & \ \
  +\sum_{1\le i \le 5} \|(\pl^i v , D_t^i p) \|_{L^2}^2
 +\sum_{1\le i \le 4} |\pl^i ( v, p) |_{L^2}^2
\le C(M, K_1)\mathscr{E}_{E}, \label{23-01-a}\\
&\sum_{i+j=3,4} \|\pl^i D_t \pl^j (v,\rho, p)   \|_{L^2}^2
+ \sum_{i+j+k=2,3}\|\pl^i D_t \pl^j  D_t \pl^k  (v,\rho, p)  \|_{L^2}^2
 \notag\\
&  \ \
+\sum_{i+j+k=3}\|D_t^i\pl D_t^j \pl D_t^k (v,\rho, p)\|_{L^2}^2
+\sum_{i+j=3,4}\|D_t^i\pl D_t^j (v,\rho,p) \|_{L^2}^2
 \notag\\
&  \ \  + \sum_{1\le i\le 5}\lt\|D_t \pl^i  s\rt\|_{L^2}^2 +\sum_{i=1,2}(|\pl^{i+1} D_t p |^2_{L^2}  +  |\pl^i D_t^2 p |^2_{L^2})
\notag\\
&  \ \ + |\pl  D_t^3 p |^2_{L^2}
\le C(M,K,K_1)(\mathscr{E}_{E}^2+\mathscr{E}_{E}),
\label{23-01-b}\\
& \sum_{i+j=5}\|D_t^i\pl D_t^j (\rho,p) \|_{L^2}^2+ \lt\| (D_t^6 p, \ D_t \pl^4 {\rm curl} v ) \rt\|_{L^2 }^2 +\|   \Da D_t^4 p  \|_{L^2}^2\notag \\
&  \ \
+\|\pl D_t^6 \rho - \rho_p D_t \pl D_t^5 p \|_{L^2}^2
\le C(M,K,K_1)(\mathscr{E}_{E}^3+\mathscr{E}_{E}^2
+\mathscr{E}_{E}),\label{23-01-c}
\\
& |\overline\pl^{i} \theta|_{L^2}^2 \le C( M, K,K_1,\ea_b^{-1})\mathscr{E}_{E}, \ \  i=0,1,2, \label{23-01-d}\\
&  |\overline\pl^{3} \theta|_{L^2}^2 \le C( M, K,K_1,\ea_b^{-1})(\mathscr{E}_{E}^2+ \mathscr{E}_{E}).\label{23-01-e}
\end{align}\end{subequations}
\end{lem}

{\em Proof}. In a similar way to deriving Lemmas \ref{22lem-1}-\ref{22lem-3}, we can prove \ef{23-01} by noting that
\begin{align*}
&|\pl^4 D_t {\rm curl} v|\les \sum_{1\le i \le 5, \ 1\le j \le 6-i} M^{6-i-j}|\pl^i (\rho,p)||\pl^j \rho|
\\
& \quad +\sum_{1\le i\le 3}|\pl^i v||\pl^{6-i} v|
+|\pl^2 \rho|^2|\pl^2 p|
 ,
\\
&|D_t^6 p|\les \sum_{1\le i\le 6} M^{6-i}D_t^i \rho + \sum_{\substack{2\le i\le 4,\ 2\le j\le 3\\  i+j\le 6}}M^{6-i-j}|D_t^i \rho||D_t^j \rho|
+|D_t^2 \rho|^3,
\\
& |\pl D_t^6 \rho - \rho_p D_t \pl D_t^5 p|
\les
|D_t^3 p||D_t^2 p||\pl D_t p|+
|D_t^2 p|^2 (|\pl D_t^2 p|
\\ &\quad
+M |(D_t^2, \pl D_t) p|)
+\sum_{0\le i\le 1, \ 1\le j \le 6,  \ i+j\le 6} M^{7-i-j}|\pl^i D_t^j p|
\\ & \quad
+ \sum_{0\le i\le 1, \ 1\le j \le 5,  \ 2\le i+j, \  2\le l \le 6-j} M^{7-i-j-l}|\pl^i D_t^j p||D_t^l p|.
\end{align*}
This finishes the proof of the lemma. \hfill $\Box$

\begin{lem} \label{23lem-2}  Let $\ka=0$ in \ef{1-b}, then it holds that
\begin{subequations}\label{23-02}
\begin{align}
& \sum_{i+j+k=4}\|D_t^i\pl D_t^j \pl D_t^k (\rho,p) \|_{L^2}^2+ \sum_{i+j=5}\|D_t^i\pl D_t^j  v\|_{L^2}^2
\notag \\ & \quad
\le C(M,K,K_1)(\mathscr{E}_{E}^3+\mathscr{E}_{E}^2
+\mathscr{E}_{E}),\label{23-02-a}\\
& \sum_{i+j+k+l=3}\|\pl^i D_t \pl^j D_t \pl^k D_t \pl^l (\rho,p) \|_{L^2}^2+\|  \Da D_t^5 p-D_t^7 \rho  \|_{L^2}^2
\notag \\ & \quad
+\|D_t^4 \Da v \|_{L^2}^2 \le C(M,K,K_1,\ea_b^{-1}) (\mathscr{E}_{E}^3+\mathscr{E}_{E}^2
+\mathscr{E}_{E}),\label{23-02-b}\\
& \sum_{i+j+k=4}\|\pl^i D_t \pl^j D_t \pl^k (\rho,p) \|_{L^2}^2+\|D_t \pl^4   {\rm div} v\|_{L^2}^2+ |\pl^5 p|_{L^2}^2
\notag \\ & \quad
\le C(M,K,K_1,\ea_b^{-1},\widetilde{L}) (\mathscr{E}_{E}^3+\mathscr{E}_{E}^2
+\mathscr{E}_{E}),\label{23-02-c}\\
&\sum_{i+j=5}\|\pl^i D_t \pl^j (\rho,p) \|_{L^2}^2+ |\pl^4 D_t p|_{L^2}^2 + |\Pi \pl^5 D_t p|_{L^2}^2+
\|\pl^5 {\rm div} v\|_{L^2}^2
\notag\\
&\quad
\le C(M,K,K_1,\ea_b^{-1},L) (\mathscr{E}_{E}^4+\mathscr{E}_{E}^3+\mathscr{E}_{E}^2
+\mathscr{E}_{E}). \label{23-02-d}
\end{align}\end{subequations}
\end{lem}

{\em Proof}. It follows from \ef{23-01}, \ef{LZ-A5a},
\ef{3-7-3} and \ef{22.3.17} that for $\ka=0$,
\be\label{23-Jan21}
\sum_{i+j+k=4}\|D_t^i\pl D_t^j \pl D_t^k (\rho,p) \|_{L^2}^2
\le C(M,K,K_1)(\mathscr{E}_{E}^3+\mathscr{E}_{E}^2
+\mathscr{E}_{E}),
\ee
where we have used
\begin{align*}
&|\pl^2 D_t^4 \rho|
\les \sum_{\substack{0\le i \le 2, \ 0\le j\le 4, \
0\le k\le 2-i\\ 0\le l\le 4-j,\ 2\le i+j, k+l}}M^{6-i-j-k-l}|\pl^i D_t^j (p,s)| |\pl^k D_t^l (p,s)|
\\
& \ \ +\sum_{\substack{0\le i \le 2, \ 0\le j\le 4 \\  1\le i+j}}
M^{6-i-j} |\pl^i D_t^j (p,s)| +|D_t^2 p|^2 |\pl^2(p,s)|
 +|D_t^2 p||\pl D_t p|^2  .
\end{align*}
In view of \ef{7.5-2-a} and  \ef{7.7}, we see that for $\ka=0,1$,
\begin{align*}
& |D_t^5 \pl v + \rho^{-1} D_t^4 \pl^2 p-\ka D_t^4 \pl^2 \phi | \les  \sum_{\substack{0\le i  \le 4 }} |  D_t^i \pl v || D_t^{4-i} \pl v|+ |D_t^2 \rho|^2 |\pl^2 p|
\\
& \quad +\sum_{0\le i\le 3} M^{4-i} |D_t^i \pl^2 p| +\sum_{\substack{0\le i\le 2 ,\
   2\le j \le 4-i}} M^{4-i-j} |D_t^i \pl^2 p| |D_t^j \rho|+\mathcal{J}_2 ,
\end{align*}
where $\mathcal{J}_2$ is given by \ef{23-1-21}. This, together with \ef{23-01}, \ef{23-Jan21} and \ef{7.7}, implies that for $\ka=0$,
\be\label{23-Jan21-1}
 \sum_{i+j=5}\|D_t^i\pl D_t^j  v\|_{L^2}^2
\le C(M,K,K_1)(\mathscr{E}_{E}^3+\mathscr{E}_{E}^2
+\mathscr{E}_{E}).
\ee
Due to  \ef{22.4.2},
 \ef{3-7-3} and \ef{22.3.17}, one has
 that for $\ka=0,1$,
\begin{align}
& |\pl^2 \Da D_t^2 p-\ka \pl^2 D_t^2 (\rho \Da \phi)|\les
|\pl^2 D_t^4 \rho|+(|\pl^2 \rho||D_t^2 \rho| +|\pl D_t\rho||D_t\pl\rho|) |\pl^2 p|\notag\\
& \quad  +\sum_{1\le i\le 4} |\pl^i v||\pl^{5-i }D_t p |+\sum_{\substack{0\le i,j,k,l  \le 2\\
i+k, j+l\le 2  }} |\pl^i D_t^j \rho||\pl^k D_t^l \pl v ||\pl^{2-i-k} D_t^{2-j-l} \pl v|
\notag\\&\quad
+\sum_{0\le i \le 2, \ 0\le j, k\le 1}
|\pl^i D_t^j \pl^{k+1} v||\pl^{2-i} D_t^{1-j}\pl^{2-k} p|+|\pl^2\rho||(\pl D_t, D_t\pl)\rho|| D_t\pl p|
\notag\\
& \quad+\sum_{\substack{0\le i,j \le 2  }}M^{5-i-j}\big(|\pl^i D_t^j \pl p|
+\sum_{\max\{0,\ 2-i-j\}\le k\le 1 }M^{1-k}|\pl^i D_t^j \pl^k \rho|\big)\notag\\
&\quad +\sum_{\substack{0\le i,j \le 2, \
1\le i+j,  \
0\le k\le 2-i \\ 0\le l \le 2-j, \ \max\{0,\ 2-k-l\}\le m\le 1}}M^{5-i-j-k-l-m} |\pl^i D_t^j \pl p||\pl^k D_t^l \pl^m \rho |\notag\\
&\quad
+\sum_{\substack{0\le i,j\le 2, \ 0\le   k\le 1, \ 0\le l \le 2-i, \  0\le m\le 2-j\\ 0\le n\le 1-k,  \ 2\le i+j+k, l+m+n \  }}M^{6-i-j-k-l-m-n} |\pl^i D_t^j \pl^k \rho||\pl^l D_t^m \pl^n \rho|, \label{23.1}
\end{align}
which means, with the help of \ef{23-01} and \ef{23-Jan21}, that for $\ka=0$,
\be\label{23-1-22}
\|\pl^2 \Da D_t^2 p\|_{L^2}^2 \le C(M,K,K_1)(\mathscr{E}_{E}^3+\mathscr{E}_{E}^2
+\mathscr{E}_{E}).
\ee
Hence, we use the same way as that of proving \ef{7.10.3} to obtain that for $\ka=0$,
\begin{align}
\sum_{i+j+k=4}\|\pl^i D_t \pl^j D_t \pl^k (\rho,p) \|_{L^2}^2
\le C(M,K,K_1,\ea_b^{-1},\widetilde{L}) (\mathscr{E}_{E}^3+\mathscr{E}_{E}^2
+\mathscr{E}_{E}),
\label{23-1-22-1}
\end{align}
where we have used
\begin{align*}
&|\pl^4 D_t^2 \rho|
\les |D_t^2 p| |\pl^2(p,s)|^2
+|\pl D_t p|^2 |\pl^2(p,s)|\\
&\ \ +\sum_{\substack{0\le i \le 4 ,\  0\le j\le 2, \ 1\le i+j}}
M^{6-i-j}|\pl^i D_t^j (p,s)| \\
& \ \  + \sum_{\substack{0\le i \le 4 ,\  0\le j\le 2, \
0\le k\le 4-i \\  0\le l\le 2-j, \
 2\le i+j,  k+l}}M^{6-i-j-k-l}|\pl^i D_t^j (p,s)| |\pl^k D_t^l (p,s)| .
\end{align*}
It follows from \ef{1-a} and \ef{7.5-3} that for $\ka=0,1$,
\begin{align*}
&|\pl^4 D_t  {\rm div} v|=|\pl^4 D_t (\rho^{-1}D_t \rho)|\les |\pl^2 \rho|^2|D_t^2 \rho|
\\
&\quad +|\pl D_t \rho|^2|\pl^2 \rho| +
\sum_{\substack{0\le i\le 4 , \ 0\le j\le 2, \
1\le i+j}} M^{6-i-j}|\pl^i D_t^j \rho|
\\
&\quad
  +\sum_{\substack{0\le i\le 4 ,\ 0\le j\le 2, \ 0\le k\le 4-i \\  0\le l \le 2-j,\
2\le i+j,  k+l}} M^{6-i-j-k-l}|\pl^i D_t^j \rho||\pl^k D_t^l \rho|
,\\
& |\pl^4 \Da p-\ka \pl^4(\rho \Da\phi)|\les |\pl^4 D_t^2 \rho|+
|\pl^2 v|^2|\pl^2 \rho|+ |\pl^2 \rho|^2|\pl^2 p|\\
&\quad+ \sum_{1\le i\le 3}|\pl^i v |(|\pl^{6-i} v|+M|\pl^{5-i}(v,\rho)|)
+\sum_{1\le i\le 5} M^{6-i}|\pl^i(\rho, p)|
\\
& \quad +\sum_{2\le i\le 4, \ 2\le j\le 6-i} M^{6-i-j}|\pl^i(\rho, p)|
|\pl^j \rho|,
\end{align*}
which, together with \ef{CL-5.28'}, \ef{3-7-3}, \ef{23-01} and \ef{23-1-22-1}, implies that for $\ka=0$,
\begin{align}
\|D_t \pl^4   {\rm div} v\|_{L^2}^2+ |\pl^5 p|_{L^2}^2
\le C(M,K,K_1,\ea_b^{-1},\widetilde{L}) (\mathscr{E}_{E}^3+\mathscr{E}_{E}^2
+\mathscr{E}_{E}).
\label{23-1-22-2}
\end{align}

It follows from \ef{22.4.2},
 \ef{3-7-3} and \ef{22.3.17} that for $\ka=0,1$,
\begin{align*}
& |\pl  \Da D_t^3 p-\ka\pl D_t^3(\rho \Da\phi)|\les
|\pl  D_t^5 \rho|+\sum_{0\le i, j,  k\le 1}
|\pl^i D_t^j \pl^{k+1} v||\pl^{1-i} D_t^{1-j}\pl^{2-k} D_t p| \\
& \quad +\sum_{1\le i\le 3} |\pl^i v||\pl^{4-i }D_t^2 p |+\sum_{0\le i,   k\le 1,\ 0\le j\le 2}
|\pl^i D_t^j \pl^{k+1} v||\pl^{1-i} D_t^{2-j}\pl^{2-k} p|
\\&\quad
+\sum_{\substack{0\le i,k  \le 1, \ 0\le j, l \le 3\\
i+k\le 1, j+l\le 3  }} |\pl^i D_t^j \rho||\pl^k D_t^l \pl v ||\pl^{1-i-k} D_t^{3-j-l} \pl v|
\\
&\quad +\sum_{\substack{0\le i  \le 1, \ 0\le j\le 3  }}M^{5-i-j}\big(|\pl^i D_t^j \pl p|
+\sum_{\max\{0,\ 2-i-j\}\le k\le 1 }M^{1-k}|\pl^i D_t^j \pl^k \rho|\big)\\
&\quad +\sum_{\substack{0\le i  \le 1, \ 0\le j\le 3, \
1\le i+j,  \
0\le k\le 1-i \\ 0\le l \le 3-j, \ \max\{0,\ 2-k-l\}\le m\le 1}}M^{5-i-j-k-l-m} |\pl^i D_t^j \pl p||\pl^k D_t^l \pl^m \rho |\\
&\quad
+\sum_{\substack{0\le i,k\le 1, \ 0\le   j\le 3, \ 0\le l \le 1-i, \  0\le m\le 3-j\\ 0\le n\le 1-k,  \ 2\le i+j+k, l+m+n \  }}M^{6-i-j-k-l-m-n} |\pl^i D_t^j \pl^k \rho||\pl^l D_t^m \pl^n \rho|\\
&
\quad +|D_t\pl p|(|\pl D_t \rho||D_t \pl \rho|+|D_t^2 \rho ||\pl^2 \rho|)+ |\pl^2 p||D_t^2 \rho||D_t \pl \rho|,
\end{align*}
which means, with the aid of \ef{23-01}, that for $\ka=0$,
$$
\|\pl  \Da D_t^3 p\|_{L^2}^2 \le C(M,K,K_1) (\mathscr{E}_{E}^3+\mathscr{E}_{E}^2
+\mathscr{E}_{E}).
$$
Then, we may employ the same way as that of deriving \ef{lem5-a} to get that for $\ka=0$,
\begin{align}
\sum_{i+j+k+l=3}\|\pl^i D_t \pl^j D_t \pl^k D_t \pl^l (\rho,p) \|_{L^2}^2
\le C(M,K,K_1,\ea_b^{-1}) (\mathscr{E}_{E}^3+\mathscr{E}_{E}^2
+\mathscr{E}_{E}),
\label{23-1-23}
\end{align}
where we have used
\begin{align*}
&|\pl^3 D_t^3 \rho|
\les \sum_{\substack{0\le i, j\le 3 ,\
0\le k\le 3-i \\
0\le l\le 3-j, \ 2\le i+j, k+l}}M^{6-i-j-k-l}|\pl^i D_t^j (p,s)| |\pl^k D_t^l (p,s)|\\
& \ \
 +\sum_{\substack{0\le i, j\le 3 ,\  1\le i+j}}M^{6-i-j}
|\pl^i D_t^j (p,s)|+|D_t^2 p||\pl D_t p||\pl^2(p,s)|
+|\pl D_t p|^3.
\end{align*}
Due to \ef{23.1.15}, one has that for $\ka=0,1$,
\begin{align*}
& |D_t^4 \Da v + \rho^{-1} D_t^3 \pl \Da p-\ka D_t^3\pl \Da\phi|
   \les
\sum_{0\le i\le 3}|D_t^i\pl v||D_t^{3-i}\pl^2 v|
\\ &
 \ +\sum_{\substack{0\le i \le 3\\ 0\le j\le 2 , \   i+j\le 4 }} M^{5-i-j} |D_t^i\pl^{j+1}p|
+\sum_{\substack{0\le i\le 3\\ 0\le j\le 2, \ 2\le  i+j }} M^{6-i-j} |D_t^i\pl^{j}\rho|
\\
&\  +\sum_{\substack{0\le i \le 3, \ 0\le j\le 2  , \ 1\le i+j \le  4\\
0\le k\le 3-i, \ 0 \le l \le 2-j , \ 2 \le k+l}} M^{5-i-j-k-l} |D_t^i\pl^{j+1}p|  |D_t^k\pl^{l}\rho|
\\
&\
+|D_t\pl p|(|D_t^2 \rho||\pl^2 \rho|+|D_t\pl \rho|^2)
+|\pl^2 p||D_t^2 \rho||D_t \pl \rho|
\\
& \ +\sum_{\substack{0\le i \le 3, \ 0\le j\le 2  , \
0\le k\le 3-i\\ 0 \le l \le 2-j , \ 2 \le i+j,k+l}} M^{6-i-j-k-l} |D_t^i\pl^{j}\rho|  |D_t^k\pl^{l}\rho|,
\end{align*}
which, together with \ef{23-01}
and
\ef{23-1-23}, implies that for $\ka=0$,
\be\label{23-1-24}
\|D_t^4 \Da v \|_{L^2}^2 \le C(M,K,K_1,\ea_b^{-1}) (\mathscr{E}_{E}^3+\mathscr{E}_{E}^2
+\mathscr{E}_{E}).
\ee
It follows from \ef{22.4.2}   and \ef{3-7-3} that
for $\ka=0,1$,
\begin{align*}
& |\pl^3 \Da D_t p-\ka\pl^3   D_t (\rho \Da \phi) |\les
|\pl^3 D_t^3 \rho|+|\pl^2\rho||(\pl D_t, D_t\pl)\rho||\pl^2 p|
+|\pl^2 \rho|^2 |D_t \pl p| \\
& \quad +\sum_{1\le i\le 5} |\pl^i v||\pl^{6-i }p |+\sum_{\substack{0\le i, k \le 3, \ 0\le j,l\le 1\\
i+k\le 3, \ j+l\le 1  }} |\pl^i D_t^j \rho||\pl^k D_t^l \pl v ||\pl^{3-i-k} D_t^{1-j-l} \pl v|\\
& \quad +\sum_{\substack{0\le i \le 3 ,\ 0\le j\le 1}}M^{5-i-j}\big(|\pl^i D_t^j \pl p|
+\sum_{\max\{0,\ 2-i-j\}\le k\le 1 }M^{1-k}|\pl^i D_t^j \pl^k \rho|\big)\\
&\quad +\sum_{\substack{0\le i \le 3, \ 0\le j\le 1, \
1\le i+j,  \
0\le k\le 3-i \\ 0\le l \le 1-j, \ \max\{0,\ 2-k-l\}\le m\le 1}}M^{5-i-j-k-l-m} |\pl^i D_t^j \pl p||\pl^k D_t^l \pl^m \rho |\\
&\quad
+\sum_{\substack{0\le i\le 3, \ 0\le j,  k\le 1, \ 2\le i+j+k, \ 2\le l \le 4-i-k}}M^{6-i-j-k-l} |\pl^i D_t^j \pl^k \rho||\pl^l \rho|,
\end{align*}
which implies, with the aid of \ef{23-01} and
\ef{23-1-23},  that for $\ka=0$,
$$
 \| \pl^3 \Da D_t p \|_{L^2}^2
\le C(M,K,K_1,\ea_b^{-1}) (\mathscr{E}_{E}^3+\mathscr{E}_{E}^2
+\mathscr{E}_{E}).$$
Therefore, we use the same way as that of proving \ef{7.10.3} to obtain that for $\ka=0$,
\begin{align}
&\sum_{i+j=5}\|\pl^i D_t \pl^j (\rho,p) \|_{L^2}^2+ |\Pi \pl^4 D_t p|_{L^2}^2 + |\pl^4 D_t p|_{L^2}^2 + |\Pi \pl^5 D_t p|_{L^2}^2
\notag\\
&\quad
+
\|\pl^5 {\rm div} v\|_{L^2}^2
\le C(M,K,K_1,\ea_b^{-1},L) (\mathscr{E}_{E}^4+\mathscr{E}_{E}^3+\mathscr{E}_{E}^2
+\mathscr{E}_{E}),
\label{23-1-24}
\end{align}
where we have used
\begin{align*}
&|\pl^5 ({\rm div} v, D_t \rho)|
= |\pl^5( (-\rho^{-1}, 1) \rho_p D_t p))|
\\ &
\les
 \sum_{0\le i\le 5} M^{5-i} |\pl^i D_t p|
+\sum_{0\le i \le 3, \ 2\le j\le 5-i }M^{5-i-j}|\pl^j(p,s)||\pl^i D_t p|
\\ &\quad
+|\pl^2(p,s)|^2(|\pl D_t p| +M^2)
+M|\pl^3(p,s)||\pl^2(p,s)| .
\end{align*}

Finally, it follows from \ef{22.4.2},
 \ef{3-7-3} and \ef{22.3.17} that for $\ka=0,1$,
\begin{align*}
&|  \Da D_t^5 p-D_t^7 \rho -\ka D_t^5(\rho \Da \phi) |
 \les \sum_{\substack{0\le i,j  \le 5, \   i+j\le 5  }} |  D_t^i \rho||  D_t^j \pl v || D_t^{5-i-j} \pl v|\\
& \ +\sum_{1\le i\le 2} |\pl^i v||\pl^{3-i }D_t^4 p |
+\sum_{0\le i, j \le 1}
| D_t^i \pl^{j+1} v||  D_t^{1-i}\pl^{2-j} D_t^3 p|
\\ & \
+\sum_{0\le i \le 2, \ 1\le j\le 2}
| D_t^i \pl^{j } v||  D_t^{2-i}\pl^{3-j} D_t^2 p|+M|D_t^4\Da v|+|D_t^4 \pl v||\pl^2 p|
\\ & \
+\sum_{0\le i  \le 3, \ 1\le j\le 2}
| D_t^i \pl^{j } v| ( |  D_t^{3-i}\pl^{3-j} D_t  p|+|D_t^{4-i}\pl^{3-j}  p|)
\\&\
+\sum_{0\le i\le 5}M^{6-i}\big(|D_t^i\pl p|+\sum_{\max\{0,\ 2-i\}\le j\le 1} M^{1-j}|D_t^i \pl^j \rho|\big)
\\ & \
+\sum_{\substack{1\le i\le 4, \ 1\le j\le 5-i \\ \max\{0, \ 2-j\}\le k\le 1}} M^{6-i-j-k} |D_t^i \pl p||D_t^j \pl^k \rho|
+|D_t^2 \pl p||D_t^2 \rho||D_t \pl \rho|
\\ & \
+|D_t \pl p|(|D_t^3\rho||D_t \pl \rho|+(|D_t^2 \pl \rho|+M|(D_t^2, D_t \pl)\rho|)|D_t^2 \rho|)
\\ & \
+\sum_{\substack{1\le i\le 4, \ 0\le j\le 1, \ 1\le k\le 5-i \\ 0\le l\le 1-j,  \  2\le i+j, k+l }} M^{7-i-j-k-l} |D_t^i \pl^j \rho||D_t^k \pl^l \rho| + M|D_t^2 \rho|^2|D_t\pl \rho|,
\end{align*}
which proves, using
 \ef{23-01}, \ef{23-Jan21}, \ef{23-Jan21-1} and \ef{23-1-24},  that for $\ka=0$,
$$\|  \Da D_t^5 p-D_t^7 \rho  \|_{L^2}^2
\le C(M,K,K_1,\ea_b^{-1}) (\mathscr{E}_{E}^3+\mathscr{E}_{E}^2
+\mathscr{E}_{E}).$$
Hence, \ef{23-02} is a conclusion of  \ef{23-Jan21}, \ef{23-Jan21-1}, and \ef{23-1-22-1}-\ef{23-1-24}. \hfill$\Box$

\subsection{Energy estimates}

\begin{prop}\label{23.1.10} Let $\ka=0$ in \ef{1-b}, then it holds that
$$
\frac{d}{dt} \mathscr{E}_{E}\le C(M, K,K_1,\ea_b^{-1}, L,  \widetilde{L}) \sum_{1\le i \le 4}  \mathscr{E}_{E}^i  .
$$
\end{prop}

{\em Proof}. It suffices to prove  that
\begin{subequations}\label{22.8.6}\begin{align}
& \frac{d}{dt}E_r\le C(M,K,K_1,\ea_b^{-1},L,\widetilde{L}) \sum_{1\le i \le 4}   \mathscr{E}_{E}^i  , \ \ 1\le r\le 5,  \label{8.6}\\
&  \frac{d}{dt}(P_r + W_r) \le
C(M,K,K_1,\ea_b^{-1},L,\widetilde{L})\sum_{1\le i\le 4}   \mathscr{E}_{E}^i, \ \  0\le r \le 5, \label{8.6'}
\end{align}
\end{subequations}
where
\begin{align*}
& E_r=  \int_{\mathscr{D}_t}
\rho \delta^{mn} \zeta^{IJ} (\pl^r_I v_m)\pl^r_J v_ndx +  \int_{\pl \mathscr{D}_t} |\Pi \pl^r p|^2 (-\pl_N p)^{-1}  ds,\\
& P_r=\int_{\mathscr{D}_t}  \lt( |D_t^{r+1} \rho|^2+  \rho_p |\pl D_t^{r}p|^2\rt) dx,\ \
 W_0= \int_{\mathscr{D}_t} \lt(\rho |v|^2
+\rho^2+p^2 + s^2   \rt)dx ,\\
& W_r= \int_{\mathscr{D}_t}
\lt(|\pl^{r-1}({\rm curl} v, {\rm div} v)|^2+ |\pl^r (\rho,p,s)|^2 \rt)dx , \ \  r\ge 1.
\end{align*}
When $1\le r\le 4$, \ef{8.6} can be shown by choosing $\ka=0$ in \ef{3.2-4}. In  the case of $0\le r\le 4$, \ef{8.6'} can be proven by choosing $\ka=0$ in \ef{8.2}. So,
it is enough to prove \ef{22.8.6} for $r=5$.

Letting $\kappa=0$ in \ef{22.6.11-3}, we have
$$
\frac{d}{dt}E_5
=2 \int_{\mathscr{D}_t} \mathcal{H}_5 dx
-   \int_{\pl \mathscr{D}_t}
 (2\mathcal{L}_5
 +  |\Pi \pl^5 p|^2 (-\pl_N p)^{-1}    N^i \pl_N v_i)  ds,
$$
where $\mathcal{H}_5$ and  $\mathcal{L}_5$ are defined by
\ef{8.6.1} and \ef{8.6.2}, respectively.
In a similar way to deriving \ef{3.2-8} and \ef{8.1-4}, one has
\begin{align*}
&\|\mathcal{H}_5\|_{L^1}\le C(M,K) \lt( \|\pl^5 (v,  \rho, p,  {\rm div} v)  \|_{L^2}^2
 +  \|H_5 \|_{L^2}^2  \rt),\\
& |\mathcal{L}_5|_{L^1 }  \le C(  \ea_b^{-1},M, L)(|\pl^5 p|_{L^2}^2+ |\Pi \pl^5 D_t p|_{L^2}^2 \notag\\
&\quad +C(K)\sum_{2\le i \le 5} \|\pl^i v\|_{L^2}^2 \sum_{2\le j\le 5} \|\pl^j p\|_{L^2}^2   ),
\end{align*}
where
$$
 H_5 =
 |\pl^4 v||\pl^2 v|
 +|\pl^3 v|^2
 +|\pl^2 \rho|^2|\pl^2 p|
 +\sum_{\substack{1\le i, j \le 4,\ i+j\le 6}}|\pl^i \rho||\pl^j(\rho, p)|.$$
This, together with \ef{23-01} and \ef{23-02}, proves \ef{8.6} for $r=5$.

In a similar way to derive \ef{8.2-2} and \ef{8.2-4}, we  use \ef{23-01} and \ef{23-02} to get
\begin{align*}
&\frac{d}{dt} P_5
\les  \int_{\mathscr{D}_t }  (|D_t^7 \rho - \Da D_t^5 p|^2 +
|\pl D_t^6 \rho - \rho_p D_t \pl D_t^5 p |^2 )dx
\notag \\
&\quad + (M+1)P_5
 \le  C(M,K,K_1,\ea_b^{-1})\sum_{1\le i\le 3}   \mathscr{E}_{E}^i,
\\
&\frac{d}{dt} W_5
 \le   \int_{\mathscr{D}_t}
\lt(|D_t\pl^{4}({\rm curl} v, {\rm div} v)|^2   +   |D_t \pl^5 (\rho,p, s)|^2   \rt)dx
 \notag\\
 &+(M+1) W_5
 \le C(M,K,K_1,\ea_b^{-1},L,\widetilde{L}) \sum_{1\le i \le 4}  \mathscr{E}_{E}^i .
\end{align*}
So, \ef{8.6'} holds for $r=5$. We finish the proof of the lemma. \hfill $\Box$

\subsection{Proof of Theorem \ref{thm-2}}

The proof  follows from Proposition \ref{prop-2}, which is stated as follows.

\begin{prop}\label{prop-2}
Let $\ka=0$ in \ef{1-b}, then there exists a continuous function
$\overline{T}>0$ such that
\begin{align*}
& 2^{-1}{\rm Vol} \mathscr{D}_0\le  {\rm Vol} \mathscr{D}_t\le 2 {\rm Vol} \mathscr{D}_0 , \label{}\\
 & 2^{-1}\underline{\varrho} \le \rho( t,x)\le 2\overline{\varrho}, \ \  |s(t,x)|\le \overline{s}, \ \  x\in \mathscr{D}_t,
 \label{}\\
 &
-\pl_N p(t,x)\ge   2^{-1}\varepsilon_1,   \ \ x\in  \pl \mathscr{D}_t,  \\
& \iota_1^{-1}(t)\le   16K_0,\ \
 \mathscr{E}_{E}(t)\le 2 \mathscr{E}_{E}(0),\\
& \|   \pl  (v, p, s )(t,\cdot) \|_{L^\infty} \le 2 \|   \pl  ( v, p, s)(0,\cdot) \|_{L^\infty}, \\
&|\theta(t,\cdot)|_{L^\iy}+\iota_0^{-1}(t)\le  C\lt(\varepsilon_1^{-1}, K_0,  \mathscr{E}_{E}(0)\rt),\\
&
|\pl_N D_t p(t,\cdot)|_{L^\iy}+ |\pl_N D_t^2 p(t,\cdot)|_{L^\iy}     \le C(  \underline{\varrho}^{-1},\overline{\varrho},
\overline{s}, \varepsilon_1^{-1}, K_0, \mathscr{E}_{E}(0),
{\rm Vol} \mathscr{D}_0)
\end{align*}
for $t\le \overline{T}(  \underline{\varrho}^{-1},\overline{\varrho},
\overline{s}, \varepsilon_1^{-1}, K_0, \mathscr{E}_{E}(0),
{\rm Vol} \mathscr{D}_0)$.
\end{prop}

{\em Proof}. The most part of the proof is the same as that of Proposition \ref{prop23.2.20}.
As shown in Proposition \ref{prop23.2.20}, we can prove that there exists a continuous function
$\widetilde{T}>0$ such that
\begin{subequations}\label{4.3.2}\begin{align}
& 2^{-1}{\rm Vol} \mathscr{D}_0\le  {\rm Vol} \mathscr{D}_t\le 2 {\rm Vol} \mathscr{D}_0 , \label{4.3.2-b}\\
 & 2^{-1} \underline{\varrho}  \le \rho(t, x) \le 2\overline{\varrho}, \ \ |s(t,x)|\le  \overline{s},
 \ \  x\in \mathscr{D}_t,\label{4.3.2-a}\\
&
-\pl_N p(t,x)\ge   2^{-1}\varepsilon_1,  \ \ x\in  \pl \mathscr{D}_t, \label{4.3.2-e} \\
& \iota_1^{-1}(t)\le   16K_0,\ \
 \mathscr{E}_{E}(t)\le 2 \mathscr{E}_{E}(0),\label{4.3.2-f}\\
& \|   \pl  (v, p, s )(t,\cdot) \|_{L^\infty} \le 2 \|   \pl  ( v, p, s)(0,\cdot) \|_{L^\infty}, \label{4.3.2-g}
\end{align}\end{subequations}
for $t\le \widetilde{T}(  K,   L, \widetilde{L}, \underline{\varrho}^{-1},\overline{\varrho},
\overline{s}, \varepsilon_1^{-1}, K_0, \mathscr{E}_{E}(0),
{\rm Vol} \mathscr{D}_0) $.

It follows from \ef{22.5.7-1}, \ef{CL-A.8}, \ef{CL-5.19'}, \ef{4.3.2-e} and \ef{4.3.2-f}  that
for $t\le \widetilde{T}$,
\begin{align*}
&|\theta(t,\cdot)|_{L^\iy}=|(\pl_N p)^{-1}(t,\cdot)|_{L^\iy}|\Pi \pl^2 p(t,\cdot)|_{L^\iy}\notag\\
& \le |(\pl_N p)^{-1}(t,\cdot)|_{L^\iy}| \pl^2 p(t,\cdot)|_{L^\iy}
\le C_1\lt(\varepsilon_1^{-1}, K_0,  \mathscr{E}_{E}(0)\rt),
\label{}
\end{align*}
which implies, with the help of \ef{lemkk1} and  \ef{4.3.2-f}, that for $t\le \widetilde{T}$,
$$
\iota_0^{-1}(t) \le \max\lt\{32K_0,\ \
C_1\lt(\varepsilon_1^{-1}, K_0,  \mathscr{E}_{E}(0)\rt) \rt\}.
$$
So, we may choose
\be\label{4.3.5}
K=64K_0+4C_1\lt(\varepsilon_1^{-1}, K_0,  \mathscr{E}_{E}(0)\rt).
\ee
In a similar way to deriving \ef{7.1-6}, we use
  \ef{CL-5.34}, \ef{23-01-a},
\ef{23-01-b}, \ef{4.3.2} and \ef{4.3.5} to get  that for $t\le \widetilde{T}$,
$$
|\pl_N D_t p(t,\cdot)|_{L^\iy}
\le C_2( \underline{\varrho}^{-1},\overline{\varrho},
\overline{s}, \varepsilon_1^{-1}, K_0, \mathscr{E}_{E}(0),
{\rm Vol} \mathscr{D}_0 ).
$$
Similarly, we use \ef{CL-5.34}, \ef{23-01-b}, \ef{23-1-22}, \ef{4.3.2} and \ef{4.3.5} to obtain that
for $t\le \widetilde{T}$,
$$
|\pl_N D_t^2 p(t,\cdot)|_{L^\iy}
\le C_3( \underline{\varrho}^{-1},\overline{\varrho},
\overline{s}, \varepsilon_1^{-1}, K_0, \mathscr{E}_{E}(0),
{\rm Vol} \mathscr{D}_0 ).
$$
This finishes the proof of the proposition
by choosing $L=2 C_2$ and $\widetilde{L}=2 C_3$. \hfill $\Box$

\section{The isentropic Euler-Poisson equations}\label{sec5}
In this section, we investigate the free boundary problem \ef{22.1} of the isentropic Euler-Poisson equations under the stability condition \ef{23.8}. The main results are given in Theorem \ref{thm-3}.
Let $h=h(\rho)=\int_1^\rho \lambda^{-1} p'(\lambda)d\lambda$, then the higher-order energy functionals are defined by
\begin{align*}
& \mathscr{E} (t)= \int_{\mathscr{D}_t}   |v|^2 dx
+\sum_{1\le r \le 5} \int_{\pl \mathscr{D}_t} |\Pi \pl^r   (\phi-h )|^2 (\pl_N (\phi-h ))^{-1}  ds
\notag\\
&
\quad +\sum_{1\le r \le 5} \int_{\mathscr{D}_t}
\lt(  \delta^{mn} \zeta^{IJ} (\pl^r_I v_m)\pl^r_J v_n+ |\pl^{r-1}{\rm curl} v|^2+|\pl^{r-1} {\rm div} v|^2\rt)dx\notag\\
&
\quad + \sum_{0\le r \le 5} \int_{\mathscr{D}_t}  (|\pl^r \rho|^2 + |\pl^r p|^2 + |\pl^r \phi|^2  +|D_t^{r+1} \rho|^2+  \rho_p |\pl D_t^{r}p|^2  )dx. \label{}
\end{align*}
In order to state the main result of the present work, we set
\begin{align*}
&\underline{\varrho}=\min_{x\in \mathscr{D}_0} \rho_0(x), \  \  \overline{\varrho}=\max_{x\in \mathscr{D}_0}\rho_0(x),
 \ \  \varepsilon_1= \min_{x\in \pl\mathscr{D}_0} \pl_N (
\phi-h)( 0,x) ,  \label{}\\
&K_0= \max_{x\in \pl\mathscr{D}_0} |\theta(0,x)|
+|{\iota_0}^{-1}(0)|,
\label{}
\end{align*}
where $h(0,x)=h(\rho_0(x))$, and $\phi(0,x)$ is determined by the Dirichlet problem \ef{1-d}. With these notations, the main results of this section are stated as follows:

\begin{thm}\label{thm-3}
Suppose that \eqref{22.conP} hold, and
$0<  \underline{\varrho}, \overline{\varrho},
\varepsilon_1,    K_0, \mathscr{E}(0), {\rm Vol} \mathscr{D}_0<\iy$.
Then there exists a continuous function
$\mathscr{T}\lt(
\underline{\varrho}^{-1},\overline{\varrho},
  \varepsilon_1^{-1},  K_0, \mathscr{E}(0),
{\rm Vol} \mathscr{D}_0\rt)>0$
such that  any smooth solution of the free boundary problem \eqref{22.1} and \eqref{1-d} for $0\le t\le T$ with $T\le \mathscr{T}$
satisfies the  following estimates:  for $0\le t\le T$,
\begin{subequations}\begin{align}
&\mathscr{E}(t)\le 2 \mathscr{E}(0),  \ \
 2^{-1}{\rm Vol} \mathscr{D}_0\le  {\rm Vol} \mathscr{D}_t\le 2 {\rm Vol} \mathscr{D}_0 , \label{}\\
&2^{-1} \underline{\varrho} \le \min_{x\in \mathscr{D}_t} \rho(t, x), \  \  \max_{x\in \mathscr{D}_t}\rho(t, x)\le 2\overline{\varrho}, \label{23.2.21-1}\\
&2^{-1}\varepsilon_1\le \min_{x\in \pl\mathscr{D}_t}\pl_N (\phi - h)( t,x) ,
\label{23.5} \\
&\max_{x\in \pl\mathscr{D}_t}|\theta(t,x)|+|\iota_0^{-1}(t)|\le C\lt(\underline{\varrho}^{-1},\overline{\varrho},
  \varepsilon_1^{-1}, K_0, \mathscr{E}(0)\rt),
\end{align}\end{subequations}
where $h=h(\rho)=\int_1^\rho \lambda^{-1} p'(\lambda)d\lambda$.
\end{thm}

\begin{rmk}
It follows from \ef{CL-5.16}  that
$\|v(t, \cdot)\|^2_{H^5(\mathscr{D}_t)}\le C \mathscr{E}(t).$
\end{rmk}

We make the following a priori assumptions: for $t\in [0,T]$,
\begin{align*}
&2^{-1} {\rm Vol} \mathscr{D}_0\le  {\rm Vol} \mathscr{D}_t \le 2  {\rm Vol} \mathscr{D}_0, & \\
&2^{-1}\underline{\varrho} \le \rho( t,x)\le 2\overline{\varrho},   \ \   |\pl (v, p ,\phi)(t,x)|      \le M       & {\rm in} \ \ \mathscr{D}_t,
\\
&|\theta(t,x)|+  \iota_0^{-1}(t)\le K \ \ {\rm and} \ \   \iota_1^{-1}(t)\le K_1 &   {\rm on} \ \  \pl \mathscr{D}_t,\\
&\pl_N (\phi-h)(t,x) \ge \ea_{b},  \ \  |\pl_N D_t p(t,x)| \le  L  &   {\rm on} \ \  \pl \mathscr{D}_t ,\\
&  |\pl_N D_t  \phi(t,x)| \le  \bar L ,  \ \ |\pl_N D_t^2 p(t,x)| \le  \widetilde{L}       &   {\rm on} \ \  \pl \mathscr{D}_t ,
\end{align*}
where $M$,  $K$, $K_1$,  $\ea_{b}$,   $L$, $\bar L$ and $\widetilde{L}$ are positive constants. As in Section \ref{sec3},
it follows from the maximal principle that for $t\in [0,T]$ and $x\in \mathscr{D}_t$,
$$\min\lt\{0, \ -\ln \max_{x\in \mathscr{D}_t} \rho(t, x)\rt\}\le \phi(t,x )\le \max\lt\{0, \ -\ln \min_{x\in \mathscr{D}_t} \rho(t, x)\rt\} .$$
We may assume without loss of generality that
$$ {\rm Vol} \mathscr{D}_0= {4\pi}/{3},   \  \
\underline{\varrho} =2^{-1},  \  \ \overline{\varrho}=2,
$$
which implies that for $t\in [0,T]$,
$$ {2\pi}/3  \le  {\rm Vol} \mathscr{D}_t \le {8\pi}/3 ,\ \  4^{-1}\le \rho(t,x)\le 4 \ {\rm and} \
  |\phi(t,x )|\le \ln 4 \   {\rm in} \  \mathscr{D}_t. $$
Similarly, one has, due to $\pl h = \rho^{-1}\pl p$, that   for
$t\in [0,T]$,
$$|D_t \rho(t,x)| +|D_t p(t,x)|+|\pl \rho(t,x)|
+|\pl h(t,x)|    \les    M \ \ {\rm in} \ \ \mathscr{D}_t.$$

\subsection{Regularity estimates}
\begin{prop}
Let $h=h(\rho)=\int_1^\rho \lambda^{-1} p'(\lambda)d\lambda$, then it holds that
\begin{align}
& \sum_{0\le r \le 5}\big(\|\pl^r (v,h)\|_{L^2}^2+
\|D_t \pl^r(\rho, p, \phi) \|_{L^2}^2 \big)+  \| \pl^5 {\rm div} v\|_{L^2}^2 +\|D_t v\|_{L^2}^2  \notag\\
& +\sum_{0\le r \le 4} \big(\|D_t \pl^r({\rm curl} v , \ {\rm div}  v )\|_{L^2}^2 +\|D_t \pl D_t^r p \|_{L^2}^2\big)
 +\| D_t^7 \rho - \Da D_t^5 p  \|_{L^2}^2  \notag\\
&
+\|\pl D_t^6 \rho - \rho_p D_t \pl D_t^5 p \|_{L^2}^2
 +\sum_{1\le r\le 5}(|\pl^r(h, \phi)|_{L^2}^2
+|\Pi\pl^r D_t(h, \phi)|_{L^2}^2)
\notag\\
&\le C(M, K,K_1,\ea_b^{-1}, L, \bar L,  \widetilde{L}) \sum_{1\le i \le 7}  \mathscr{E}^i.\notag
\end{align}
\end{prop}

The proof consists of the following two lemmas, Lemmas \ref{23lem5.1} and \ref{23lem5.2}.

\begin{lem}\label{23lem5.1} Let $h=h(\rho)=\int_1^\rho \lambda^{-1} p'(\lambda)d\lambda$, then it holds that
\begin{subequations}\label{lem5.1}
\begin{align}
& \sum_{i=1,2,3}  \|\pl^i ( v,\rho, p, \phi) \|_{L^\infty}^2
 + \sum_{i+j=0,1,2} \|  \pl^i D_t \pl^j (v,\rho,p)  \|_{L^\infty}^2
+\|   \pl h \|_{L^\iy}^2  \notag \\
 & \ \
 +\sum_{1\le i\le 3} \|\pl^i h \|^2_{L^2}
 +\sum_{1\le i \le 5}  \|(\pl^i v, \ D_t^i p ) \|_{L^2}^2
+\sum_{i=1,2} (\|\pl^i D_t h\|_{L^2}^2   \notag\\
  & \ \  + |\pl^i h |_{L^2}^2 )     +\sum_{1\le i \le 4} |\pl^i ( v, p,\phi) |_{L^2}^2
 +|\pl D_t ( p ,h)  |^2_{L^2}
\le C(M,  K_1)\mathscr{E} , \label{lem5.1-d'}\\
&\sum_{i=4,5}( \| \pl^{i-2} h \|_{L^\iy}^2 + \|\pl^{i} h \|^2_{L^2}+ |\pl^{i-1} h |_{L^2}^2  )
 \le C(M,K_1)(\mathscr{E}^2+\mathscr{E}),
\label{lem5.1-e'}\\
&
\|  D_t  \phi \|_{L^\iy}^2 +\sum_{i+k=2,\ j=0,1}
  \|D_t^i \pl^j D_t^k ( \rho, p) \|_{L^\iy}^2+\|  D_t^3(\rho,p)  \|_{L^\iy}^2\notag \\
 & \ \
 +
\sum_{i+j=1,2}\|\pl^i D_t \pl^j  \phi\|_{L^2}^2
     +\sum_{i+j=2}\|D_t^i \pl D_t^j v\|_{L^2}^2
 +|\pl D_t   \phi |^2_{L^2} \notag \\
 & \ \
\le C(M, K, K_1)\mathscr{E} , \label{lem5.1-d}\\
&\|D_t^2\phi \|^2_{L^\iy} +\sum_{i+j=3,4} \|\pl^i D_t \pl^j (v,\rho, p)   \|_{L^2}^2
+\sum_{i=3,4} \|\pl^i D_t h\|_{L^2}^2
 \notag\\
&  \ \ +\sum_{i+j+k=1,2}\|\pl^i D_t \pl^j  D_t \pl^k  \phi  \|_{L^2}^2
+ \sum_{i+j+k=2,3}\|\pl^i D_t \pl^j  D_t \pl^k  (\rho, p)  \|_{L^2}^2
\notag\\
&  \ \ +\sum_{i+j=3}\|D_t^i\pl D_t^j (v,\rho,p) \|_{L^2}^2  +\sum_{i+j=4}\|D_t^i\pl D_t^j (\rho,p) \|_{L^2}^2  +
 |\pl  D_t^2 \phi |^2_{L^2} \notag\\
&  \ \  +\sum_{i=1,2}(
  |\pl^{i+1} D_t (p,h) |^2_{L^2}  +  |\pl^i D_t^2 p |^2_{L^2})
\le C(M,K,K_1)(\mathscr{E}^2+\mathscr{E}),
\label{lem5.1-e}\\
& \lt\| (D_t^6 p, \ D_t \pl^4 {\rm curl} v ) \rt\|_{L^2 }^2
+\|\pl D_t^6 \rho - \rho_p D_t \pl D_t^5 p \|_{L^2}^2
\notag \\
&  \ \ \le C(M,K,K_1)(\mathscr{E}^3+\mathscr{E}^2
+\mathscr{E}) , \label{lem5.1-g} \\
& |\overline\pl^{i} \theta|_{L^2}^2 \le C( M, K,K_1,\ea_b^{-1})\mathscr{E}, \ \  i=0,1, \label{lem5.1-a}\\
& |\overline\pl^{2} \theta|_{L^2}^2 \le C( M, K,K_1,\ea_b^{-1})(\mathscr{E}^2+ \mathscr{E}), \label{lem5.1-b}\\
&  |\overline\pl^{3} \theta|_{L^2}^2 +|\Pi \pl^5 ( p, \phi)|_{L^2}^2 + |\pl^5 \phi|_{L^2}^2 \le C( M, K,K_1,\ea_b^{-1})\sum_{1\le i\le 4}\mathscr{E}^i,\label{lem5.1-c}\\
&\sum_{i+j=1,2} \|  \pl^i D_t \pl^j \phi  \|_{L^\iy}^2 +\sum_{i+j=2}\|D_t^i \pl D_t^j v\|_{L^\iy}^2 +\sum_{i+j=3,4} \| \pl^i D_t  \pl^j \phi \|_{L^2}^2
\notag\\
&  \ \ + \sum_{i+j+k=2,3}\|\pl^i D_t \pl^j  D_t \pl^k  v  \|_{L^2}^2
 +\sum_{i+j+k=3}\|D_t^i\pl D_t^j \pl D_t^k (\rho, p)\|_{L^2}^2
 \notag\\
 &  \ \  +  \sum_{i=2,3,4}|\pl^i D_t  \phi|_{L^2}^2 + |\pl^2  D_t^2 \phi |^2_{L^2}
  + |\pl  D_t^3 p |^2_{L^2}
\notag\\
 &  \ \  \le C( M,   K, K_1, \ea_b^{-1}, \bar L)(\mathscr{E}^2+\mathscr{E} ), \label{lem5.1-f}
\\
& \sum_{i+j= 2}\| D_t^i \pl   D_t^j \phi  \|_{L^\iy}^2+\|D_t^3 \phi\|_{L^\iy}^2 +  \sum_{i+j+k=3}\big(\|\pl^i D_t \pl^j  D_t \pl^k  \phi  \|_{L^2}^2 \notag\\
 &  \ \ + \|D_t^i\pl D_t^j \pl D_t^k (v,\phi)\|_{L^2}^2\big)
+\sum_{i+j=3}\|D_t^i\pl D_t^j \phi\|_{L^2}^2
 +\sum_{i+j=4}\|D_t^i\pl D_t^j v \|_{L^2}^2 \notag\\
 &  \ \ +\sum_{i+j=5}\|D_t^i\pl D_t^j (\rho,p) \|_{L^2}^2
  +  |\Pi\pl^3 D_t^2 \phi |^2_{L^2}\notag\\
 &  \ \  \le C( M,   K, K_1, \ea_b^{-1}, \bar L)(\mathscr{E}^3+\mathscr{E}^2+ \mathscr{E}), \label{lem5.1-h} \\
& \|  D_t^4\phi \|_{L^\iy}^2 +\sum_{i+j=5} \| \pl^i D_t  \pl^j \phi \|_{L^2}^2
+\sum_{i+j=4}\|D_t^i\pl D_t^j \phi\|_{L^2}^2
 \notag\\
 &  \ \
+\sum_{i+j+k=4}\|D_t^i \pl D_t^j \pl D_t^k \phi\|_{L^2}^2 + \|\Da D_t^4 p\|_{L^2}^2
  +  |\Pi\pl^5 D_t  \phi|_{L^2}^2
  \notag\\
 &  \ \  \le
  C( M,   K, K_1, \ea_b^{-1}, \bar L)(\mathscr{E}^4+\mathscr{E}^3+\mathscr{E}^2+ \mathscr{E}), \label{lem5.1-i}\\
  & \| D_t^5 \phi\|_{L^2}+ \|\pl D_t^5 \phi \|_{L^2}^2 \le C( M,   K, K_1, \ea_b^{-1}, \bar L)(\mathscr{E}^5+\mathscr{E}^4+\mathscr{E}^3+\mathscr{E}^2+ \mathscr{E}).\label{lem5.1-j}
\end{align}\end{subequations}
\end{lem}

{\em Proof}. Since we can prove this lemma by the same idea of Lemmas \ref{22lem-1}-\ref{22lem-4} and \ref{23lem-1},  we only present the additional calculations needed here, especially on the higher-order derivatives of $h$ and $\phi$.
The bound for $\|\pl^i h \|_{L^2}$ $(1\le i \le 5)$ can be derived easily from  $\pl h= \rho^{-1} \pl p$ and the bound obtained for $(\rho,p)$, based on which we can get the desired bound for $\|\pl^i h\|_{L^\iy}$ $(1\le i\le 3)$ and $|\pl^i h |_{L^2}$ $(1\le i \le 4)$. Thus, the bound for $|\pl^i \theta|_{L^2}$ $(0\le i\le 3)$ follows from \ef{hb2}.
So, we can use \ef{hb1} to get the desired bound for $|\Pi\pl^5 p|_{L^2}$ and  $|\Pi\pl^5 \phi|_{L^2}$. The bound for $|\pl^5 \phi|_{L^2}$ comes   from
$$
|\pl^4 \Da \phi| \le |\pl^4 e^{-\phi}|
+|\pl^4 \phi|
\les |\pl^4 \rho|+\sum_{1\le i\le 4}M^{4-i}|\pl^i \phi|+|\pl^2 \phi|^2.
$$
The bound for $\|\pl^i D_t h \|_{L^2}$ $(1\le i \le 4)$ follows from $D_t h=\rho^{-1}D_t p$ and the bound obtained for $(\rho,p)$, and
that for $|\pl^i D_t h |_{L^2}$ $(1\le i \le 3)$ from \ef{CL-5.19'}.

The bound for $|\pl^4 D_t \phi|_{L^2}$, $|\Pi\pl^5 D_t \phi|_{L^2}$ and $\|\pl^5 D_t \phi\|_{L^2}$
follows  from
$$\|\pl^3 \Da D_t\phi \|_{L^2}^2 \le C( M,   K, K_1, \ea_b^{-1}, \bar L)(\mathscr{E}^2 +\mathscr{E}),$$
which is due to
\begin{align*}
&|\pl^3 \mathfrak{G}_1|\les
\sum_{1\le i\le 5}|\pl^i v||\pl^{6-i}\phi|  \ \ {\rm and}  \ \
 |\pl^3 ( e^{-\phi} D_t   \phi)|\les
\sum_{0\le i\le 3}M^{3-i}|\pl^i D_t \phi|
\\
& \qquad +(|\pl D_t\phi| + M |D_t    \phi|)
|\pl^2 \phi|+|D_t    \phi|
|\pl^3 \phi|.
\end{align*}
In a similar way to deriving \ef{lem5-a}, we can obtain the desired
bound for $|\pl^2 D_t^2 \phi|_{L^2}$, $|\Pi\pl^3 D_t^2 \phi|_{L^2}$ and $\|\pl^3 D_t^2 \phi\|_{L^2}$, by noticing that
$$
\|\pl \Da D_t^2 \phi \|_{L^2}^2 \le C( M,   K, K_1, \ea_b^{-1}, \bar L)(\mathscr{E}^2 +\mathscr{E}),
$$
which follows from
\begin{align*}
&|\pl  ( e^{-\phi} D_t^2   \phi)|\les |\pl D_t^2 \phi|+M|D_t^2 \phi| \ \ {\rm and} \ \  |\pl \mathfrak{G}_2|\les \ \sum_{1\le i\le 3}|\pl^i v||\pl^{4-i}D_t\phi|\\
&\quad
+\sum_{0\le i,j,k\le 1}|\pl^i D_t^j \pl^{k+1}v||\pl^{1-i}D_t^{1-j}\pl^{2-k}\phi|
+ |\pl D_t\phi| |D_t\phi| +M|D_t \phi|^2
 .
\end{align*}
Finally, \ef{lem5.1-j} can be shown with the aid of
\begin{align*}
&|\mathfrak{G}_5|\les \sum_{0\le r\le 3, \ 0\le i\le r,\ 0\le j\le 1}|D_t^i \pl^{j+1} v ||D_t^{r-i}\pl^{2-j}D_t^{4-r} \phi|
\\
&\ \ +\sum_{0\le i\le 4}(|D_t^i \Da v||D_t^{4-i}\pl \phi| + |D_t^i \pl v||D_t^{4-i }\pl^2 \phi|)
+|D_t^4\phi||D_t\phi|
\\
& \ \ +|D_t^3 \phi|(|D_t^2 \phi|
+|D_t \phi|^2)+|D_t^2 \phi|^2|D_t \phi|
+|D_t^2 \phi||D_t \phi|^3
+|D_t \phi|^5.
\end{align*}
This finishes the proof of the lemma.
\hfill$\Box$

\begin{lem}\label{23lem5.2} It holds that
\begin{subequations}\label{lem5.2}
\begin{align}
& \sum_{i+j+k+l=3}\|\pl^i D_t \pl^j D_t \pl^k D_t \pl^l (\rho,p) \|_{L^2}^2
+\|D_t^4 \Da v \|_{L^2}^2
\notag \\ & \quad
 \le
C( M,   K, K_1, \ea_b^{-1}, \bar L)(\mathscr{E}^3+\mathscr{E}^2+ \mathscr{E}) ,\label{lem5.2-a}\\
& \sum_{i+j+k=4}\|D_t^i\pl D_t^j \pl D_t^k (\rho,p) \|_{L^2}^2+ \sum_{i+j=5}\|D_t^i\pl D_t^j  v\|_{L^2}^2
\notag \\ & \quad
 \le
C( M,   K, K_1, \ea_b^{-1}, \bar L)(\mathscr{E}^4+\mathscr{E}^3+\mathscr{E}^2+ \mathscr{E}) ,\label{lem5.2-b}\\
& \|  \Da D_t^5 p-D_t^7 \rho  \|_{L^2}^2
 \le
C( M,   K, K_1, \ea_b^{-1}, \bar L)(\mathscr{E}^5+\mathscr{E}^4+\mathscr{E}^3+\mathscr{E}^2+ \mathscr{E}) ,
\label{lem5.2-e}\\
& \sum_{i+j+k=4}\|\pl^i D_t \pl^j D_t \pl^k (\rho,p) \|_{L^2}^2+\|D_t \pl^4   {\rm div} v\|_{L^2}^2+ |\pl^5 (p, h)|_{L^2}^2
\notag \\ & \quad
\le  C( M,   K, K_1, \ea_b^{-1}, \bar L,\widetilde{L})
(\mathscr{E}^4+\mathscr{E}^3+\mathscr{E}^2+ \mathscr{E}),\label{lem5.2-c}\\
&\sum_{i+j=5}\|\pl^i D_t \pl^j (\rho,p) \|_{L^2}^2+
\|\pl^5 {\rm div} v\|_{L^2}^2 + |\pl^4 D_t p|_{L^2}^2 + |\Pi \pl^5 D_t p|_{L^2}^2
\notag\\
&\quad
\le C(M,K,K_1,\ea_b^{-1}, \bar L,L) \sum_{1\le i\le 5}\mathscr{E}^i, \label{lem5.2-d}\\
&  \|\pl^5 D_t h \|_{L^2}^2
+|\pl^4 D_t h |_{L^2}^2
\le C(M,K,K_1,\ea_b^{-1}, \bar L,L) \sum_{1\le i\le 5}\mathscr{E}^i, \label{lem5.2-d'}\\
&  |\Pi\pl^5 D_t h |_{L^2}^2
\le C(M,K,K_1,\ea_b^{-1}, \bar L,L) \sum_{1\le i\le 7}\mathscr{E}^i. \label{lem5.2-f}
\end{align}
\end{subequations}
\end{lem}

{\em Proof}. In a similar way to proving Lemma \ref{23lem-2}, we can obtain \ef{lem5.2} by virtue of the following additional calculations.  To get
the bound for $\|\pl^4 D_t^2 p\|_{L^2}$ and $|\pl^5 p|_{L^2}$ in \ef{lem5.2-c}, we need the following estimates, respectively.
\begin{align}
& |\pl^2 D_t^2 (\rho \Da\phi )|\les   |\pl^2 D_t^2 (\rho e^{-\phi})|+|\pl^2 D_t^2 \rho^2| \notag\\
& \quad \les  \sum_{0\le i, j \le 2}  |\pl^i D_t^j (\rho,\phi)||\pl^{2-i} D_t^{2-j} (\rho,\phi)|+
M^2|(D_t^2,\pl D_t)(\rho, \phi)|\notag \\
& \quad
+(|\pl^2(\rho, \phi)|+M^2)|D_t \phi|^2
+M ( | (\pl D_t, \pl^2 )(\rho, \phi)|   +M^2)|D_t\phi|, \label{23.2} \\
& |\pl^4 (\rho \Da\phi )|    \les \sum_{1\le i\le 4}M^{4-i}|\pl^i(\rho, \phi)|
+|\pl^2(\rho, \phi)|^2  .\notag
\end{align}
The following estimates are needed to derive
\ef{lem5.2-a}.
\begin{align*}
& |\pl  D_t^3 (\rho \Da\phi )|\les   |\pl  D_t^3 (\rho e^{-\phi})|+|\pl  D_t^3 \rho^2| \\
& \ \ \les  \sum_{0\le i\le 1, \ 0\le j\le 3}
|\pl^i D_t^j (\rho, \phi)||\pl^{1-i}D_t^{3-j}(\rho, \phi)|
+M^2|D_t^2\phi|\\
& \ \  +M|(D_t^2, \pl D_t)(\rho, \phi)||D_t \phi|
+(|\pl D_t(\rho,\phi)|+M^2)|D_t\phi|^2+M|D_t\phi|^3
,\\
& |D_t^3 \pl \Da \phi|
\les |D_t^3 \pl e^{-\phi}|+|D_t^3 \pl \rho|
\les |D_t^3 \pl (\rho,\phi)|
+M|D_t^2\phi||D_t\phi|
\\
&  \ \   +\sum_{0\le i\le 3, \ 0\le j\le 1,\ 1\le i+j\le 3} |D_t^i \pl^j \phi||D_t^{3-i}\pl^{1-j} \phi|+|D_t\pl \phi|||D_t \phi|^2
+M|D_t\phi|^3.
\end{align*}
We use  the following estimates to prove \ef{lem5.2-d} and \ef{lem5.2-e}, respectively.
\begin{align*}
& |\pl^3 D_t (\rho \Da\phi )|    \les
|\pl^3  D_t (\rho e^{-\phi})|+|\pl^3  D_t \rho^2|
\les M^2|\pl D_t(\rho, \phi)|
+M^3|D_t (\rho, \phi)|
\\
& \ \ + \sum_{0\le i\le 3, \ 0\le j\le 1}
|\pl^i D_t^j (\rho, \phi)|
|\pl^{3-i}D_t^{1-j}(\rho, \phi)|
+M|\pl^2(\rho, \phi)||D_t(\rho,\phi)|
,\\
& |D_t^5 (\rho \Da\phi )|    \les
| D_t^5 (\rho e^{-\phi})|+|D_t^5 \rho^2|
\les \sum_{0\le i\le 5}
|  D_t^i (\rho, \phi)|
| D_t^{5-i}(\rho, \phi)|
\\
&
\ \ +
|  D_t^3 (\rho, \phi)||  D_t  (\rho, \phi)|^2
+|  D_t^2 (\rho, \phi)|(|D_t^2 \phi||  D_t (\rho, \phi)|+|  D_t  (\rho, \phi)|^3)\\
& \ \ +|D_t\phi|^4|  D_t  (\rho, \phi)|.
\end{align*}

For $|\pl^5 h|_{L^2}$, it follows from
$\Da (h-\phi)=-{\rm div} D_t v=-D_t {\rm div} v - (\pl_i v^k)\pl_k v^i$ that
$$|\pl^4 \Da (h-\phi)| \les |D_t \pl^4 {\rm div} v| + \sum_{1\le i\le 3} |\pl^i v||\pl^{6-i} v|,
$$
which implies, due to  $|\Pi \pl^5 (h-\phi) |_{L^2}^2\les M \mathscr{E}$, that
$$|\pl^5(h-\phi)|_{L^2}^2\le C( M,   K, K_1, \ea_b^{-1}, \bar L,\widetilde{L})
(\mathscr{E}^4+\mathscr{E}^3+\mathscr{E}^2+ \mathscr{E}).$$
This, together with the bound obtained for $|\pl^5 \phi|_{L^2}$ in \ef{lem5.1-c}, gives the desired bound for $|\pl^5 h|_{L^2}$. The bound for $\|\pl^5 D_t h \|_{L^2}^2$ (and hence $|\pl^4 D_t h |_{L^2}$) follows from
\begin{align*}
& |\pl^5 D_t h|\les \sum_{1\le i\le 5} M^{5-i} (|\pl^i D_t p|+M|\pl^i \rho|)+|\pl^2 \rho|^2 (|\pl D_t p|
\\
& \quad +M^2) +|\pl^3 \rho||\pl^2 \rho|
 +\sum_{1\le i \le 3, \ 2\le j\le 5-i} M^{5-i-j}|\pl^i D_t p| |\pl^j \rho|.
\end{align*}
Note that on $ \pl \mathscr{D}_t$,
\be\label{23.2.16}
|\pl_N D_t h| =| (\pl_N \rho^{-1}) D_t p + \rho^{-1} |\pl_N D_t p|=|   \rho^{-1} |\pl_N D_t p|\le 4L
\ee
then we can use \ef{hb1}, and the bound obtained for $|\pl^i D_t h|_{L^2}$ $(1\le i\le 4)$ and $|\overline\pl^{i} \theta|_{L^2}$ $(1\le i\le 3)$
to prove \ef{lem5.2-f}. This finishes the proof of the lemma.
\hfill$\Box$

\subsection{Energy estimates}

\begin{prop}It holds that
\be\label{23.2.13}
\frac{d}{dt} \mathscr{E} \le C(M, K,K_1,\ea_b^{-1}, L, \bar L,  \widetilde{L}) \sum_{1\le i \le 7}  \mathscr{E}^i.
\ee
\end{prop}

{\em Proof}. It follows from \ef{8.8} that for $r\ge 1$,
$$
D_t \pl^r v + \pl^{r+1}(h(\rho)-\phi)= [D_t, \pl^r] v,
$$
which implies that
\begin{align}
&2^{-1} D_t(\da^{mn} \zeta^{IJ} (\pl^r_I v_m) \pl^r_J v_n)\notag \\
=&\da^{mn} \zeta^{IJ} (  D_t \pl^r_I v_m) \pl^r_J v_n +2^{-1}  \da^{mn} (D_t \zeta^{IJ}) (\pl^r_I v_m) \pl^r_J v_n\notag \\
=&\mathcal{A}_r- {\rm div} \lt(\zeta^{IJ}(  \pl^r_I (h(\rho)-\phi))\pl^r_J v\rt)  ,
\label{8.8-3}
\end{align}
where
\begin{align*}
& \mathcal{A}_r= \zeta^{IJ}(\pl^r_I (h(\rho)-\phi) ) \pl^r_J {\rm div} v
 +\da^{mn}\{2^{-1} (D_t \zeta^{IJ}) \pl^r_I v_m \notag\\
 & \ \  + (\pl_m \zeta^{IJ}) \pl^r_I  (h(\rho)-\phi)   + \zeta^{IJ}  [D_t, \pl^r_I] v_m      \}\pl^r_J v_n.
\end{align*}
Due to
$\pl_m (h(\rho)-\phi)  = N_m \pl_N (h(\rho)-\phi) $ on $\pl \mathscr{D}_t$, we have on $\pl \mathscr{D}_t$,
\begin{align}
&2^{-1}D_t\lt((\pl_N (h(\rho)-\phi)) ^{-1}  \zeta^{IJ}(\pl^r_I (h(\rho)-\phi) )\pl^r_J (h(\rho)-\phi) \rt)\notag\\
=& (\pl_N (h(\rho)-\phi) )^{-1} \zeta^{IJ} (\pl^r_I (h(\rho)-\phi) )D_t \pl^r_J (h(\rho)-\phi)
\notag
\\
& + 2^{-1}\lt(D_t((\pl_N (h(\rho)-\phi))^{-1} \zeta^{IJ})\rt) (\pl^r_I (h(\rho)-\phi) )\pl^r_J (h(\rho)-\phi)
\notag\\
=& \mathcal{B}_r-(\pl_N (h(\rho)-\phi))^{-1}
\zeta^{IJ} (\pl^r_I (h(\rho)-\phi) )
(\pl^r_J v^m) \pl_m (h(\rho)-\phi)
\notag \\
= & \mathcal{B}_r-   N_m  \zeta^{IJ} (\pl^r_I (h(\rho)-\phi) )
\pl^r_J v^m,\label{8.8-4}
\end{align}
 where
\begin{align*}
&\mathcal{B}_r=(\pl_N (h(\rho)-\phi))^{-1} \zeta^{IJ}(\pl^r_I (h(\rho)-\phi) ) \{D_t \pl^r_J (h(\rho)-\phi)
\notag\\
& \quad
- \pl^r_J D_t (h(\rho)-\phi) + (\pl^r_J v^m) \pl_m (h(\rho)-\phi)  \}
\notag\\
& \quad +(\pl_N (h(\rho)-\phi))^{-1}  \zeta^{IJ}(\pl^r_I (h(\rho)-\phi) ) \pl^r_J D_t (h(\rho)-\phi) \notag\\
& \quad +
2^{-1}\lt(D_t((\pl_N (h(\rho)-\phi))^{-1}\zeta^{IJ})\rt) (\pl^r_I (h(\rho)-\phi) )\pl^r_J (h(\rho)-\phi).
\end{align*}
In view of $\int_{\mathscr{D}_t}  {\rm div} \lt(\zeta^{IJ}(  \pl^r_I (h(\rho)-\phi))\pl^r_J v\rt) dx
= \int_{\pl \mathscr{D}_t}   N_m  \zeta^{IJ} (\pl^r_I (h(\rho)-\phi) )
\pl^r_J v^m
ds$, \ef{3.2-1},   \ef{8.8-3} and \ef{8.8-4},
we see that for $r\ge 1$,
\begin{align}
&
\frac{d}{dt}
E_r=\int_{\mathscr{D}_t}
\lt\{ D_t\lt( \delta^{mn} \zeta^{IJ} (\pl^r_I v_m)\pl^r_J v_n \rt)+
  \delta^{mn} \zeta^{IJ} (\pl^r_I v_m)(\pl^r_J v_n) {\rm div} v \rt\}dx
\notag \\
& +  \int_{\pl \mathscr{D}_t} \{ D_t\lt( |\Pi \pl^r   (\phi-h(\rho))|^2 (\pl_N (\phi-h(\rho)))^{-1} \rt)
\notag\\
& -|\Pi \pl^r   (\phi-h(\rho))|^2 (\pl_N (\phi-h(\rho)))^{-1} N^i \pl_N v_i  \}  ds
  =\mathcal{X}_r, \label{23.2.14}
\end{align}
where
\begin{align*}
& E_r=\int_{\mathscr{D}_t}
  \delta^{mn} \zeta^{IJ} (\pl^r_I v_m)\pl^r_J v_ndx
 \\
  & \quad +  \int_{\pl \mathscr{D}_t} |\Pi \pl^r   (\phi-h(\rho))|^2 (\pl_N (\phi-h(\rho)))^{-1}  ds
,\\
& \mathcal{X}_r=
\int_{\mathscr{D}_t}
\lt\{  2 \mathcal{A}_r+
  \delta^{mn} \zeta^{IJ} (\pl^r_I v_m)(\pl^r_J v_n) {\rm div} v \rt\}dx
  - \int_{\pl \mathscr{D}_t} \{ 2 \mathcal{B}_r \\
&\quad  +|\Pi \pl^r   (\phi-h(\rho))|^2 (\pl_N (\phi-h(\rho)))^{-1} N^i \pl_N v_i  \}  ds
.\end{align*}

In a similar way to proving \ef{3.2-8} and \ef{8.1-4}, we use \ef{lem5.1}-\ef{23.2.16} to obtain that for $1\le r\le 5$,
\begin{align}\label{23-1-17}
\|\mathcal{A}_r\|_{L^1}
+|\mathcal{B}_r|_{L^1}
 \le C(M, K,K_1,\ea_b^{-1}, L, \bar L,  \widetilde{L}) \sum_{1\le i \le 7}  \mathscr{E}^i,
\end{align}
where we have used
\begin{align*}
& |\mathcal{A}_r|
\les |\pl^r(h-\phi)||\pl^r {\rm div} v|
+\big(M|\pl^r v|+K|\pl^r(h-\phi)|
\\
& \quad +\sum_{1\le i\le r}|\pl^i v||\pl^{r+1-i}v|  \big)|\pl^r v|,\\
& |\mathcal{B}_r|\le C(  \ea_b^{-1},M, L, \bar L)|\pl^r (h-\phi)|  (  |\Pi \pl^r D_t (h-\phi)| \notag\\
&\quad +\sum_{1\le i\le r-1} |\pl^i v||\pl^{r+1-i} (h-\phi)| +|\pl^r (h-\phi)| \big).
\end{align*}
Substitute \ef{23-1-17} into \ef{23.2.14} to give
that for $1\le r\le 5$,
$$
\frac{d}{dt}
E_r\le C(M, K,K_1,\ea_b^{-1}, L, \bar L,  \widetilde{L}) \sum_{1\le i \le 7}  \mathscr{E}^i. $$
This proves \ef{23.2.13}, since the rest terms contained in $\mathscr{E}$ can be dealt with by the same idea of Lemma \ref{lem3.10} and Proposition \ref{23.1.10}.
\hfill $\Box$

\subsection{Proof of Theorem \ref{thm-3}}

The proof  follows from Proposition \ref{prop-3}, which is stated as follows.

\begin{prop}\label{prop-3}
There exists a continuous function
$\overline{T}>0$ such that
\begin{align*}
& 2^{-1}{\rm Vol} \mathscr{D}_0\le  {\rm Vol} \mathscr{D}_t\le 2 {\rm Vol} \mathscr{D}_0 , \label{}\\
 & 2^{-1}\underline{\varrho} \le \rho( t,x)\le 2\overline{\varrho},  \ \   x\in \mathscr{D}_t,
 \label{}\\
 &
-\pl_N p(t,x)\ge   2^{-1}\varepsilon_1,   \ \ x\in  \pl \mathscr{D}_t,  \\
& \iota_1^{-1}(t)\le   16K_0,\ \
 \mathscr{E}(t)\le 2 \mathscr{E}(0),\\
& \|   \pl  (v, p, s )(t,\cdot) \|_{L^\infty} \le 2 \|   \pl  ( v, p, s)(0,\cdot) \|_{L^\infty}, \\
&|\theta(t,\cdot)|_{L^\iy}+\iota_0^{-1}(t)
 \le C\lt(\underline{\varrho}^{-1},\overline{\varrho},
  \varepsilon_1^{-1}, K_0, \mathscr{E}(0)\rt),\\
&
 |\pl_N D_t p(t,\cdot)|_{L^\iy}+ |\pl_N D_t^2 p(t,\cdot)|_{L^\iy}   \le C(  \underline{\varrho}^{-1},\overline{\varrho},
\overline{s}, \varepsilon_1^{-1}, K_0, \mathscr{E}(0),
{\rm Vol} \mathscr{D}_0)
\end{align*}
for $t\le \overline{T}(  \underline{\varrho}^{-1},\overline{\varrho},
  \varepsilon_1^{-1}, K_0, \mathscr{E}(0),
{\rm Vol} \mathscr{D}_0)$.
\end{prop}

{\em Proof}. The most part of the proof is the same as that of Proposition \ref{prop23.2.20}.
As shown in Proposition \ref{prop23.2.20}, we can prove that there exists a continuous function
$\widetilde{T}>0$ such that
\begin{subequations}\label{5.3.2}\begin{align}
& 2^{-1}{\rm Vol} \mathscr{D}_0\le  {\rm Vol} \mathscr{D}_t\le 2 {\rm Vol} \mathscr{D}_0 , \label{}\\
 & 2^{-1} \underline{\varrho}  \le \rho(t, x) \le 2\overline{\varrho},
 \ \  x\in \mathscr{D}_t,\label{}\\
&
\pl_N (\phi -h) (t,x)\ge   2^{-1}\varepsilon_1,  \ \ x\in  \pl \mathscr{D}_t, \label{} \\
& \iota_1^{-1}(t)\le   16K_0,\ \
 \mathscr{E}(t)\le 2 \mathscr{E}(0),\label{5.3.2-d}\\
& \|   \pl  (v, p, \phi )(t,\cdot) \|_{L^\infty} \le 2 \|   \pl  ( v, p, \phi )(0,\cdot) \|_{L^\infty}, \label{}
\end{align}\end{subequations}
for $t\le \widetilde{T}(  K,   L, \bar L, \widetilde{L}, \underline{\varrho}^{-1},\overline{\varrho},
  \varepsilon_1^{-1}, K_0, \mathscr{E}(0),
{\rm Vol} \mathscr{D}_0) $, where we have used
\ef{23.2.16}.

It follows from \ef{22.5.7-1}, \ef{CL-A.8}, \ef{lem5.1-d'}, \ef{lem5.1-e'} and \ef{5.3.2}  that
for $t\le \widetilde{T}$,
\begin{align*}
|\theta(t,\cdot)|_{L^\iy}
 & \le |(\pl_N (\phi -h))^{-1}(t,\cdot)|_{L^\iy}|  \pl^2 (\phi -h)(t,\cdot)|_{L^\iy}\\
& \le C_1\lt(\underline{\varrho}^{-1},\overline{\varrho},
  \varepsilon_1^{-1}, K_0, \mathscr{E}(0)\rt),
\label{}
\end{align*}
which implies, with the help of \ef{lemkk1} and  \ef{5.3.2-d}, that for $t\le \widetilde{T}$,
$$
\iota_0^{-1}(t) \le \max\lt\{32K_0,\ \
C_1\lt(\underline{\varrho}^{-1},\overline{\varrho},
  \varepsilon_1^{-1}, K_0, \mathscr{E}(0)\rt) \rt\}.
$$
So, we may choose
\be\label{5.3.3}
K=64K_0+4C_1\lt(\underline{\varrho}^{-1},\overline{\varrho},
  \varepsilon_1^{-1}, K_0, \mathscr{E}(0)\rt).
\ee
In a similar way to deriving \ef{7.1-6}, we use
  \ef{CL-5.34}, \ef{lem5.1-d},
\ef{22-3-9-3}, \ef{5.3.2} and \ef{5.3.3} to obtain  that for $t\le \widetilde{T}$,
$$
|\pl_N D_t \phi(t,\cdot)|_{L^\iy}
\le C_2( \underline{\varrho}^{-1},\overline{\varrho},
  \varepsilon_1^{-1}, K_0, \mathscr{E}(0),
{\rm Vol} \mathscr{D}_0).
$$
Similarly, we use \ef{CL-5.34},
   \ef{lem5.1-d'},
\ef{lem5.1-e}, \ef{5.3.2} and \ef{5.3.3} to get  that for $t\le \widetilde{T}$,
$$
|\pl_N D_t p(t,\cdot)|_{L^\iy}
\le C_3( \underline{\varrho}^{-1},\overline{\varrho},
  \varepsilon_1^{-1}, K_0, \mathscr{E}(0),
{\rm Vol} \mathscr{D}_0).
$$
Again, we use  \ef{CL-5.34}, \ef{lem5.1-e}, \ef{23.1}, \ef{23.2}, \ef{4.3.2} and \ef{4.3.5} to obtain that
for $t\le \widetilde{T}$,
$$
|\pl_N D_t^2 p(t,\cdot)|_{L^\iy}
\le C_4( \underline{\varrho}^{-1},\overline{\varrho},
  \varepsilon_1^{-1}, K_0, \mathscr{E}(0),
{\rm Vol} \mathscr{D}_0).
$$
This finishes the proof of the proposition
by choosing $\bar L=2 C_2$, $ L=2 C_3$ and
$\widetilde{L}=2 C_4$. \hfill $\Box$

\section*{Aknowledgement}  Luo's research is supported by a grant from the Research Grants Council of the Hong Kong Special Administrative Region, China (Project No. 11306621). Trivisa
gratefully acknowledges the support of the National Science Foundation under the awards DMS-1614964
and DMS-2008568. Zeng's research is supported in part by NSFC  Grants  12171267 and 11822107.
 This research was initiated during Luo's visit to the University of Maryland, College Park  and Tsinghua University.

\bibliographystyle{plain}

\noindent{Tao Luo}\\
Department of Mathematics,\\
	City University of Hong Kong,
	Hong Kong\\
E-mail: taoluo@cityu.edu.hk\\

\noindent{Konstantina Trivisa}\\
Department of Mathematics,\\
	University of Maryland, College Park,
	USA \\
E-mail: trivisa@math.umd.edu\\

\noindent{Huihui Zeng}\\
Department of Mathematical Sciences, \\
Tsinghua University, \\
Beijing, 100084, China\\
Email : hhzeng@mail.tsinghua.edu.cn

\end{document}